\documentclass[a4paper,11pt]{scrartcl}

\usepackage{graphicx}
\usepackage[utf8]{inputenc}
\usepackage[english]{babel}
\usepackage{amsmath,amssymb,url,tabularx,stmaryrd,mathrsfs,nicefrac}
\usepackage{dsfont}
\usepackage[inline]{enumitem}

\usepackage{todonotes,comment,afterpage}
\usepackage{subcaption}
\usepackage{framed}
\usepackage{mathtools}
\usepackage{fullpage}
\usepackage[title]{appendix}
\usepackage{amsthm}
\usepackage{graphicx}
\usepackage{hyperref}
\usepackage{soul}

\theoremstyle{plain}
\newtheorem{theorem}{Theorem}[section]
\newtheorem{lemma}[theorem]{Lemma}
\newtheorem{corollary}[theorem]{Corollary}
\newtheorem{proposition}[theorem]{Proposition}

\newtheorem{assumption}[theorem]{Assumption}

\newtheorem{remark}[theorem]{Remark}

\numberwithin{equation}{section}
\numberwithin{figure}{section}
\numberwithin{table}{section}

\newcommand{\R}{\mathbb{R}}

\newcommand{\N}{\mathbb{N}}

\newcommand{\bigO}{\mathcal{O}}
\newcommand{\Omegaext}{\Omega^{\operatorname{ext}}}

\newcommand{\abs}[1]{\left|#1\right|}

\newcommand{\lefttriplenorm}{\ensuremath{\left| \! \left| \! \left|}}
\newcommand{\righttriplenorm}{\ensuremath{\right| \! \right| \! \right|}}
\newcommand{\energynorm}[1]{\lefttriplenorm #1 \righttriplenorm_{\dGnorm}}
\newcommand{\energynormp}[1]{\lefttriplenorm #1 \righttriplenorm_{\dGnormp}}
\newcommand{\bigtriplenorm}[1]{\ensuremath{#1| \! #1| \! #1|}}

\newcommand{\smalltriplenorm}[2]{\ensuremath{#1| \! #1| \! #1| #2 #1| \! #1 | \!#1|}}

\newcommand{\enorm}[1]{\smalltriplenorm{}{#1}_k}
\newcommand{\AAA}{\boldsymbol {\mathcal A}}

\newcommand{\x}{\mathbf{u}}
\newcommand{\xh}{\mathbf{u}_h}
\newcommand{\yh}{\mathbf{y}_h}

\newcommand{\support}{\operatorname{supp}}

\newcommand{\laplace}{\Delta}

\newcommand{\eremk}{\hbox{}\hfill\rule{0.8ex}{0.8ex}}

\def\TT{\mathcal{T}}
\newcommand{\TTC}[1]{{\mathcal{T}_{#1}}}
\newcommand{\TTDG}[1]{{\mathcal{T}_{#1}^{DG}}}
\def\SS{\mathcal{S}}
\def\PP{\mathscr{P}}

\def\ii{i}

\def\fdiv{\operatorname{div}}
\def\id{\mathrm{I}}

\def\LL{\mathcal{L}}

\def\bs{\mathbf}

\renewcommand{\Re}{\operatorname{Re}}
\renewcommand{\Im}{\operatorname{Im}}
\renewcommand{\i}{\ii}

\usepackage{tikz}
\usetikzlibrary{shapes}
\usepackage{pgfplots}
\usetikzlibrary{calc}

\usetikzlibrary{positioning}
\usetikzlibrary{arrows}
\usetikzlibrary{calc}
\usetikzlibrary{intersections, pgfplots.fillbetween,patterns}
\usepgfplotslibrary{colorbrewer}

\pgfplotsset{compat=1.13}
\iffalse 

\usetikzlibrary{external}
\tikzexternaldisable

\else

\fi

\newcommand{\Ad}{\nu}
\newcommand{\kn}{k n}
\newcommand{\nGamma}{\mathbf n_{\Gamma}} 

\newcommand{\gammaz}{\gamma_0} 
\newcommand{\gammazint}{\gammaz^{\mathrm{int}}} 
\newcommand{\gammazext}{\gammaz^{\mathrm{ext}}} 
\newcommand{\gammao}{\partial_{\nGamma}} 
\newcommand{\gammaoint}{\partial_{\nGamma}^{\mathrm{int}}} 
\newcommand{\gammaoext}{\partial_{\nGamma}^{\mathrm{ext}}} 
\newcommand{\gammanu}[2][]{\partial_{\nGamma}^{\nu,\mathrm{#1}} #2 }

\newcommand{\djmp}[1]{\llbracket #1 \rrbracket_\Gamma}
\newcommand{\njmp}[1]{\llbracket \partial^{\nu}_{\nGamma} #1 \rrbracket_\Gamma}

\newcommand{\aclass}{\boldsymbol {\mathcal A}}

\newcommand{\vext}{v^{\mathrm{ext}}} 
\newcommand{\lambdah}{\lambda_h} 
\newcommand{\uext}{u^\mathrm{ext}} 
\newcommand{\uhext}{\uext_h} 
\newcommand{\vhext}{\vext_h} 
\newcommand{\uh}{u_h} 
\newcommand{\vh}{v_h} 
\newcommand{\mh}{m_h} 

\newcommand{\psim}{\psi_m} 
\newcommand{\psitilde}{\psi_{\mathrm{ext}}} 
\newcommand{\psiext}{\psi_{\mathrm{ext}}} 
\newcommand{\phim}{\phi_m} 
\newcommand{\phiext}{\phi_{\mathrm{ext}}} 

\newcommand{\psiF}{\psi^{\mathcal{F}}} 
\newcommand{\psiA}{\psi^{\mathcal{A}}} 
\newcommand{\psimF}{\psim^{\mathcal{F}}} 
\newcommand{\psimA}{\psim^{\mathcal{A}}} 
\newcommand{\psiextF}{\psiext^{\mathcal{F}}} 
\newcommand{\psiextA}{\psiext^{\mathcal{A}}} 

\newcommand{\V}{\mathcal V} 
\newcommand{\Vz}{\V_0} 
\newcommand{\Vtilde}{\widetilde \V} 
\newcommand{\Vk}{\V_k} 

\newcommand{\Vtildek}{\Vtilde_k} 
\newcommand{\Vtildez}{\Vtilde_0} 
\newcommand{\K}{\mathcal K} 
\newcommand{\Kz}{\K_0} 
\newcommand{\Kprime}{\K^{\!\top}} 
\newcommand{\Kprimez}{\Kprime_0} 
\newcommand{\Ktilde}{\widetilde \K} 
\newcommand{\Ktildez}{\Ktilde_0} 
\newcommand{\W}{\mathcal W} 
\newcommand{\Wz}{\W_0} 
\newcommand{\Kk}{\K_k} 
\newcommand{\Kprimek}{\Kprime_k} 
\newcommand{\Ktildek}{\Ktilde_k} 
\newcommand{\Wk}{\W_k} 
\newcommand{\A}{\mathcal A} 
\newcommand{\Aprime}{\mathcal A'} 
\newcommand{\B}{\mathcal B} 
\newcommand{\Aprimek}{\Aprime_k} 
\newcommand{\Bk}{\B_k} 

\newcommand{\Vh}{V_h}
\newcommand{\Wh}{W_h}
\newcommand{\Zh}{Z_h}
\newcommand{\VhDG}{V_h^{DG}}
\newcommand{\Vhc}{V_h^{C}}

\newcommand{\Rm}{R_m}
\newcommand{\Rext}{R_{ext}}

\newcommand{\hff}[2][\eta]{H_{#2,#1}}
\newcommand{\lff}[2][\eta]{L_{#2,#1}}

\newcommand{\hffg}[2][\eta]{\hff[#1]{\Gamma}^{#2}}
\newcommand{\lffg}[2][\eta]{\lff[#1]{\Gamma}^{#2}}

\newcommand{\F}{\mathcal{F}}
\newcommand{\SV}{\SS_{\V}}
\newcommand{\SK}{\SS_{\K}}
\newcommand{\SKprime}{\SS_{\Kprime}}
\newcommand{\SW}{\SS_{\W}}

\newcommand{\AVtilde}{\widetilde{\A}_{\V}}
\newcommand{\AKtilde}{\widetilde{\A}_{\K}}
\newcommand{\T}{\TT}

\newcommand{\dGnorm}{dG(\Omega)}
\newcommand{\dGnormp}{dG^+(\Omega)}

\newcommand{\ahOmega}{a_h^\Omega}
\newcommand{\bhGamma}{b_h^\Gamma}
\newcommand{\taun}{\Omega_h}
\newcommand{\Fh}{\mathcal F_h}
\newcommand{\FhI}{\Fh^I}
\newcommand{\FhB}{\Fh^B}

\newcommand{\msf}{{\mathfrak h}}
\newcommand{\ldc}{\{\!\!\{}
\newcommand{\rdc}{\}\!\!\}}

\newcommand{\afrak}{\mathfrak{a}}
\newcommand{\bfrak}{\mathfrak{b}}
\newcommand{\dfrak}{\mathfrak{d}}

\newcommand{\bx}{\mathbf x}
\newcommand{\by}{\mathbf y}
\newcommand{\bn}{\mathbf n}

\newcommand{\sTheta}{s^{\prime\prime}}

\newcommand{\rhsm}{g}
\newcommand{\rhsext}{h}

\title{FEM-BEM coupling for the high-frequency Helmholtz problem\footnote{This research was funded by the Austrian Science Fund (FWF) projects
10.55776/F65 (JMM, IP), 10.55776/P33477 (IP), and 10.55776/P36150 (AR).}}

\author{Jens Markus Melenk\thanks{Institut f\"ur Analysis und Scientific Computing, TU Wien, 1040 Vienna, Austria (melenk@tuwien.ac.at)}, Ilaria Perugia\thanks{Fakult\"at f\"ur Mathematik, Universit\"at Wien, 1090 Vienna, Austria (ilaria.perugia@univie.ac.at)} and Alexander Rieder\thanks{Institut f\"ur Analysis und Scientific Computing, TU Wien, 1040 Vienna, Austria (alexander.rieder@tuwien.ac.at)}}
\date{\today}

\begin{document}
\maketitle
\begin{abstract}
We present a wavenumber-explicit analysis of FEM-BEM coupling methods for time-harmonic Helmholtz problems proposed in~\cite{MMPR20} for conforming discretizations and in~\cite{dgfembem} for discontinuous Galerkin (DG) volume discretizations. We show that the conditions that $kh/p$ be sufficiently small and that $\log(k) / p$ be bounded imply quasi-optimality of both conforming and DG-method, where $k$ is the wavenumber, $h$ the mesh size, and $p$ the approximation order. The analysis relies on a $k$-explicit regularity theory for a three-field coupling formulation.
\end{abstract}


\section{Introduction}
\label{sec:intro}
Many acoustic or electromagentic wave propagation problems in inhomogeneous media are posed in time-harmonic form and are additionally
naturally posed in full space ${\mathbb R}^d$. Numerically, the treatment of such problems is challenging for two reasons. First, time-harmonic 
problems at large wavenumbers $k$ are notoriously beset by dispersion errors (also known as pollution errors) in which 
the gap between the best approximation error and the actual numerical error widens as the wavenumber $k$ increases. 
For both Helmholtz and Maxwell problems, it has now become clear that high order methods are much better suited to 
control dispersion errors than low order methods; see, e.g., \cite{MS11,MPS13,bernkopf,lafontaine-spence-wunsch22,galkowski-lafontaine-spence-wunsch23,chaumont-frelet-gallistl-nicaise-tomezyk22,chaumont-frelet-nicaise20,nicaise-tomezyk20,galkowski2023sharp} and references therein where a mathematical 
analysis is put forward for this observation, and it is shown that the conditions 
\begin{equation}
\label{eq:intro-scale-resolution}
\frac{kh}{p} \mbox{ sufficient small} \qquad \mbox{ and } \qquad p = O(\log k)
\end{equation}
are sufficient to suppress dispersion errors, where $h$ is the mesh size and $p$ the approximation order. 
Second, the treatment of unbounded domains requires one to resort to a form of transparent boundary conditions 
such as PML \cite{berenger94,michielssen97,bramble-pasciak08,chaumont-frelet-gallistl-nicaise-tomezyk22}, 
absorbing boundary conditions \cite{engquist-majda77,bayliss-gunzburger-turkel82,givoli91} (and the surveys 
in \cite{givoli04}, \cite[Sec.~{3.3}]{ihlenburg98}), pole condition, \cite{hohage-nannen09,nannen-schaedle11,halla16}, 
or FEM-BEM coupling, \cite{costabel88a,sayas09,steinbach11,aurada-feischl-fuehrer-karkulik-melenk-praetorius13,carstensen-funken99,erath12},
the latter providing great geometric flexibility. 
The purpose of the present work is to provide a $k$-explicit analysis of the FEM-BEM coupling strategy for Helmholtz problems that 
has been proposed in our earlier works~\cite{MMPR20} for a conforming volume discretization and~\cite{dgfembem} for a $hp$-DGFEM volume discretization. 
For both discretizations strategies, our $k$-explicit analysis shows that the scale resolution condition~\eqref{eq:intro-scale-resolution} ensures quasi-optimality of the FEM-BEM coupling method. 

The problem under consideration is the Helmholtz equation in heterogeneous media with piecewise analytic coefficients, analytic interfaces between
the different materials, and an analytic coupling boundary separating the computational domain from a homogeneous unbounded medium. 
This setting prevents the appearance of corner or edge singularities and allows us to focus on the influence of the wavenumber $k$. The coupling strategy studied by us is taken from our previous work
\cite{MMPR20,dgfembem}. It is particularly suitable for smooth coupling boundaries $\Gamma$ and, in the limit $k = 0$, is related to the symmetric
coupling \cite{costabel88a,han90} for the Poisson problem. For a more detailed discussion of the features of our coupling strategy, we refer to \cite{MMPR20}. 
For numerical examples for coupling strategies, we refer to \cite{MMPR20} (conforming $hp$-FEM) and \cite{dgfembem} ($hp$-DGFEM).


\section{Notation and Model problem}\label{sec:notation}
Let $\Omega \subset \R^d$, $d=2$, $3$, be a bounded Lipschitz domain with analytic boundary $\Gamma:=\partial\Omega$.
We set $\Omegaext:=\R^d \setminus \overline{\Omega}$.
Assume that~$\Omega$ can be decomposed as $\overline{\Omega} = \bigcup_{j=0}^{L}{\overline{P_j}}$,
where every subdomain $P_j$ with $0\le j\le L$ is open, Lipschitz and has an
  analytic boundary, and the subdomains $P_j$, $j=0,\ldots,L$, are
  pairwise disjoint. 
  We write $\PP:=\{ P_j, \; j=0,\dots,L\}$. 

We consider the Helmholtz equation with Sommerfeld radiation condition:
\begin{subequations}
  \label{eq:model_problem}
\begin{align}
  -\fdiv(\nu \nabla u) - (\kn)^2 u &=  f   \qquad \text{in $\R^d$,}\\
  \lim_{\abs{\bx} \to \infty} |\bx|^{\frac{d-1}{2}} \big(\partial_{\abs{\bx}} - i k\big) u(\bx) &=0,
\end{align}
\end{subequations}
where $k\in\R$ with $k\ge k_0>0$, and the coefficients $n\in L^\infty(\R^d,\R)$,
$\nu\in L^{\infty}(\R^d,\R^{d\times d})$
are 
such that their restrictions to $P_j$ admit analytic extensions to $\overline{P}_j$,
$j=0,\ldots, L$ (see \eqref{eq:ass_on_nu_n} for the precise statement). We assume that on $\Omegaext$ we have $\nu_{|_{\Omegaext}}={\mathbb I}$ the identity matrix and
$n_{|_{\Omegaext}}=1$, that $\support(f)\subset \overline{\Omega}$ with $f\in L^2(\Omega)$, and that $\nu$ is symmetric and uniformly positive definite in $\Omega$, i.e., 
\[
\exists\,\nu_0>0\ \text{s.t.} \qquad \by^T \nu(\bx) \by\ge \nu_0 \abs{\by}^2 \qquad\forall \by\in \R^d,\ \text{for a.e. } \bx\in\Omega.
\]
In addition, we assume
  that $\nu{}_{|_{\Omega}}$ is analytic in a neighborhood of $\Gamma$.

 The following two relevant model problems fit into our setting:
\begin{enumerate*}[label=(\roman*)]
   \item the case of a homogeneous underlying medium, represented by $P_0$ with $\nu=n=1$ in $P_0$, and multiple
   scatterers occupying $P_j$, $1\le j\le L$;
 \item the case of a single scatterer $P_0 = \Omega$ occupying the whole~$\Omega$ with possibly variable, smooth $\nu$, $n$. 
\end{enumerate*}

We employ standard (fractional) Sobolev spaces in domains or on manifolds as introduced in 
\cite{mclean} and
introduce the following shorthands:
Given $s \in \R$ and a domain $\mathcal{O}$, with either
$\mathcal{O} \subseteq \R^d$ or $\mathcal{O} \subseteq \Gamma$,
we write 
\begin{align}
\label{eq:weighted-norms}
  \|u\|^2_{s,\mathcal{O}}:=\|u\|^2_{H^{s}(\mathcal{O})} \qquad \text{and, for $s \ge 0$, } \qquad
  \|u\|^2_{s,k,\mathcal{O}}:=|u|^2_{H^{s}(\mathcal{O})} +  k^{2s} \|u\|^2_{L^2(\mathcal{O})} .
\end{align}

We also use  the broken Sobolev norms, for $s\ge 2$,
\begin{align*}
  \|u\|^2_{s,\PP}:=\sum_{j=0}^{L}{\|u\|^2_{s,P_j}} \qquad \text{and}\qquad
    \|u\|^2_{s,k,\PP}:=\sum_{j=0}^{L}{\|u\|^2_{s,k,P_j}},
\end{align*}
and define the broken Sobolev space
\[
H_{\PP}^s(\Omega):=\{u\in H^1(\Omega):\ 
\|u\|_{s,k,\PP}< +\infty\},
\]
endowed with the $\|\cdot\|^2_{s,k,\PP}$ norm.

In most of the paper, we work with triples of Sobolev spaces of varying smoothness.
  We
  start with the energy space
  $\mathcal{V}:=H^1(\Omega) \times H^{-1/2}(\Gamma) \times H^{1/2}(\Gamma)$.
  For $s\geq 0$
  and $\bs u=(u,m,\uext) \in H^{1+s}_\PP(\Omega) \times H^{-1/2+s}(\Gamma) \times H^{1/2+s}(\Gamma)$,
  we define the norms with $s$ orders of extra smoothness as
  \begin{align*}
    \smalltriplenorm{}{\bs u}_{k,\mathcal{V},s}
    &:=\|u\|_{1+s, \PP}+k^{s+1}\|u\|_{0,\Omega}
      + \|m\|_{s-1/2, \Gamma}+k^{s}\|m\|_{-1/2,\Gamma}
      + \|\uext\|_{s+1/2, \Gamma}+k^s\|\uext\|_{1/2,\Gamma}.
  \end{align*}
  We note that the low order terms carry a $k$-weight for each additional order of Sobolev
  regularity beyond the energy space $\mathcal{V}=H^1(\Omega) \times H^{-1/2}(\Gamma) \times H^{1/2}(\Gamma)$.
  This makes them well-suited to measure the regularity of solutions of Helmholtz problems.

  We also need a second scale of norms, which measure the regularity of right-hand sides with appropriately
  $k$-weighted lower order terms: 
  For $\bs r=(r,\Rm,\Rext) \in L^2(\Omega) \times H^{3/2}(\Gamma) \times H^{1/2}(\Gamma)$
  we define
  \begin{align*}
    \smalltriplenorm{}{\bs r}_{k,\mathcal{V}',1}
    &:=\|r\|_{L^2(\Omega)}
      + \|\Rm\|_{3/2, \Gamma}+k\|\Rm\|_{1/2,\Gamma}
      + \|\Rext\|_{1/2, \Gamma}+k\|\Rext\|_{-1/2,\Gamma}.
  \end{align*}
  Here, the natural base space is $\mathcal{V}'=(H^1(\Omega))' \times H^{1/2}(\Gamma) \times H^{-1/2}(\Gamma)$, and
  the additional powers of $k$ reflect this.
  For both spaces, we use unweighted versions of the norms by supressing the subindex $k$, i.e., 
$\smalltriplenorm{}{\bs u}_{\mathcal{V},s} :=\|u\|_{1+s, \PP}
      + \|m\|_{s-1/2, \Gamma}
      + \|\uext\|_{s+1/2, \Gamma}$
and 
$\smalltriplenorm{}{\bs u}_{\mathcal{V}',1} :=\|u\|_{0, \Omega}
      + \|m\|_{3/2, \Gamma}
      + \|\uext\|_{1/2, \Gamma}$, which 
corresponds
  to the ``natural'' unweighted  norms on the products $H^{1+s}_\PP(\Omega) \times H^{-1/2+s}(\Gamma) \times H^{1/2+s}(\Gamma)$
  and $L^2(\Omega) \times H^{3/2}(\Gamma) \times H^{1/2}(\Gamma)$, respectively.
  
  For $s \in \R$, we use the notation $\langle \cdot,\cdot\rangle_{\Gamma}$ for the continuous extension
  of the $L^2$-inner product to $H^{-s}(\Gamma) \times H^{s}(\Gamma)$.
  We write $\big\langle \cdot,\cdot \big\rangle$ for the duality pair on the product space
$L^2(\Omega) \times H^{s}(\Gamma) \times H^{s-1}(\Gamma)$ and its dual.

We make the following assumption throughout the paper.
\begin{assumption}
  \label{ass:existence_global}
There exists a Lipschitz domain $\widetilde{\Omega} \supset \overline{\Omega}$
  such that for each right-hand side $f \in L^2(\R^d)$
  with $\operatorname{supp}(f) \subseteq \widetilde{\Omega}$, 
  problem~\eqref{eq:model_problem} has a unique solution $u$ 
  and there exists a constant $\beta_{0} \ge 0$ such that
  \begin{align*}
    \|u\|_{1,k,\widetilde{\Omega}} &\lesssim C(\widetilde{\Omega})\,k^{\beta_0} \|f\|_{0,\widetilde{\Omega}}.
  \end{align*}
\end{assumption}

We employ trace operators: 
The interior and exterior Dirichlet trace operators on $\Gamma$ are denoted by $\gammazint$, $\gammazext$ respectively; 
we write $\gammaoint$ and $\gammaoext$ for the interior and exterior normal derivatives on $\Gamma$, namely, $\gammaoint u:=\nabla u_{|_{\Omega}}\cdot \nGamma$ and $\gammaoext u:=\nabla u_{|_{\Omegaext}}\cdot \nGamma$. We also define $\gammanu[int] u:=\nu\nabla u_{|_{\Omega}}\cdot \nGamma$.  
%
%
For the jumps across $\Gamma$, we write 
\begin{align*}
  \djmp{u}:=\gammazint u-\gammazext u, \qquad \njmp{u}:=\gammanu[int] u -\gammaoext u . 
\end{align*}

We introduce two classes of analytic functions on a domain
$\mathcal{O} \subseteq \R^d$ (to lighten the notation, we omit the dependence on
  $k$ in the names of these spaces):
\begin{align*}
\aclass(C,\vartheta,\mathcal{O}) &:= \{v \in C^\infty(\mathcal{O})\,|\,
                                    \|\nabla^n v\|_{0,\mathcal{O}} \leq
                                    C \vartheta^n \max\{n+1,k\}^{n}
                                   \quad \forall n \in {\mathbb N}_0\}, \\
  \aclass^{\infty}(C,\vartheta,\mathcal{O}) &:= \{v \in C^\infty(\mathcal{O})\,|\,
                                     \|\nabla^n v\|_{L^{\infty}(\mathcal{O})} \leq
                                     C \vartheta^n 
                                     (n+1)^{n}
\quad \forall n \in {\mathbb N}_0\}.
\end{align*}
Here, the pointwise norm of the tensor $\nabla^n v$ is defined as $|\nabla^n v|^2 = \sum_{\alpha \in {\mathbb N}_0^{d}: |\alpha| = n} 
\frac{n!}{\alpha!} 
|D^\alpha v|^2$.
For matrix-valued functions, the derivatives are taken component-wise and $|\cdot|$
  is taken to be the spectral norm.
We also allow $\mathcal{O}=\PP$ for the piecewise
analytic case.
For example, we  
  we will assume that the coefficients in \eqref{eq:model_problem} satisfy
  $  \nu \in \aclass^{\infty}\big(C_{\nu}, \vartheta_\nu, \PP\big)$
  and $n \in \aclass^{\infty}\big(C_{n}, \vartheta_n, \PP\big)$.
  In other words,
\begin{align}
\label{eq:ass_on_nu_n}
  \forall P_j \in \PP
\colon \qquad 
  \nu{}_{|_{P_j}} \in \aclass^{\infty}\big(C_{\nu}, \vartheta_\nu, P_{j}\big) \quad \text{and } \quad
  n{}_{|_{P_j}} \in \aclass^{\infty}\big(C_{n}, \vartheta_n, P_{j}\big).
\end{align}
On the boundary $\Gamma$, we
just use the corresponding trace space. Namely,
for a fixed tubular neighborhood $\mathcal{O} \subset \Omega$ of $\Gamma$, we write
$$
\aclass(C,\vartheta,\mathcal{O},\Gamma) := \big\{ \gammazint v, v \in \aclass(C,\vartheta,\mathcal{O}) \big\} \quad \text{and} \quad
\aclass^{\infty}(C,\vartheta,\mathcal{O},\Gamma) := \big\{ \gammazint v, v \in \aclass^{\infty}(C,\vartheta,\mathcal{O}) \big\}.
$$





\subsection{Roadmap}
\label{sec:roadmap}
Since the $k$-explicit convergence analysis for the FEM-BEM coupling is rather 
involved, let us indicate the main steps for the case of a conforming discretization. 
The related, somewhat more involved DG discretization will be discussed in detail in Section~\ref{sec:convergence}.  

The FEM-BEM coupling takes the form of a three-field formulation with a volume variable 
$u \in H^1(\Omega)$ and two boundary variables $m \in H^{-1/2}(\Gamma)$, 
$\uext \in H^{1/2}(\Gamma)$, 
where the auxiliary variable $m=\gammanu[int] u + i k u$ is the impedance trace and $\uext$ 
  is the exterior Dirichlet trace, which coincides with the interior trace $\gammazint u$.
A triplet $(u,m,\uext)$ is collected in the vector ${\mathbf u} 
\in {\mathcal V}= H^1(\Omega) \times H^{-1/2}(\Gamma) \times H^{1/2}(\Gamma)$ .
The analysis of the numerical methods is carried out in the $k$-dependent norm 
\begin{equation}
\label{eq:enorm}
\enorm{\mathbf u}^2:=
\|\nu^{1/2} \nabla u\|^2_{L^2(\Omega)} + k^2 \|u\|^2_{L^2(\Omega)} + \|m\|^2_{H^{-1/2}(\Gamma)} + \|\uext\|^2_{H^{1/2}(\Gamma)}
\sim \lefttriplenorm{\mathbf u }\righttriplenorm^2_{k,\mathcal{V},0}. 
\end{equation}
Elements ${\mathbf u}$ of ${\mathcal V}$ may also be piecewise analytic. In this Subsection~\ref{sec:roadmap} we use 
the analyticity class $\AAA(M):= \aclass(M,\vartheta, \Omega\setminus\partial\PP) \times \aclass(k M,\vartheta,{\mathcal{O}},\Gamma) \times \aclass(M,\vartheta,{\mathcal{O}},\Gamma)$, 
where we track only the dependence on $M$ and suppress the dependence on $\vartheta$ and ${\mathcal O}$. That is, ${\mathbf u} \in \AAA(M)$ implies that the first component $u$ is piecewise analytic and the
second and third components admit analytic extensions to a fixed neighborhood of $\Gamma$.

The three-field variational formulation of (\ref{eq:model_problem}) reads: 
Find ${\mathbf u}  = (u,m,\uext) \in {\mathcal V}$ such that 
\begin{equation*}
\TTC{k}({\mathbf u},{\mathbf v})  = \ell({\mathbf v}) \quad \forall {\mathbf v} \in {\mathcal V}. 
\end{equation*}
for a 
sesquilinear form $\TTC{k}$ and a linear functional $\ell$ incorporating the given data  (see (\ref{form:T:for:Helmholtz})). 
The conforming Galerkin discretization based on the space ${\mathcal V}_h^C \subset {\mathcal V}$ (cf.\ Section~\ref{sec:conforming-fem}) 
then 
leads to an approximation ${\mathbf u}_N$ and the Galerkin error ${\mathbf e}_N:= {\mathbf u} - {\mathbf u}_N$. The analysis 
proceeds in the following steps. 
\begin{enumerate}[start=1,label=(\Roman*)]
\item 
\label{item:intro-1}
(G{\aa}rding inequality for $\TTC{k}$ and the sesquilinear form $\TT_+$) The form $\TTC{k}$ can be understood as a perturbation of the case $k = 0$ and then reads 
\begin{align*}
\TTC{k}({\mathbf u},{\mathbf v}) & = \TTC{0}({\mathbf u},{\mathbf v}) - k^2 (n^2 u,v)_{L^2(\Omega)} + i k (u,v)_{L^2(\Gamma)} 
+ \langle \boldsymbol{\mathcal{K}}_k  - \boldsymbol{\mathcal{K}}_0){\mathbf u},{\mathbf v} \rangle, 
\end{align*}
where the operator $\boldsymbol{\mathcal K}_k$ involves boundary integral operators that realize the coupling, and 
$\boldsymbol{\mathcal K}_0$ is  the corresponding one for $k = 0$;
  see (\ref{form:T:for:Helmholtz}) for the precise definition. 
The choice of the coupling between the PDE on $\Omega$ and on $\Omegaext$ is such that sequilinear form $\TTC{0}$ 
reproduces a variant of the classical symmetric coupling of Costabel \cite{costabel88a} and Han \cite{han90}. In particular, $\TTC{0}$ 
is positive semi-definite: $\TTC{0}({\mathbf u},{\mathbf u}) \gtrsim |u|^2_{H^1(\Omega)} + \|u\|^2_{-1/2,\Gamma} + |\uext|^2_{1/2,\Gamma}$. 
The difference $\boldsymbol{\mathcal K}_k - \boldsymbol{\mathcal K}_0$ will turn out to be a compact operator 
(see Prop.~\ref{prop:decomposition} and specifically Lemma~\ref{lemma:introducing_theta}),
so that the sesquilinear form $\TTC{k}$
satisfies the G{\aa}rding inequality (\ref{eq:intro-gaarding-neu}) below. 
We introduce the operator $\Theta$ (cf.\ Lemma~\ref{lemma:introducing_theta}) and the sesquilinear form $\TT_+$ by
\begin{align}
\label{eq:intro-def-Theta}
\langle {\mathbf u}, \Theta {\mathbf v} \rangle & = 2 k^2 (n^2 u,v)_{L^2(\Omega)} 
- \langle \boldsymbol {\mathcal K}_k  - \boldsymbol {\mathcal K}_0){\mathbf u},{\mathbf v} \rangle, \\ 
\label{eq:intro-def-TTx+}
\TT_+({\mathbf u}, {\mathbf v}) &:= \TTC{k}({\mathbf u}, {\mathbf v}) + \langle {\mathbf u}, \Theta {\mathbf v}\rangle\\
&= \TTC{0}({\mathbf u}, {\mathbf v}) +  k^2 (n^2 u, v)_{L^2(\Omega)} + i k (u,v)_{L^2(\Gamma)}
+ \langle \uext,1\rangle_{\Gamma}\langle 1, \vext \rangle_\Gamma
\end{align}
so that we have, with a constant $C > 0$ independent of $k$,  the estimate 
\begin{align}
\label{eq:intro-gaarding-neu}
\Re \big( \TT_+({\mathbf u}, {\mathbf u}) \big)= 
\Re \big(\TTC{k}({\mathbf u}, {\mathbf u}) +  \langle {\mathbf u}, \Theta{\mathbf u}\rangle \big)  \ge C \enorm{\mathbf u}^2. 
\end{align} 
This coercivity is not unexpected in view of the fact that, up to the purely imaginary term $\ii k (u,v)_{L^2(\Gamma)}$ and the
non-negative term $\langle \uext,1\rangle_{\Gamma}\langle 1, \vext \rangle_\Gamma$,
  the sesquilinear form $\TT_+$ corresponds to a symmetric FEM-BEM coupling 
for the operator $-\operatorname{div} (\nu \nabla u) + k^2 n^2 u$ in $\Omega$ and the Laplacian in $\Omegaext$. The sesquilinear form $\TT_+$ is 
also uniformly continuous: 
\begin{equation} 
\label{eq:intro-TT+-continuous}
|\TT_+({\mathbf u},{\mathbf v})| \lesssim \enorm{\mathbf u} \enorm{\mathbf v}. 
\end{equation} 
\item 
\label{item:intro-2}
(continuity properties of $\TTC{k}$) 
Due to the presence of the boundary integral operators in $\TTC{k}$, the sesquilinear form $\TTC{k}$ is not bounded uniformly in $k$. Instead, one 
has 
\begin{equation}
|\TTC{k}({\mathbf u},{\mathbf u})| \leq C k^{\mu_{stab}} \enorm{\mathbf u}\enorm{\mathbf v} , 
\qquad \text{for some } \mu_{stab} \leq 4,
\end{equation}
(cf.\ Corollaries~\ref{cor:TDC_is_bounded_with_k} and \ref{cor:TDG_is_bounded_with_k}).
However, the terms responsible for the $k$-dependence arise from the operators in $\boldsymbol{\mathcal K}_k$, which are captured by the operator $\Theta$. 
It it convenient to decompose $\Theta$ as 
\begin{equation*}
\Theta = \Theta^{\mathcal F} + \Theta^{\mathcal A}, 
\end{equation*}
where the operator $\Theta^{\mathcal F}$ is a smoothing operator of finite negative order with $k$-dependence that matches the order (cf.\ Lemma~\ref{lemma:introducing_theta}); 
in particular $\langle \cdot, \Theta^{\mathcal F}\cdot \rangle $ is bounded uniformly in $k$ with respect to $\enorm{\cdot}$. 
The operator ${\mathbf v} \mapsto \Theta^{\mathcal A}{\mathbf v}$ maps into the class of analytic functions $\AAA(C k^{\mu_{stab}} \enorm{\mathbf v})$ 
by Lemma~\ref{lemma:introducing_theta}.  In view of the uniform-in-$k$ boundedness of $\TTC{k} + \langle \cdot,\Theta \cdot\rangle$, a triangle inequality shows the uniform-in-$k$ continuity assertions
\begin{equation}
|\langle {\mathbf u}, \Theta^{\mathcal F} {\mathbf v}\rangle | + 
|\TTC{k}({\mathbf u}, {\mathbf v}) + \langle {\mathbf u}, \Theta^{\mathcal A} {\mathbf v} \rangle | 
\lesssim \enorm{\mathbf u} \enorm{\mathbf v}. 
\end{equation}
\item
\label{item:intro-3}
(dual problems and adjoint approximation property) 
For the compact operators $\Theta^{\mathcal F}$ and $\Theta^{\mathcal A}$, we introduce for given ${\mathbf w}$ the dual solutions 
${\boldsymbol \Psi}^{\mathcal F}_{\mathbf w}$, 
${\boldsymbol \Psi}^{\mathcal A}_{\mathbf w}$, 
${\boldsymbol \Psi}_{\mathbf w}$ by 
\begin{align} 
\label{eq:intro-dual-solutions}
\TTC{k}(\cdot,{\boldsymbol \Psi}^{\mathcal F}_{\mathbf w}) & = \langle \cdot, \Theta^{\mathcal F} {\mathbf w}\rangle, & 
\TTC{k}(\cdot,{\boldsymbol \Psi}^{\mathcal A}_{\mathbf w}) & = \langle \cdot, \Theta^{\mathcal A} {\mathbf w}\rangle, &
\TTC{k}(\cdot,{\boldsymbol \Psi}_{\mathbf w}) & = \langle \cdot,{\mathbf w}\rangle, 
\end{align}
which exist by our Assumption~\ref{ass:existence_global} (see Corollary~\ref{cor:existence_dual})
and the adjoint approximation properties 
\begin{align}
\label{eq:intro-adjoint-approximation-property}
  \eta^{(\mathcal F)}
  & :=
    \sup_{\bs w \in \mathcal{V}\setminus \{0\}}\inf_{{\boldsymbol \Psi}_N \in {\mathcal V}_h^C} \frac{\enorm{ {\boldsymbol \Psi}^{\mathcal F}_{\mathbf w} - {\boldsymbol \Psi}_N} }{\enorm{\mathbf w}}, 
&  
\eta^{(\mathcal A)}& :=\sup_{\bs w \in \mathcal{V}\setminus\{0\}} \inf_{{\boldsymbol \Psi}_N \in {\mathcal V}_h^C} \frac{\enorm{ {\boldsymbol \Psi}^{\mathcal A}_{\mathbf w} - {\boldsymbol \Psi}_N}}{\enorm{\mathbf w}}. 
\end{align}
\item 
\label{item:intro-4}
(Nitsche trick; cf.\ proof of Thms.~\ref{thm:dg_is_quasioptimal}, \ref{thm:fem_is_quasioptimal}) For arbitrary discrete ${\mathbf v}_N$, the G{\aa}rding inequality, the Galerkin orthogonality, and the 
uniform-in-$k$ boundedness of $\TTC{k} + \langle \cdot, \Theta^{\mathcal A} \cdot\rangle $ yield  
\begin{align}
\label{eq:intro-gaarding}
\enorm{{\mathbf e}_N}^2
& \lesssim \Re \big[\TTC{k}({\mathbf e}_N,{\mathbf u}_N) + \langle {\mathbf e}_N, \Theta {\mathbf e}_N\rangle  \big]
=  \Re \big[\TTC{k}({\mathbf e}_N,{\mathbf u} - {\mathbf v}_N) + \langle {\mathbf e}_N, \Theta {\mathbf e}_N\rangle \big] \\
\nonumber 
&= \Re \left[ \TTC{k}({\mathbf e}_N,{\mathbf u} - {\mathbf v}_N ) + \langle {\mathbf e}_N, \Theta^{\mathcal A} ({\mathbf u} - {\mathbf v}_N) \rangle 
    - \langle {\mathbf e}_N, \Theta^{\mathcal A} ({\mathbf u} - {\mathbf v}_N) \rangle 
    + \langle {\mathbf e}_N, \Theta {\mathbf e}_N\rangle   
             \right]\\
\nonumber 
&\lesssim   \enorm{{\mathbf e}_N} \enorm{ {\mathbf u} - {\mathbf v}_N }
+ | \langle {\mathbf e}_N, \Theta^{\mathcal A} ({\mathbf u} - {\mathbf v}_N) \rangle|  
+ |\langle {\mathbf e}_N, \Theta^{\mathcal F} {\mathbf e}_N\rangle |
+ |\langle {\mathbf e}_N, \Theta^{\mathcal A} {\mathbf e}_N\rangle|. 
\end{align}
The three terms 
$\langle {\mathbf e}_N, \Theta^{\mathcal A}( {\mathbf u} - {\mathbf v}_N) \rangle$, 
$\langle {\mathbf e}_N, \Theta^{\mathcal F} {\mathbf e}_N\rangle$, and  
$\langle {\mathbf e}_N, \Theta^{\mathcal A} {\mathbf e}_N\rangle$ are 
treated with duality arguments, which exploit the smoothing properties of $\Theta^{\mathcal F}$ and $\Theta^{\mathcal A}$.
We start with the term $\langle {\mathbf e}_N, \Theta^{\mathcal A}  {\mathbf e}_N \rangle$. 
With the notation (\ref{eq:intro-dual-solutions}) for dual solutions and Galerkin orthogonality satisfied by ${\mathbf e}_N$  
for arbitrary ${\boldsymbol \Psi}^{\mathcal A}_N$, we estimate rather generously with the $k$-dependent continuity of $\TTC{k}$ 
\begin{align*}
\left| \langle {\mathbf e}_N,\Theta^{\mathcal A} {\mathbf e}_N\rangle 
\right|
& = 
\left| 
\TTC{k}({\mathbf e}_N, {\boldsymbol \Psi}^{\mathcal A}_{{\mathbf e}_N} - {\boldsymbol \Psi}^{\mathcal A}_{N})  
\right|
\lesssim k^{\mu_{stab}} \enorm{{\mathbf e}_N} \enorm{ {\boldsymbol \Psi}^{\mathcal A}_{{\mathbf e}_N} - {\boldsymbol \Psi}^{\mathcal A}_N}
\end{align*}
Hence, with (\ref{eq:intro-adjoint-approximation-property})
\begin{align}
\label{eq:intro-analytic-estimate-1}
\left| 
\langle {\mathbf e}_N,\Theta^{\mathcal A} {\mathbf e}_N\rangle 
\right|
& \lesssim k^{\mu_{stab}} \eta^{(\mathcal A)}\enorm{ {\mathbf e}_N}^2. 
\end{align}
Completely analogously, we arrive at 
\begin{align}
\label{eq:intro-analytic-estimate-2}
\left| 
\langle {\mathbf e}_N,\Theta^{\mathcal A} ({\mathbf u} - {\mathbf v}_N)\rangle 
\right|
& \lesssim k^{\mu_{stab}} \eta^{(\mathcal A)} \enorm{{\mathbf e}_N}  \enorm{ {\mathbf u} - {\mathbf v}_N}. 
\end{align}
For the term $\langle {\mathbf e}_N, \Theta^{\mathcal F} {\mathbf e}_N\rangle$, we get 
with Galerkin orthogonality for arbitrary discrete ${\boldsymbol \Psi}^{{\mathcal F}}_{N}$ 
\begin{align*}
\left| 
\langle {\mathbf e}_N,\Theta^{\mathcal F} {\mathbf e}_N\rangle 
\right|
& = 
\left| 
\TTC{k}({\mathbf e}_N, {\boldsymbol \Psi}^{\mathcal F}_{{\mathbf e}_N} - {\boldsymbol \Psi}^{\mathcal F}_{N})   
\right|
\\
& = 
\left| 
\TTC{k}({\mathbf e}_N, {\boldsymbol \Psi}^{\mathcal F}_{{\mathbf e}_N} - {\boldsymbol \Psi}^{\mathcal F}_{N})   + 
\langle {\mathbf e}_N, \Theta^{\mathcal A} ( {\boldsymbol \Psi}^{\mathcal F}_{{\mathbf e}_N} - {\boldsymbol \Psi}^{\mathcal F}_{N}) \rangle  
- \langle {\mathbf e}_N, \Theta^{\mathcal A} ( {\boldsymbol \Psi}^{\mathcal F}_{{\mathbf e}_N} - {\boldsymbol \Psi}^{\mathcal F}_{N}) \rangle   
\right| \\
&\lesssim 
\enorm{ {\mathbf e}_N}  \enorm{ {\boldsymbol \Psi}^{\mathcal F}_{{\mathbf e}_N} - {\boldsymbol \Psi}^{\mathcal F}_N}
+| \langle {\mathbf e}_N, \Theta^{\mathcal A} ( {\boldsymbol \Psi}^{\mathcal F}_{{\mathbf e}_N} - {\boldsymbol \Psi}^{\mathcal F}_{N}) \rangle|. 
\end{align*}
The term $\langle {\mathbf e}_N, \Theta^{\mathcal A} ( {\boldsymbol \Psi}^{\mathcal F}_{{\mathbf e}_N} - {\boldsymbol \Psi}^{\mathcal F}_{N}) \rangle$ 
can again be treated with a duality argument analogous to (\ref{eq:intro-analytic-estimate-1}), (\ref{eq:intro-analytic-estimate-2}) to yield 
$|\langle {\mathbf e}_N, \Theta^{\mathcal A} ( {\boldsymbol \Psi}^{\mathcal F}_{{\mathbf e}_N} - {\boldsymbol \Psi}^{\mathcal F}_{N}) \rangle|
\lesssim k^{\mu_{stab}} \eta^{(\mathcal A)} \|{\mathbf e}_N\|_k \|{\boldsymbol \Psi}^{\mathcal F}_{{\mathbf e}_N} - {\boldsymbol \Psi}^{\mathcal F}_{N}\|_k$. Hence, 
since ${\boldsymbol \Psi}^{\mathcal F}_N \in {\mathcal V}_h^C$ is arbitrary,
\begin{align*}
\left| \langle {\mathbf e}_N,\Theta^{\mathcal F} {\mathbf e}_N\rangle 
\right|
& \lesssim \inf_{{\boldsymbol \Psi}^{\mathcal F}_N \in {\mathcal V}_h^C}
\left( 1  + k^{\mu_{stab}} \eta^{(\mathcal A)} \right) \enorm{ {\mathbf e}_N}  \enorm{ {\boldsymbol \Psi}^{\mathcal F}_{{\mathbf e}_N} - {\boldsymbol \Psi}^{\mathcal F}_N} 
\lesssim 
\left( 1 + k^{\mu_{stab}} \eta^{(\mathcal A)} \right) \eta^{(\mathcal F)} \enorm{ {\mathbf e}_N }^2.
\end{align*}
\item 
\label{item:intro-5}
(quasi-optimality under abstract scale resolution; cf.\ Theorems~\ref{thm:dg_is_quasioptimal}, \ref{thm:fem_is_quasioptimal})
Inserting the above estimates into (\ref{eq:intro-gaarding}) yields 
\begin{align}
\label{eq:intro-quasi-optimality}
\|{\mathbf e}_N\|^2_k &\lesssim (1 + k^{\mu_{stab}} \eta^{(\mathcal A)}) \enorm{ {\mathbf e}_N}  \enorm{ {\mathbf u} - {\mathbf v}_N}  
+ \left( ( 1+ k^{\mu_{stab}}\eta^{(\mathcal A)}) \eta^{(\mathcal F)} + k^{\mu_{stab}} \eta^{(\mathcal A)}\right) \enorm{ {\mathbf e}_N }^2. 
\end{align}
\end{enumerate}
Estimate (\ref{eq:intro-quasi-optimality}) shows that quasi-optimality is achieved if $\eta^{(\mathcal F)}$ and $k^{\mu_{stab}} \eta^{(\mathcal A)}$ are small. 
To quantify the adjoint approximation properties in terms of $k$, it is essential to understand the regularity of the dual solutions 
${\boldsymbol \Psi}^{\mathcal F}_{\mathbf w}$, ${\boldsymbol \Psi}^{\mathcal A}_{\mathbf w}$ of (\ref{eq:intro-dual-solutions}) for arbitrary ${\mathbf w} \in {\mathcal V}^1$ 
and is the core of the present work. 
This is achieved using the idea of ``regularity by decomposition'' from \cite{MS10,MS11, bernkopf} adapted to the present FEM-BEM setting. Structurally, we follow the procedure 
outlined in \cite{bernkopf} for scalar problems. This ``regularity by decomposition'' provides a $k$-explicit regularity for the solution ${\boldsymbol \Psi}^{\mathcal F}_{\mathbf w}$
of (\ref{eq:intro-dual-solutions}) in \ref{item:intro-11} that allows us to estimate $\eta^{(\mathcal F)}$ in \ref{item:intro-12} below. 
\begin{enumerate}[start=6,label=(\Roman*)]
\item 
\label{item:intro-6}
The regularity theory for the solutions 
${\boldsymbol \Psi}^{\mathcal A}_{\mathbf w}$, ${\boldsymbol \Psi}^{\mathcal F}_{\mathbf w}$ of (\ref{eq:intro-dual-solutions}) relies on \emph{a priori} bounds of the form 
\begin{equation}
\label{eq:intro-polynmomial-well-posedness}
\enorm{ {\boldsymbol \Psi}_{\mathbf w}}  \leq C k^\beta \left[\|w\|_{0,\Omega} + \|w^{m}\|_{-1/2,\Gamma} + \|w^{ext}\|_{1/2,\Gamma}\right]
\end{equation}
for the solution ${\boldsymbol \Psi}_{\mathbf w}$ of (\ref{eq:intro-dual-solutions}) 
where the constants $C$, $\beta \ge 0$ are  independent of $k$.
This solvability and stability assertion for the dual problem follows from the stipulated stability of  primal problem in Assumption~\ref{ass:existence_global} and is formalized in Lemma~\ref{lemma:existence_primal} and Corollary~\ref{cor:existence_dual}.
\item 
\label{item:intro-7}
(tool: high- and low-pass filters; cf.\ Section~\ref{sect:filters}) 
As a tool, we will require filter operators, the low-pass filter ${\mathbf L}_{\eta }$ and the high-pass filter ${\mathbf H}_{\eta }$. The tuning parameter $\eta \in(0, 1)$ plays an
important role to ensure contractivity in step \ref{item:intro-10} below. Key features of the operators are:  
\begin{enumerate}[label=\alph*)]
\item The low-pass filter ${\mathbf L}_{\eta }$ maps into a class of analytic functions.
\item The high-pass filter ${\mathbf H}_{\eta }$ acts componentwise and 
${\mathbf H}_{\eta } (w,w^m,w^{ext})=: (\widetilde w, \widetilde w^m, \widetilde w^{ext})$ satisfies for $0 \leq s' \leq s$ 
\begin{subequations}
\begin{align}
\label{eq:intro-properties-of-H}
\|\widetilde w^m\|_{{s'},\Omega} & \lesssim (\eta/k)^{s-s'} \|w\|_{s,\Omega}, 
& 
\|\widetilde w^m\|_{-1/2+s',\Gamma} & \lesssim (\eta/k)^{s-s'} \|w^m\|_{-1/2+s,\Gamma}  \\
\|\widetilde w^{ext}\|_{1/2+s',\Gamma} & \lesssim (\eta/k)^{s-s'} \|w^{ext}\|_{{1/2+s},\Gamma}. 
\end{align}
From this follow the two estimates
\begin{align}
\label{eq:into-application-of-H}
k \|\widetilde w^m \|_{1/2,\Gamma} +  
\|\widetilde w^m \|_{3/2,\Gamma} &\lesssim \|w^m\|_{3/2,\Gamma}, & 
k \|\widetilde w^{ext} \|_{-1/2,\Gamma} +  
\|\widetilde w^{ext} \|_{1/2,\Gamma} &\lesssim \|w^{ext}\|_{1/2,\Gamma},
\end{align}
\end{subequations}
which will be needed below. 
\end{enumerate}

\item 
\label{item:intro-8}
(regularity I: analyticity for analytic data; cf.\ Lemma~\ref{lemma:analytic_regularity}) 
Due to the piecewise analyticity of the data (coefficients, interfaces, boundary), the solution ${\boldsymbol \Psi}_{\mathbf w}$ is piecewise analytic if the right-hand side 
${\mathbf w}$ is. This is captured with the analyticity classes ${\mathcal A}$, and in fact, if ${\mathbf w} \in \AAA(M)$, then the solution 
${\boldsymbol \Psi}_{\mathbf w} \in \AAA(C M k^{\beta})$, where we suppress the fact that the parameters $\vartheta$ may change. 
The operator $\Theta^{\mathcal A}$ maps into the class of piecewise analytic functions (cf.\ Lemma~\ref{lemma:introducing_theta}) 
so that the corresponding solution ${\boldsymbol \Psi}^{\mathcal A}_{\mathbf w}$ of (\ref{eq:intro-dual-solutions}) is piecewise analytic and in an analyticity
class $\AAA(C k^{\beta+3} \enorm{\mathbf w})$.
\item
\label{item:intro-9}
(regularity II: regularity of the problem with the ``good'' sign)  The sequilinear form $\TT_+$ is (uniformly-in-$k$) bounded and coercive in $\enorm{\cdot}$ 
so that the Lax-Milgram Lemma and elliptic regularity theory can be brought to bear. The solution ${\boldsymbol \Psi}^+_{\mathbf w}$ of 
\begin{equation}
\label{eq:intro-plus-problem}
\TT_+(\cdot,{\boldsymbol \Psi}^+_{\mathbf w})  = \langle \cdot, {\mathbf w}\rangle 
\end{equation}
satisfies by elliptic regularity theory (see Lemma~\ref{lemma:mod_aprioi_h2} for details) 
\begin{align}
\label{eq:intro-plus-problem-regularity}
   \lefttriplenorm{\boldsymbol \Psi}^+_{\mathbf w}\righttriplenorm_{k,\mathcal{V},1}  &\lesssim
    \smalltriplenorm{}{\mathbf{w}}_{k,\mathcal{V}',1}.
\end{align}
\item
\label{item:intro-10}
(regularity by decomposition, cf.\ Theorem~\ref{thm:regularity_splitting}) 
For right-hand sides ${\mathbf w}$ with finite regularity, the solution ${\boldsymbol \Psi}_{\mathbf w}$ of 
(\ref{eq:intro-dual-solutions}) is decomposed as ${\boldsymbol \Psi}_{\mathbf w} = {\boldsymbol \Psi}_{H^2} + {\boldsymbol \Psi}_{\mathcal A}$, where 
the term ${\boldsymbol \Psi}_{\mathcal A}$ is again in an analyticity class and 
${\boldsymbol \Psi}_{H^2} \in H^2_{\PP}(\Omega) \times H^{1/2}(\Gamma) \times H^{3/2}(\Gamma)$ with 
\begin{equation}
\lefttriplenorm{\boldsymbol \Psi}_{H^2}\righttriplenorm_{k,\mathcal{V},1} 
\lesssim \smalltriplenorm{}{\mathbf w}_{\mathcal{V'},1}. 
\end{equation}
This decomposition is achieved as follows. We set $(r_0,r^m_0,r^{ext}_0):= {\mathbf r}_0:= {\mathbf w}$ and define, using the solution operators 
(\ref{eq:intro-plus-problem}), (\ref{eq:intro-dual-solutions}), the functions 
${\boldsymbol \Psi}_{H^2,0}:= {\boldsymbol \Psi}^+_{{\mathbf H}_{\eta }({\mathbf r}_0)}$ and 
${\boldsymbol \Psi}_{{\mathcal A},0}:= {\boldsymbol \Psi}_{{\mathbf L}_{\eta }({\mathbf r}_0)}$. 
The functions 
${\boldsymbol \Psi}_{{\mathcal A},0}$ is in an analyticity class by \ref{item:intro-8} and the function 
${\boldsymbol \Psi}_{H^2,0}$ satisfies the regularity (we write ${\mathbf H}_{\eta k} {\mathbf r}_0 = (\widetilde r, \widetilde r^m, \widetilde r^{ext})$) 
\begin{align*} 
 \lefttriplenorm{\boldsymbol \Psi}_{H^2,0}\righttriplenorm_{k,\mathcal{V},1} 
  & \stackrel{(\ref{eq:intro-plus-problem-regularity}) } {\lesssim}
     \smalltriplenorm{}{{\mathbf H}_{\eta k} {\mathbf r}_0 }_{k,\mathcal{V}',1}  
   \stackrel{(\ref{eq:into-application-of-H})} {\lesssim }\|w\|_{0,\Omega} + \|r^m \|_{3/2,\Gamma} + \|r^{ext} \|_{1/2,\Gamma}
  =\smalltriplenorm{}{\bs w}_{\mathcal{V}',1}.  
\end{align*} 
The key observation in the contraction argument of Lemma~\ref{lemma:decomposition_single_step} is that the difference ${\boldsymbol \Psi}_{\mathbf w} -  \left( 
{\boldsymbol \Psi}_{H^2,0} + {\boldsymbol \Psi}_{{\mathcal A},0}\right)$ leads to a residual ${\mathbf r}_1 = (r_1,r^m_1, r^{ext}_1) $  defined by 
$$
\langle \cdot, {\mathbf r}_1\rangle = \TT_k( \cdot, {\boldsymbol \Psi}_{\mathbf w} -  \left( {\boldsymbol \Psi}_{H^2,0} + {\boldsymbol \Psi}_{{\mathcal A},0}\right)) 
$$
with the contraction property 
\begin{align*}
  \smalltriplenorm{}{{\mathbf r}_1}_{\mathcal{V}',1}
  \leq C \eta^{1/4} \smalltriplenorm{}{{\mathbf r}_0}_{\mathcal{V}',1}. 
\end{align*}
Thus, by making $\eta$ sufficiently small, one obtains  the estimate
$\|{\mathbf r}_1\|_{L^2(\Omega) \times H^{3/2}(\Gamma) \times H^{1/2}(\Gamma)} \leq \frac{1}{2} \|{\mathbf r}_0\|_{L^2(\Omega) \times H^{3/2}(\Gamma) \times H^{1/2}(\Gamma)} $. 
The above argument can be repeated for ${\mathbf r}_1$ and thus, iteratively, one arrives at a convergent series 
${\boldsymbol \Psi}_{H^2}:= \sum_{i} {\boldsymbol \Psi}_{H^2,i}$ and 
${\boldsymbol \Psi}_{\mathcal A} = \sum_{i} {\boldsymbol \Psi}_{{\mathcal A},i}$. 
\item
\label{item:intro-11}
(regularity of ${\boldsymbol \Psi}^{\mathcal F}_{\mathbf w}$) The regularity result \ref{item:intro-10} implies a decomposition 
of ${\boldsymbol \Psi}^{\mathcal F}_{\mathbf w} = 
{\boldsymbol \Psi}^{\mathcal F} _{H^2}  + 
{\boldsymbol \Psi}^{\mathcal F} _{\mathcal A}$, where   
${\boldsymbol \Psi}^{\mathcal F} _{\mathcal A}$ is in an analyticity class. For the first term, 
${\boldsymbol \Psi}^{\mathcal F} _{H^2}$, we observe from Lemma~\ref{lemma:introducing_theta} that the first component 
of $\Theta^{\mathcal F} {\mathbf w} = (\widetilde w, \widetilde w^m, \widetilde w^{ext})$ 
is given by $\widetilde w = 2 k^2 n^2 w$ and the other two components can be bounded by 
\begin{align*}
\|\widetilde w^m\|_{3/2,\Gamma} + \|\widetilde w^{ext} \|_{1/2,\Gamma} &\lesssim k \left[ \|w^m\|_{-1/2,\Gamma} + \|w^{ext}\|_{1/2,\Gamma} \right].
\end{align*}
so that 
$$
\smalltriplenorm{}{\Theta^{\mathcal F} {\mathbf w}}_{\mathcal{V}',1}
\lesssim k \|w\|_{0,\Omega} + k \|w^m \|_{-1/2,\Gamma} + k \|w^{ext}\|_{1/2,\Gamma}
\lesssim\smalltriplenorm{}{\bs w}_{k,\mathcal{V},0}. 
$$
We conclude that from \ref{item:intro-10}
\begin{equation*}
  \smalltriplenorm{}{{\boldsymbol \Psi}^{\mathcal F}_{\mathbf w}}_{k,\mathcal{V},1}
\lesssim k \|w\|_{0,\Omega} + k \|w^m \|_{-1/2,\Gamma} + k \|w^{ext}\|_{1/2,\Gamma} 
\lesssim k \enorm{\mathbf w}.  
\end{equation*}

\item
\label{item:intro-12}
(estimating $\eta^{(\mathcal F)}$ and $\eta^{(\mathcal A)}$)
In this last step, we relate the adjoint approximation properties $\eta^{(\mathcal F)}$ and $\eta^{(\mathcal A)}$, which involve the solution 
operators of adjoint problems, by standard approximation properties of piecewise polynomial space for Sobolev functions. With the space 
(cf.\ Section~\ref{sec:convergence}) 
${\mathcal V}^s:= H^{s+1}_{\PP}(\Omega) \times H^{-1/2+s}(\Gamma) \times H^{1/2+s}(\Gamma)$ (equipped with the natural norm) we set 
\begin{align}
\label{eq:intro-approximation-properties}
  \eta^{(s)}_C
  & :=
    \sup_{0 \ne {\mathbf u} \in {\mathcal V}^s} \inf_{{\mathbf v}_N \in {\mathcal V}_h^C} \frac{\enorm{{\mathbf u} - \mathbf v_N}}{
    \smalltriplenorm{}{\mathbf u}_{k,{\mathcal V},s}}, 
& 
\eta^{(\mathcal A)}_C& := \sup_{{\mathbf u} \in {\mathcal A}(1)} \inf_{{\mathbf v}_N \in {\mathcal V}_h^C} \enorm{{\mathbf u} - \mathbf v_N}. 
\end{align}
In view of items \ref{item:intro-8} and \ref{item:intro-11}, we see that we may estimate 
\begin{align*}
\eta^{(\mathcal A)} & \leq C k^{\beta+4} \eta^{(\mathcal A)}_C, & 
\eta^{(\mathcal F)} & \leq C k \eta^{(1)}_C + C k^{\beta+5}\eta^{(\mathcal A)}_C, 
\end{align*}
so that estimating $\eta^{(\mathcal A)}$ and $\eta^{(\mathcal F)}$ of (\ref{eq:intro-adjoint-approximation-property}) is reduced to controlling
the more established quantities (\ref{eq:intro-approximation-properties}). 
In the present work, the space ${\mathcal V}_h$ consists of piecewise polynomials of degree $p$ (or $p-1$)
which leads to exponential (in $p$) convergence for $\eta^{(\mathcal A)}_C$ and algebraic rates $(h/p)^{s}$ for $\eta^{(s)}_C$. The condition that 
$k \eta^{(1)}_C$ be small then leads to the condition that $kh/p$ be small, and the condition that
$k^{\mu_{stab}} \eta^{(\mathcal A)} \leq k^{\beta+8}\eta^{(\A)}_C$ and $k^{\beta+5} \eta^{(\A)}_C$ be small leads to the 
condition $p \gtrsim \log k$, which is the scale resolution condition (\ref{eq:intro-scale-resolution}). Formal statements can be found in 
Cor.~\ref{cor:dg_explicit_convergence} for $hp$-DGFEM and in Cor.~\ref{cor:conforming_explicit_convergence} for the $hp$-FEM discretization.  
\end{enumerate}


\section{Boundary integral operators and mapping properties}
In this section, we introduce the pertinent boundary integral operators for the
Helmholtz equation and recall their mapping properties. As
$n\equiv 1$ outside of $\Gamma$,
we start with the free space Green's function for the Helmholtz
operator with wavenumber $k$
\[
G_k(\bx,\by):=\begin{cases}
 \frac{i}{4}H^{(1)}_0(k |\bx-\by|), & d=2,\\[0.2cm]
  \displaystyle{\frac{e^{\ii k |\bx-\by|}}{4 \pi |\bx-\by|}},& d=3,
\end{cases}  
\]
where $H^{(1)}_0$ denotes the Hankel function of the first kind of
order $0$.
In the special case $d=2$ and $k=0$ the Green's function is instead given
  by $G_0(\bx,\by):=-\frac{1}{2\pi}{\ln|\bx-\by|}$.
We define the single and double layer potentials by
\begin{align}
  \label{eq:def_potentials}
  \Vtildek \varphi(\bx)&:=\int_{\Gamma}{G_{k}(\bx-\by) \varphi(\by)\,ds(\by)} \quad \text{and}\quad
  \Ktildek \varphi(\bx):=\int_{\Gamma}{\partial_{\nGamma(\by)}G_{k}(\bx-\by) \varphi(\by)\,ds(\by)}.
\end{align}

These potentials induce the four boundary integral operators $\Vk$, $\Kk$, $\Wk$, and $\Kprimek$:
\begin{alignat*}{4}
  \text{``single layer'': } \;\;\Vk \varphi&:=\gammazint (\Vtildek\varphi),  \quad
  &\text{``double layer'': } \;\;   &(\nicefrac{-1}{2} + \Kk) \varphi:= \gammazint (\Ktildek\varphi),  \\
  \text{``hypersingular'': } \;\;\Wk \varphi&:=-\gammaoint (\Ktildek\varphi),  \quad
  &\text{``adjoint double layer'': }&(\nicefrac{1}{2} + \Kprimek) \varphi:= \gammaoint (\Vtildek\varphi).
\end{alignat*}


The following jump conditions are valid \cite[Thm.~{6.11}]{mclean}:
\begin{equation}
  \label{eq:jump_conditions}
  \begin{alignedat}{4}
  \djmp{\Vtildek \varphi}&=0,  \qquad \text{and}\qquad & \djmp{\gammao \Vtildek \varphi}&=\varphi, \\
  \djmp{\Ktildek \varphi}&=-\varphi, \qquad \text{and}\qquad & \djmp{\gammao \Ktildek \varphi}&=0.
\end{alignedat}
\end{equation}

For the boundary integral operators associated with the
Laplacian, i.e., in the special case $k=0$, the following mapping
properties hold true.
\begin{proposition}[{\cite[Thm.~{7.2}]{mclean}}]
  \label{prop:mapping properties_laplace_smooth}
  Let $\Gamma$ be analytic. Then, for any $s \in \R$:
  \begin{equation} \label{eq:mapping_laplace_smooth}
    \begin{split}
      & \Vz: H^{-1+s}(\Gamma) \rightarrow H^{s}(\Gamma), \quad\quad \quad \quad
      \Kz: H^{s}(\Gamma) \rightarrow H^{s}( \Gamma),\\
      & \Kprimez: H^{-s}(\Gamma) \rightarrow H^{-s}(\Gamma), \quad \quad\quad \quad
      \Wz: H^{s}(\Gamma) \rightarrow H^{-1+s}( \Gamma).\\        
    \end{split}
  \end{equation}
\end{proposition}

The Laplace potentials have  the following mapping properties.
\begin{proposition}[{\cite[Cor.~{6.14}]{mclean}}]
\label{prop:mapping_properties_laplace_potentials}
Let $\Gamma$ be analytic. Then for
$ s > -\nicefrac{1}{2}$ and $\mathcal{O} \subseteq \R^d$ open and bounded
\begin{align*}
&\Vtildez: H^{s-1/2}(\Gamma) \rightarrow H^{1+s}(\mathcal{O}\setminus \Gamma)
\quad \mbox{ and } \quad 
        \Ktildez: H^{s+1/2}(\Gamma) \rightarrow H^{1+s}(\mathcal{O} \setminus \Gamma).
\end{align*}
\end{proposition}

We will rely on the following decomposition result, splitting
boundary integral operators into operators with good $k$-dependence,
and an analytic remainder.
\begin{proposition}[{\cite[Lem.~{A.1}]{MMPR20}}]
\label{prop:decomposition}
Let $\Gamma$ be analytic and $k \ge k_0 > 0$.
Then there are
bounded linear operators 
$\SV$, $\SK$, $\SKprime$, $\SW$ and linear maps 
$\AVtilde: H^{-\frac32}(\Gamma) \rightarrow C^\infty(\overline\Omega)$, 
$\AKtilde: H^{-\frac12}(\Gamma) \rightarrow C^\infty(\overline\Omega)$
such that 
\begin{subequations}
\label{eq:differences}
\begin{align}
\label{eq:differences-a}
\Vk - \Vz & = \SV+ \gammazint \AVtilde, \\
\label{eq:differences-b}
\Kprimek - \Kprimez & = \SKprime + \gammaoint \AVtilde, \\ 
\label{eq:differences-c}
\Kk - \Kz & = \SK + \gammazint \AKtilde,  \\
\label{eq:differences-d}
\Wk - \Wz & = \SW - \gammaoint \AKtilde.
\end{align}
\end{subequations}
For $s \ge -1$ and for constants $C_{s,s'}$, $C_{\V}$, $C_{\K}$, $\vartheta_{\V}$,
  $\vartheta_{\K}> 0$ independent of $k \ge k_0$,
the operators $\SV$, $\SK$, $\SKprime$, $\SW$, $\AVtilde$, $\AKtilde$ 
have the 
mapping properties 
\begin{subequations}
\label{eq:estimates-S}
\begin{align}
\label{eq:estimates-S-a}
\|\SV \|_{H^{-\frac12+s'}(\Gamma)\leftarrow H^{-\frac12+s}(\Gamma)} 
\leq C_{s,s'} k^{-(1+s-s')}, \qquad 
1/2 < s' \leq  s+3, \\
\label{eq:estimates-S-b}
\|\SKprime \|_{H^{-\frac32+s'}(\Gamma)\leftarrow H^{-\frac12+s}(\Gamma)}
\leq C_{s,s'} k^{-(1+s-s')},  \qquad 
3/2 < s' \leq  s +3, \\
\label{eq:estimates-S-c}
\|\SK \|_{H^{-\frac12+s'}(\Gamma)\leftarrow H^{+\frac12+s}(\Gamma)} 
\leq C_{s,s'} k^{-(1+s-s')}, \qquad 
1/2 < s' \leq  s+3, \\
\label{eq:estimates-S-d}
\|\SW \|_{H^{-\frac32+s'}(\Gamma)\leftarrow H^{+\frac12+s}(\Gamma)} 
\leq C_{s,s'} k^{-(1+s-s')}, \qquad 
3/2 < s' \leq  s +3, \\ 
\label{eq:estimates-S-e}
\AVtilde \varphi   \in \aclass(C_{\V} k \|\varphi\|_{-\frac{3}{2},\Gamma}, \vartheta_{\V}, \Omega) 
\qquad \forall \varphi \in H^{-\frac{3}{2}}(\Gamma), \\
\label{eq:estimates-S-f}
\AKtilde \psi   \in \aclass(C_{\K} k \|\psi\|_{-\frac{1}{2},\Gamma}, \vartheta_{\K}, \Omega)
\qquad \forall \psi \in H^{-\frac{1}{2}}(\Gamma). 
\end{align}
\end{subequations}
\end{proposition}

\section{Filter operators}
\label{sect:filters}
In this section, we collect some results on filter operators that decompose
functions into high- and low-frequency contributions.
Versions of these operators have already been introduced in~\cite{MS11}, \cite{M12}, and~\cite{bernkopf}.
We will use the constructions of~\cite{M12} as they closely match our requirements.

\begin{proposition}[{frequency filters on domains, \cite[Prop.~{3.2}]{bernkopf}}]
  \label{prop:ff_volume}
  Let $\eta \in (0,1)$ and $\Omega$ be a bounded Lipschitz domain.
  Then, there exist linear operators $\hff{\Omega}$ and $\lff{\Omega}$
  defined on $L^2(\Omega)$ with the following properties:
  \begin{enumerate}[label=(\roman*)]
  \item 
\label{item:prop:ff_volume-i}
$\hff{\Omega} + \lff{\Omega}=\id$.
  \item $\|\hff{\Omega}f\|_{s',\PP} \leq C_{s,s'} (\eta
    k^{-1})^{s-s'}\|f\|_{s,\PP}$, where $0\leq s' \leq s$,\\[0.2cm]
    $\|\hff{\Omega}f\|_{\left(H^\tau(\Omega)\right)'}\le C_\tau(\eta
    k^{-1})^{\tau}\|f\|_{0,\Omega} 
    $, where $0\le \tau< 1/2$.
  \item  
\label{item:prop:ff_volume-ii}
$\hff{\Omega} + \lff{\Omega}=\id$.
$\lff{\Omega} f \in \A(C \|f\|_{0,\Omega},\vartheta_\eta,\PP)$ is a piecewise analytic function. 
  \end{enumerate}
  Here $C_{s,s'}$ and $C$ are independent of $k \geq k_0$ and $\eta$. The constant $\vartheta_\eta$
  is independent of $k \geq k_0$ but may depend on $\eta$.
\end{proposition}

\begin{proposition}[{frequency filters on surfaces, \cite[Lem.~{6.2}]{M12}}]
  \label{prop:ff_gamma1}
  Let $\Omega$ be a bounded Lipschitz domain with smooth boundary $\Gamma$.
  Let $s>0$ and $\eta \in (0,1)$. Then, there exist operators
  $\hffg{+}: H^{s}(\Gamma) \to H^{s}(\Gamma)$ and
  $\lffg{+}: H^{s}(\Gamma) \to H^{1/2+s}(\R^d)$ with the following properties:
  \begin{enumerate}[label=(\roman*)]
    \item 
  \label{item:prop:ff_gamma1-i}
$\hffg{+} + \gamma^{\mathrm{int}}_0 \lffg{+}=\id$.
    \item 
  \label{item:prop:ff_gamma1-ii}
$\|\hffg{+}f\|_{s',\Gamma} \leq C_{s,s'} (\eta k^{-1})^{s-s'}\|f\|_{s,\Gamma}$, where $0\leq s' \leq s$. \hfill\refstepcounter{equation}\textup{(\theequation)}%
      \label{eq:ff_gamma1_2}
  \item  
  \label{item:prop:ff_gamma1-iii}
$\lffg{+} f$ is an entire function on $\R^d$
    with $\|\lffg{+}f\|_{0,\R^d}\le C \|f\|_{s,\Gamma}$ and
    $$
    \|\nabla^p \lffg{+} f\|_{0,\R^d}
    \leq C_\eta (\vartheta_\eta k )^{p-(1/2+s)} \|f\|_{s,\Gamma} \qquad \forall p \in \N_0,
     \ p \geq s+ 1/2.
    $$
  \end{enumerate}
  Here $C_{s,s'}$ and $C$ are independent of $k \geq k_0$ and $\eta$. The constants
  $C_{\eta}$ and $\vartheta_\eta$
  are independent of $k \geq k_0>0$ but depend on $\eta$.
\end{proposition}


Finally, we also need a filter on the boundary that works in negative norms.
\begin{proposition}[{frequency filters on surfaces, negative norms, \cite[Lem.~{6.3}]{M12}}]
  \label{prop:ff_gamma2}
  Let $\Omega$ be a bounded Lipschitz domain with analytic boundary $\Gamma$.
  Fix $\eta \in (0,1)$. Then, there exist linear operators
  $\hffg{-}$ and $\lffg{-}$ defined on $H^{-1}(\Gamma)$ with the following properties:
  \begin{enumerate}[label=(\roman*)]
  \item $\hffg{-} +  \lffg{-}=\id$.
    \label{it:hffgn_sum}
  \item
    \label{it:hffgn_est}
    $\|\hffg{-} f\|_{s',\Gamma} \leq C_{s,s'} (\eta k^{-1})^{s-s'}\|f\|_{s,\Gamma}$, 
      where $-1\leq s'\leq s \leq 1$. \hfill\refstepcounter{equation}\textup{(\theequation)}%
      \label{eq:ff_gamma2_2}
    \item
      \label{it:lffgn_est}
      $\lffg{-} f$ is the restriction to $\Gamma$ of a function that is
      analytic in a tubular neighborhood $T$ of~$\Gamma$ and satisfies
    $$
    \|\nabla^p \lffg{-} f\|_{0,T}
    \leq C_\eta k^{d/2}\vartheta_\eta^{p}\max(k,p)^p \|f\|_{-1/2,\Gamma} \qquad \forall p \in \N_0.
    $$
  \item
    \label{it:hffgn_mean}
    The function $\hffg{-} f $ has vanishing integral mean, i.e.,
    $ (\hffg{-}f,1 )_{0,\Gamma}=0$.    
  \end{enumerate}
  Here $C_{s,s'}$ is independent of $k \geq k_0$ and $\eta$. The constants
  $C_{\eta}$ and $\vartheta_\eta$  are independent of $k \geq k_0>0$ but may depend on $\eta$.
\end{proposition}
\begin{proof}
  Items~\ref{it:hffgn_sum}--\ref{it:lffgn_est} are proven in~\cite[Lem.~{6.3}]{M12}.
  Property~\ref{it:hffgn_mean} follows by inspection of the construction given there.
  To give more details: The construction is done via the spectral decomposition of the Laplace-Beltrami operator,
  with low-frequency modes 
  included in $\lffg{-}f$ and high-frequency modes 
  in~$\hffg{-}f$.

  If $\{(\lambda_j,\varphi_j)\}_{j=0}^{\infty}$ are the eigenpairs of the Laplace-Beltrami operator 
  with $\{\varphi_j\}_{j=0}^{\infty}$ being an $L^2(\Gamma)$-orthonormal basis of eigenfunctions,
  we get that either $\lambda_j=0$ or $(1,\varphi_j)_{0,\Gamma}=0$. This follows from the
  property
  \begin{align*}
    0&=-\Delta_{\Gamma} 1 = \sum_{j=0}^{\infty}{\lambda_j (1,\varphi_j)_{0,\Gamma} \varphi_j},
  \end{align*}
  and thus, by orthogonality, it holds true that $\lambda_j (1,\varphi_j)_{0,\Gamma} =0$ for all $j \in \N_0$.
  Since the eigenfunctions corresponding to $\lambda_j=0$ are
  always included in $\lffg{-}f$, we get that the integral mean of the high-frequency part vanishes.
\end{proof}

We remark that Proposition~\ref{prop:ff_gamma2}\ref{it:lffgn_est} can also be written as
  \[
    \lffg{-} f\in
    \aclass(C k^{d/2} \|f\|_{-1/2,\Gamma}, \vartheta_\eta,
    \Gamma)\qquad\forall f\in H^{-1/2}(\Gamma),
  \]
where $C_{\eta}$ and $\vartheta_\eta$
are independent of $k \geq k_0>0$ but depend on $\eta$.

\section{Discretization}
\label{sec:discretization}
We follow~\cite{MMPR20} and rewrite problem~\eqref{eq:model_problem} 
using a mortar formulation and boundary integrals.
We obtain three coupled problems for~$u: \Omega \rightarrow \mathbb C$ and $\uext$, $m:\Gamma \rightarrow \mathbb C$: 
\begin{subequations}
    \label{system:Helmholtz}
\begin{align}
& \begin{cases}
- \fdiv(\Ad \nabla  u)    - (\kn)^2 u = f & \text{in } \Omega,\\
\gammanu[int]{u} + i k  u - m = 0 & \text{on } \Gamma,\\
\end{cases} \label{first:system:Helmholtz}\\
& \begin{cases}
u - \left[ \left( \nicefrac{1}{2} + \Kk  \right) \uext - \Vk (m - \ii k \uext)   \right] = \rhsm\quad\, \text{on }\Gamma,
\end{cases} \label{third:system:Helmholtz} \\
  & \begin{cases}
    \Bk  \uext + i k \Aprimek (\uext) -   \Aprimek m = \rhsext & \, \;\quad\text{ on } \Gamma
\end{cases},\label{second:system:Helmholtz}
\end{align}
\end{subequations}
with boundary operators $\Bk$ and $\Aprimek$ defined as
\begin{equation} \label{definition:Bk:and:Akprime}
\Bk := -\Wk - \ii k \left( \nicefrac{1}{2} - \Kk \right),\quad \quad \Aprimek := \nicefrac{1}{2} + \Kprimek + \ii k  \Vk.
\end{equation}
Here, in the interest of generality, we allow for right-hand sides $\rhsm \in H^{1/2}(\Gamma)$ and
  $\rhsext \in H^{-1/2}(\Gamma)$. 
We will usually collect the three fields $u$, $m$, $\uext$ in the vector valued quantity $\bs u :=(u,m,\uext)$ 
and we will write $\bs u=:S_{k}^{-}(f,\rhsm,\rhsext)$ for the solution of the problems.
The equivalence of \eqref{system:Helmholtz} and \eqref{eq:model_problem} is valid for $\rhsm=\rhsext=0$.

This system of equations can be discretized in multiple ways.
In~\cite{MMPR20}, a combination of conforming FEM and BEM was discussed.
Building on that work, \cite{dgfembem} considered a discontinuous
Galerkin approximation for the interior subproblem.
In the present article, we refine the analysis of \cite{MMPR20} and \cite{dgfembem} in that 
we give a $k$-explicit analysis of both discretization schemes.  
To that end, we introduce some notation concerning general
discretizations with finite element and boundary element spaces 
in Sections~\ref{sec:fem-bem-spaces}--\ref{sec:DG}, following 
closely the notation of~\cite{dgfembem}.

  Before we can prove a result on the well-posedness of problem~\eqref{system:Helmholtz}
  with polynomial bounds, we need the following preparatory lemma.
  \begin{lemma}
    \label{lemma:neumann_trace}
    For $u \in H^1(\Omega)$ with $-\laplace u - k^2 u =0$, we have the following
    trace estimates:
    \begin{align*}
      \|\nabla_{\Gamma} u\|_{-1/2,\Gamma}+ \|\gammaoint u\|_{-1/2,\Gamma}  + \|\gammanu[int]u\|_{-1/2,\Gamma}
      &\lesssim \sqrt{k} \|u\|_{1,k,\Omega}.
    \end{align*}
  \end{lemma}
  \begin{proof}
    The estimate of~$\|\gammaoint u\|_{-1/2,\Gamma}$
    follows from \cite[Lem.~{15}]{LS06} using $s= \ii k$.
      To estimate the $\nu$-weighted flux 
      $\gammanu[int]u$,
      we decompose it into the normal and tangential part:
      $$
      \nGamma \cdot \nu \nabla u
      = (\nGamma \cdot \nabla u) (\nGamma^T\nu \nGamma)
      + (\boldsymbol{\tau}_{\Gamma} \cdot \nabla u) (\nGamma^T\nu {\boldsymbol{ \tau}_{\Gamma}}).
      $$
      Since the geometry and $\nu$ are smooth and using the fact that the tangential derivative
      is a differential operator of order one, we have
      \begin{align*}
        \|\nGamma \cdot \nu \nabla u\|_{-1/2,\Gamma}
        &\lesssim \| \gammaoint u\|_{-1/2,\Gamma} + \|\nabla_{\Gamma} u\|_{-1/2,\Gamma} \\
       & \lesssim \| \gammaoint u\|_{-1/2,\Gamma} + \|u\|_{1/2,\Gamma}
        \lesssim \sqrt{k} \|u\|_{1,k,\Omega},\end{align*}
      where in the last step we used the already estalished bound on $\gammaoint u$ and
      a standard trace theorem.
      This also implies the estimate of the tangential gradient. 
  \end{proof}
  
\begin{lemma}[polynomial well-posedness]
  \label{lemma:existence_primal}
  There exists a constant $\beta \ge 0$ such that, for any right-hand
  side $(f,\rhsm,\rhsext) \in L^2(\Omega) \times
  H^{1/2}(\Gamma)
  \times H^{-1/2}(\Gamma)$, the problem 
  \eqref{system:Helmholtz}
  has a unique solution $\bs u = (u,m,\uext) \in {\mathcal V} = H^1(\Omega) \times H^{-1/2}(\Gamma) \times H^{1/2}(\Gamma)$, which 
  satisfies the estimate
  \begin{align*}
    \|u\|_{1,k,\Omega} + \|m\|_{-1/2,\Gamma} + \|\uext\|_{1/2,\Gamma}
    &\lesssim k^{\beta} \Big(\|f\|_{0,\Omega}
      +\|\rhsm\|_{1/2,\Gamma} +\|\rhsext\|_{-1/2,\Gamma} \Big).
  \end{align*}
  The constant $\beta$ can be bounded by $\beta\leq\beta_0+4$, where
  $\beta_0$ is given by Assumption~\ref{ass:existence_global}.
\end{lemma}
\begin{proof}
  We use a  construction similar to that 
  in the proof of~\cite[Thm.~{3.11}]{MMPR20}.
  Existence of the solution $(u,m,\uext)$ follows from standard Fredholm theory as the sesquilinear form satisfies a G{\aa}rding inequality 
and  uniqueness is implied by Assumption~\ref{ass:existence_global}. We thus focus on the $k$-explicit stability estimate.


  By assumption, $\nu|_{\Omega}$ is smooth near $\Gamma$. Let $\mathcal{N}(\Gamma)\subset \Omega$ be a tubular neighborhood of $\Gamma$ such that $\nu|_{\Omega}$ is smooth on $\mathcal{N}(\Gamma)\cap \Omega$. Let $\chi_{\mathcal{ N}(\Gamma)} \in C^\infty(\R^3)$ satisfy $\operatorname{supp} \chi_{\mathcal{N}(\Gamma)} \subset {\mathcal N}(\Gamma) \cup \Omegaext$ and $\chi_{\mathcal{N}(\Gamma)} \equiv 1$ in a neighborhood of $\Gamma \cup \Omegaext$. 
  We  use the following
  auxiliary functions
  \begin{align*}
    \daleth&:=\Vtildek (-m +  \ii k \uext) + \Ktildek \uext
             \qquad \text{and  } \qquad
             \LL:=u \mathds{1}_{\Omega}+ \daleth \chi_{\mathcal{N}(\Gamma)}.
  \end{align*}
  By the jump conditions (\ref{eq:jump_conditions}) satisfied by the operators, we deduce that
  \begin{align}
\label{eq:lemma:existence_primal-10}
    \djmp{\gammao \daleth + \ii k  \daleth}=-m, \qquad \djmp{\daleth}=-\uext.
  \end{align}
  Applying trace operators and using the equations~\eqref{system:Helmholtz}  we obtain that $\daleth$ solves
  \begin{align*}
    -\laplace \daleth - k^2 \daleth &=0\quad \text{in $\R^d\setminus\Gamma$}, \; \quad \gammaoint{\daleth} + \ii k  \gammazint\daleth
                                      = \rhsext \quad \text{and}\quad
    \gammazext \daleth = -\rhsm + \gammazint u\quad \text{on $\Gamma$}.                 
  \end{align*}
  Restricted to $\Omega$, this is a standard Robin boundary value problem, for which
  $k$-explicit bounds are available \cite[Cor.~{1.10}]{BSW16}:
  \begin{align}
\label{eq:lemma:existence_primal-20}
    \|\daleth\|_{1,k,\Omega} &\lesssim  k \|\rhsext\|_{-1/2,\Gamma}.
  \end{align}
  
  It is then an easy calculation that $\LL$ solves the following transmission problem:
  \begin{align*}
    - \fdiv(\nu \nabla \LL) - (\kn)^2 \LL
    &= f-\big[\fdiv(\nu \nabla (\daleth \, \chi_{\mathcal{N}(\Gamma)})  - (\kn)^2
      \daleth \, \chi_{\mathcal{N}(\Gamma)}  \big] \mathds{1}_{\Omega}, \\
    \njmp{\LL+ i k \LL}&= \gammanu[int]{\daleth} - \gammaoint{\daleth}
        \quad \text{and} \quad
    \djmp{\LL}=\gammazint \daleth + \rhsm,
  \end{align*}
  
  The impedance jump
  $R_1:=\gammanu[int]{\daleth} - \gammaoint{\daleth}$,
  the Dirichlet-jump $R_0:=\rhsm + \gammazint{\daleth}$,
  and the right-hand side $\tilde f$
  satisfy by Lemma~\ref{lemma:neumann_trace} and standard trace estimates:
  \begin{align*} 
    \|R_1\|_{-1/2,\Gamma}
    &\stackrel{\text{L.~\ref{lemma:neumann_trace}}}{ \lesssim}
      k^{1/2}\|\daleth\|_{1,k,\Omega} 
      \stackrel{\eqref{eq:lemma:existence_primal-20}}{\lesssim} k^{3/2}\|\rhsext\|_{-1/2,\Gamma},\\
    \|R_0\|_{1/2,\Gamma}&\lesssim
                          \|\daleth\|_{1,k,\Omega} +\|\rhsm\|_{1/2,\Gamma}  
\stackrel{\eqref{eq:lemma:existence_primal-20}}{\lesssim }
                          k\|\rhsext\|_{-1/2,\Gamma} + \|\rhsm\|_{1/2,\Gamma}.
  \end{align*}
  We now lift these boundary jump functions. Namely, 
  we consider the function $\gimel \in H^1(\R^d \setminus \Gamma)$
  such that
  \begin{align*}
   -\fdiv{(\nu\nabla \gimel)} + (kn)^2 \gimel &= 0\quad\text{and}\quad
   \njmp{\gimel + i k \gimel}=R_1, \;\,\djmp{\gimel}=R_0.
  \end{align*}
  Such a function exists and satisfies (see, e.g., \cite[Prop.~{9}]{LS06}):
  $$
  \|\gimel\|_{1,k,\R^d \setminus \Gamma}
  \lesssim k \big( \|R_1\|_{-1/2,\Gamma} +   \|R_0\|_{1/2,\Gamma}\big)
  \lesssim k^2(1+k^{1/2}) \|\rhsext\|_{-1/2,\Gamma} + k\|\rhsm\|_{1/2,\Gamma}.
  $$
  We consider another smooth cutoff function $\chi_{\widetilde{\Omega}} \in C_{0}^{\infty}(\R^d)$
  with $\chi_{\widetilde{\Omega}} \equiv 1$ in a neighborhood of $\overline{\Omega}$ and 
  $\support{\chi_{\widetilde{\Omega}}}\subseteq \widetilde{\Omega}$ for $\widetilde{\Omega}$
  as in Assumption~\ref{ass:existence_global}.
  Then $\LL-\gimel \chi_{\widetilde{\Omega}} $ solves in $\R^d$
  \begin{align*}
    - \fdiv(\nu \nabla (\LL - \gimel \chi_{\widetilde{\Omega}} ) - (\kn)^2 (\LL-\gimel  \chi_{\widetilde{\Omega}})
    &=
    f-\big[\fdiv(\nu \nabla ( \daleth \,\chi_{\mathcal{N}(\Gamma)})  - (\kn)^2
      \daleth \,  \chi_{\mathcal{N}(\Gamma)} \big] \mathds{1}_{\Omega}\\
     & \qquad  +\big[\fdiv(\nu \nabla (\gimel  \chi_{\widetilde{\Omega}}) ) + (\kn)^2
      \gimel  \chi_{\widetilde{\Omega}}  \big] \\
     &  =: \tilde{f}
  \end{align*}
  with Sommerfeld radiation condition and $\support{\widetilde{f}} \subseteq \widetilde{\Omega}$. 
  For the new right-hand side we get, since $\nu$ is smooth in $\mathcal{N}(\Gamma) \cap \Omega$, the estimate:
  \begin{align*}
   \|\widetilde{f}\|_{0,\widetilde{\Omega}}&\lesssim
                      k \|\daleth\|_{1,k,\Omega} + k\|\gimel\|_{1,k,\widetilde{\Omega}\setminus\Gamma}+ \|f\|_{0,\Omega}
                      \lesssim
                      k^2 \|\rhsext\|_{-1/2,\Gamma}
                      +k(k^{5/2}\|\rhsext\|_{-1/2,\Gamma}+ k\|\rhsm\|_{1/2,\Gamma})
                      + \|f\|_{0,\Omega}.
  \end{align*}
  Assumption~\ref{ass:existence_global} then gives for any ball $B_R$:
  \begin{align*}
    \|\LL - \gimel  \chi_{\widetilde{\Omega}}\|_{1,k,B_R} \lesssim k^{\beta_0}\|\tilde f\|_{0,\Omega} .
  \end{align*}
  This allows us to bound $\LL$ by:
  \begin{align*}
    \|\LL\|_{1,k,B_R} \lesssim k^{\beta_0}\big( \|f\|_{0,\Omega}
    + k^{7/2} \|\rhsext\|_{-1/2,\Gamma} + k^{2} \|\rhsm\|_{1/2,\Gamma} \big) + k^{5/2}\|\rhsext\|_{-1/2,\Gamma} + k \|\rhsm\|_{1/2,\Gamma}.
  \end{align*}
  From this the statement follows by $u=\LL - \daleth \chi_{N({\Gamma})}$ in $\Omega$
  and the previously established estimate (\ref{eq:lemma:existence_primal-20}).  The functions $m$ and $\uext$ can
  be bounded as the jumps of $\daleth$ (cf.\ (\ref{eq:lemma:existence_primal-10})).  The dominant power of $k$ then
  stems from bounding $m$ via Lemma~\ref{lemma:neumann_trace}.
\end{proof}

\subsection{Finite element and boundary element spaces}
\label{sec:fem-bem-spaces}
\subsubsection*{Meshes and element maps}
The finite element spaces on $\Omega$ and $\Gamma$ 
are based on standard regular (i.e., no ``hanging nodes''), $\gamma$-shape regular meshes $\taun$. That is, the (open) elements $K \in \Omega_h$ 
are images of a fixed reference simplex $\widehat K$ under bijective element maps 
$\Phi_K$, they partition $\Omega$, and the element maps of elements sharing 
a $j$-face ($0 \leq j \leq d-1$) induce the same parametrization on that common facet.
A formal definition is given in \cite[Def.~{2.2}]{LiMelenkWohlmuthZou}
or \cite[Sec.~{8.1}]{MS23} for $d=  3$, and \cite[Def.~{2.4.1}]{melenk_book} for $d = 2$. 
In addition, we require the mesh $\taun$ to be compatible with the 
decomposition $\overline{\Omega} = \bigcup_{P \in \PP} \overline{P}$, i.e., 
for each $K \in \taun$ there is a unique $P \in \PP$ with $K \subset \PP$. 
For $K \in \Omega_h$, we write $h_K = \operatorname{diam} K$ and set $h:= \max_{K \in \Omega_h} h_K$. 
Certain approximation results will require analytic element maps $\Phi_K$. It will be 
convenient to make the following assumption. 
\begin{assumption}
[normalizable regular triangulation]\label{def:element-maps} Each element map
$\Phi_{K}$ can be written as $\Phi_{K}=R_{K}\circ A_{K}$, where $A_{K}$ is an
\emph{affine} map, and the maps $R_{K}$ and $A_{K}$ satisfy for constants
$C_{\operatorname*{affine}}$, $C_{\operatorname{metric}}$, $\gamma>0$
independent of $K$:
\begin{align*}
&  \Vert A_{K}^{\prime}\Vert_{L^{\infty}(\widehat{K})}\leq
C_{\operatorname*{affine}}h_{K},\qquad\Vert(A_{K}^{\prime})^{-1}%
\Vert_{L^{\infty}(\widehat{K})}\leq C_{\operatorname*{affine}}h_{K}^{-1},\\
&  \Vert(R_{K}^{\prime})^{-1}\Vert_{L^{\infty}(\widetilde{K})}\leq
C_{\operatorname{metric}},\qquad\Vert\nabla^{n}R_{K}\Vert_{L^{\infty
}(\widetilde{K})}\leq C_{\operatorname{metric}}\gamma^{n}n!\qquad\forall
n\in{\mathbb{N}}_{0}.
\end{align*}
Here, $\widetilde{K}:=A_{K}(\widehat{K})$.  
\eremk
\end{assumption}
It is worth highlighting that Assumption~\ref{def:element-maps} implies the shape 
regularity of the element maps, viz., 
$h_K^{-1}\|\Phi_K \|_{L^\infty(\widehat K)} + h_K \|\Phi_K^{-1} \|_{L^\infty(\widehat K)} 
\leq C$
for some $C > 0$ depending only the constants $C_{\operatorname{affine}}$, $C_{\operatorname{metric}}$, 
$\gamma$. 

By taking traces on $\Gamma$, a mesh $\taun$ on $\Omega$ induces a mesh $\Gamma_h$ on $\Gamma$ 
with element maps induced by the element maps $\Phi_K$ 

\subsubsection*{Approximation spaces}
When working with discontinuous approximation spaces, it is
useful to also work in broken Sobolev spaces.
Given a 
mesh $\taun$ on $\Omega$ and  $s\geq 0$, we define the space
\begin{equation}
\label{eq:broken-space}
H^s_{\mathrm{pw}}(\taun):=\{v\in L^2(\Omega):\, v{}_{|_K}\in H^s(K) \quad \forall K\in \taun\}.
\end{equation}
On the mesh $\taun$ we define the approximation spaces of standard mapped piecewise polynomials:
Denoting by $\mathbb P_\ell(\omega)$ the space of
polynomials of degree at most $\ell$ on the domain $\omega$, we set for $\ell\ge 1$ and $r=0$, $1$
\[
\mathcal S^{\ell, r}(\taun)=\{v\in H^r (\Omega):\,
     v{}_{|_K}\circ \Phi_K\in\mathbb P_\ell(\widehat K)\quad \forall K\in \taun
     \}.
   \]
Analogously, we define on $\Gamma_h$ for $r =0$, $1$
  \[
\mathcal S^{\ell, r}(\Gamma_h)=\{v\in H^r(\Gamma):\,
     v{}_{|_F}\circ \Phi_{K_F}\in\mathbb P_\ell(\widehat F)\quad \forall F\in \Gamma_h
     \},
   \] 
where $K_{F}$ is the element of $\Omega_h$ with $F \in \Gamma_h$ as a $(d-1)$-facet and $\widehat F = \Phi_{K_F}^{-1}(F)$. 

\subsubsection*{DG-related notation}
Let $\mathbf{n}_{K}$ denote the outward pointing unit vector normal to~$\partial K$. 
$(d-1)$-facets (simply called facets) are push-forwards of $(d-1)$-facets of 
the reference simplex $\widehat K$. We collect the (open) internal facets (i.e., the facets lying in $\Omega$) 
and boundary facets (i.e., the facets lying on $\Gamma$) of $\taun$ in the sets~$\FhI$
and $\FhB$, respectively. 
For a mesh $\taun$, we define the mesh size function $\msf:\overline\Omega\to \mathbb R$
by $\msf_{|_K}=h_K$ for any~$K\in \taun$,
$\msf=\min\{h_{K_1},h_{K_2}\}$ on each facet $F\in\FhI$ shared by
$K_1$ and~$K_2$. On each boundary facet $F\in\FhB$, we set
$\msf=h_{K_F}$, where $K_F$ is the element having $F$ as a facet.
We may fix~$\msf$ arbitrarily at mesh vertices and, for $d=3$, on edges, because we will not need it there.
Finally, we recall that the following polynomial inverse inequality from~\cite[Lem.~{A.1}]{MR3667020}:
\begin{align}
  \label{eq:inverse}
  \Vert  \msf^{1/2} p^{-1} \lambda_h \Vert_{0,\Gamma}\leq c_{inv}^G \Vert\lambda_h \Vert_{-\frac{1}{2} ,\Gamma} \qquad \forall \lambda_h\in S^{p, 0}(\Gamma_h),\ p\ge 1.
\end{align}
with a constant $c_{inv}^G$ independent of the mesh size $h$ and the polynomial degree $p$.

\subsection{Conforming FEM}
\label{sec:conforming-fem}
We start with the simpler case of the conforming FEM discretization.
We pick discrete spaces
$$
{\mathcal V}^C_h:= \Vhc \times \Wh \times \Zh, \qquad 
\Vhc :=S^{p,1}(\Omega_h), \qquad \Wh :=S^{p-1,0}(\Gamma_h), \qquad \Zh:=S^{p,1}(\Gamma_h).
$$
We use the superscript $C$ for the $H^1$-\emph{conforming} space $\Vhc$ to
  distinguish it from the \emph{discontinuous Galerkin} space $\VhDG$ defined below.

For the weak form of~\eqref{form:T:for:Helmholtz} we introduce the sesquilinear form
\begin{align}
  \begin{split}
    \TTC{k}(	 \bs u ,  \bs v )&:=
      \TTC{k}(	 (u, m, \uext) , (v,\lambda,\vext) )
  \\&=
 (\Ad\nabla u, \nabla v )_{0,\Omega} -  ((\kn)^2 \,u, v)_{0,\Omega} + \ii k ( u, v)_{0,\Gamma} - \langle m, v \rangle_\Gamma \\
	& \qquad- \langle  (-\Wk- \ii k ( \nicefrac{1}{2} - \Kk  )  + \ii k   ( \nicefrac{1}{2} +\Kprimek + \ii k \Vk ) )\uext 
        -  ( \nicefrac{1}{2} +\Kprimek + \ii k \Vk ) m, \vext   \rangle_\Gamma \\
        &\qquad +  \langle u, \lambda \rangle -  \langle  (\nicefrac{1}{2} + \Kk) \uext - \Vk (m - \ii k \uext), \lambda \rangle_\Gamma. 
\end{split}
 \label{form:T:for:Helmholtz}
\end{align}
The weak form of~\eqref{system:Helmholtz}
for $\rhsm=\rhsext=0$ then reads
\begin{equation} \label{equivalent:weak:formulation}
\begin{cases}
\text{Find } \bs u:=(u, m, \uext) \in {\mathcal V}  = H^1(\Omega) \times H^{-\frac{1}{2}}(\Gamma) \times H^{\frac{1}{2}}(\Gamma)\text{ such that}\\
\TTC{k}( \bs u , \bs v ) = (f, v)_{0,\Omega}
\quad \forall \bs v:=(v,\lambda,\vext) \in {\mathcal V} = H^1(\Omega) \times H^{-\frac{1}{2}}(\Gamma) \times H^{\frac{1}{2}}(\Gamma).
\end{cases}
\end{equation}
The discrete problem is given by taking $\xh:=(\uh,\mh,\uhext)\in {\mathcal V}^C_h$ and  
also restricting the test functions to $\bs v^h:=(\vh,\lambdah,\vhext)\in {\mathcal V}^C_h$. 

For this conforming discretization,
a partial $k$-explicit analysis has already been developed in
\cite[Appendix]{MMPR20}. We present here the most important results and
modify the presentation slightly in order to make it more convenient
for dealing with the discontinuous Galerkin discretization below.

The natural norm for analyzing this problem is given by the energy norm \eqref{eq:enorm}: 
  \begin{align*}
    \smalltriplenorm{}{\bs u}_{k}^2
    &=    \smalltriplenorm{}{(u,m,\uext)}_{k}^2
      := \|\Ad^{\frac12} \nabla u\|^2_{0,\Omega} + \|\kn\, u\|^2_{0,\Omega}+ \|m\|^2_{-\frac{1}{2},\Gamma} + \|\uext\|^2_{\frac{1}{2},\Gamma},
  \end{align*}
 which is equivalent to the previously introduced norm $\smalltriplenorm{}{\bs u}_{k,\mathcal{V},0}$.

The operators given in Lemma~\ref{lemma:introducing_theta} below will play an important role in the analysis of both the
conforming FEM and the discontinuous Galerkin method.  In particular, the operator $\Theta$
captures all the lower order terms which spoil the coercivity
of the  sesquilinear forms $\TTC{k}$ and $\TTDG{k}$; see the G{\aa}rding inequality of Proposition~\ref{prop:conforming} (see~\eqref{eq:TC_after_split}) and of Lemma~\ref{lemma:dg_garding} (see~\eqref{eq:tdg_minus_theta}) below.
For the most part, the terms correspond to the difference between the current sesquilinear
  form and the symmetric coupling of the Laplace problem; see~\eqref{eq:diff_of_blfs}.
In essence they were already presented 
in~\cite[Thm.~{A.2}]{MMPR20}. Here, we just make a slight modification involving the filter operators
in order to get uniform boundedness for the finite-regularity part. In addition,
we add the identity term to $\Theta_{\uext,\vext}$ to get easier coercivity.

\begin{lemma}
  \label{lemma:introducing_theta}
  Let $\Gamma$ be analytic and $k\geq k_0 >0$.
  Define the operators
  \begin{alignat}{6}    
    \Theta_{m,\lambda}&:=(\Vk- \Vz)^*, \qquad &
    \Theta_{\uext,\vext}&:=(\Wk -\Wz)^* + \ii k (\Kk + \Kprimek)^* +
    k^2 \Vk^* +
    1,
    \\
    \Theta_{m,\vext}&:=(\Kprimek -\Kprimez)^* - \ii k \Vk^*  & \qquad
    \Theta_{\uext,\lambda}&:=-(\Kk -\Kz)^* + \ii k \Vk^*, 
  \end{alignat}
  as well as $ \Theta_{u,v} := 2(\kn)^2$.
  Define the operator $$\Theta: \quad L^2(\Omega) \times H^{-1/2}(\Gamma) \times H^{1/2}(\Gamma)
    \to L^2(\Omega) \times H^{1/2}(\Gamma) \times H^{-1/2}(\Gamma)$$ by
  \begin{align*}
    & \langle (u,m,\uext),\Theta (v,\lambda,\vext)\rangle
      =  (  u, \Theta_{u,v} v)_{0,\Omega} \\
    \nonumber
    &\qquad\qquad\qquad
      \mbox{} -
      \langle \uext, \Theta_{\uext,\vext}\vext\rangle_\Gamma -
      \langle m, \Theta_{m,\vext} \vext\rangle_\Gamma -
      \langle \uext, \Theta_{\uext,\lambda} \lambda \rangle_\Gamma - 
      \langle m, \Theta_{m,\lambda}\lambda \rangle_\Gamma.
  \end{align*}
  Then, the operator $\Theta$ can be split into a finite regularity part $\Theta^\F$
  and an analytic remainder $\Theta^\A$ as $\Theta=\Theta^\F +\Theta^\A$.
  For $\sTheta \in [0,1]$, the finite regularity part satisfies the bound
    \begin{multline}
\label{eq:continuity-ThetaF}
    |\langle (u,m,\uext),\Theta^\F (v,\lambda,\vext)\rangle|
    \lesssim
    k^{\sTheta}
    \Big(\|u\|_{1-\sTheta,k,\Omega}+ \|m\|_{-1/2-\sTheta,\Gamma} + \|\uext\|_{1/2-\sTheta,\Gamma}\Big) \\
       \times \Big(\|v\|_{1,k,\Omega}+ \|\lambda\|_{-1/2,\Gamma} + \|\vext\|_{1/2,\Gamma}\Big)
    \end{multline}
and also 
    \begin{align}
\label{eq:continuity-ThetaF-foo}
      \smalltriplenorm{}{\Theta^{\F} (v,\lambda,\vext)}_{\mathcal{V}',1}
      &\lesssim k\smalltriplenorm{}{(v,\lambda,\vext)}_{k,\mathcal{V},0}. 
    \end{align}
%
The analytic part $\Theta^{\A}$ can be further decomposed as 
$\Theta^{\A} = \Theta^\A_{m,\lambda} + \Theta^\A_{m,\vext} + \Theta^\A_{\uext,\lambda} + \Theta^\A_{\uext,\vext}$ with
  \begin{alignat*}{6}
    \Theta^\A_{m,\lambda} \lambda
    &\in \aclass\big(C_{m,\lambda} k \|\lambda\|_{-3/2,\Gamma},\vartheta_{m,\lambda}, \Omega, \Gamma \big), \quad &
    \Theta^\A_{m,\vext} \vext
    &\in \aclass\big(C_{m,\vext} \, k^2 \, \|\vext\|_{-1/2,\Gamma},\vartheta_{m,\vext}, \Omega, \Gamma \big), \\
    \Theta^\A_{\uext,\lambda} \lambda
    &\in \aclass\big(C_{\uext,\lambda} k^{(d+2)/2} \|\lambda\|_{-3/2,\Gamma},\vartheta_{\uext,\lambda}, \Omega, \Gamma \big), \quad & 
    \Theta^\A_{\uext,\vext} \vext
    &\in \aclass\big(C_{\uext,\vext} k^{3} \|\vext\|_{-1/2,\Gamma},\vartheta_{\uext,\vext}, \Omega, \Gamma \big),
  \end{alignat*}
  as well as $\Theta_{u,v}^\A:=0$.
with constants $C_{m,\lambda}$, $C_{m,\vext}$, $C_{\uext,\lambda}$, $C_{\uext,\vext}$
$\vartheta_{m,\lambda}$, $\vartheta_{m,\vext}$, $\vartheta_{\uext,\lambda}$, $\vartheta_{\uext,\vext}$ depending solely on $\Omega$.
\end{lemma}
\begin{proof}
We will not explicitly prove (\ref{eq:continuity-ThetaF-foo}) as it follows from arguments similar to those 
for (\ref{eq:continuity-ThetaF}). 

For (\ref{eq:continuity-ThetaF}), we start with the volume term $\Theta_{u,v}$. Since the analytic part $\Theta_{u,v}^\A$ is zero
we bound for $\sTheta\geq 0$
\begin{align*}
  \big|\big( u, \Theta_{u,v}^{\F} v\big)_{0,\Omega}\big|
  &=\big|\big( u, \Theta_{u,v} v\big)_{0,\Omega}\big|
    = \big|\big( u, (\kn)^2 v\big)_{0,\Omega}\big|
    \lesssim k\|u\|_{0,\Omega} k\|v\|_{0,\Omega}\\
    &\leq k^{\sTheta} \|u\|_{1-\sTheta,k,\Omega} \|v\|_{1,k,\Omega}.
\end{align*}
We move on to the boundary operators.
As in \cite[Lem.~{3.8}]{MMPR20}, we can compute the adjoints of the BEM operators
  in the following way: for an operator $\mathcal{A}$, we define
  $\mathcal{A}'\varphi:=\overline{\mathcal{A} \overline{\varphi}}$ and get for the adjoints
  \begin{subequations}
    \label{eq:some_adjoints}
  \begin{alignat}{6}
    \Vk^* \varphi&=(\Vk)' \varphi, \qquad &
    \Wk^* \varphi&=(\Wk)' \varphi, \\
    (\Kprimek)^*\varphi&= (\Kk)'\varphi, \qquad & (\Kk)^*\varphi&= (\Kprimek)'\varphi.
  \end{alignat}
\end{subequations}
  Thus, also the splittings from Proposition~\ref{prop:decomposition} carry over to the
  adjoints, we just have to add $'$ to the operators and exchange $\Kk$ and $\Kprimek$.

  Overall, using Proposition~\ref{prop:decomposition} and the filter operators from
  Proposition~\ref{prop:ff_gamma1}, we can split the operators
  like
  \begin{align*}
    \Theta_{m,\lambda}&=\Theta^{\F}_{m,\lambda} +  \Theta^{\A}_{m,\lambda} 
                        := \SV' + \gammazint \AVtilde',\\
    \Theta_{m,\vext}&=\Theta^{\F}_{m,\vext} +  \Theta^{\A}_{m,\vext} \\
                      &:= \big(\SK' - \ii k \hffg{+} \Vz - \ii k \SV'  \big)
                        +  \big( \gammazint\AKtilde' - \ii k \lffg{+}
                        \Vz - \ii k  \gammazint \AVtilde'
                        \big),
    \\
    \Theta_{\uext,\lambda}&=\Theta^{\F}_{\uext,\lambda} +  \Theta^{\A}_{\uext,\lambda}  \\
                      &:= \big(-\hffg{-}\SKprime' + \ii k \hffg{-}\Vz + \ii k \SV'\big)
                        + (-\lffg{-}\SKprime'
                        -\gammaoint \AVtilde' + \ii k \lffg{-} \Vz + \ii k  \gammazint \AVtilde'), 
    \\
    \Theta_{\uext,\vext} &=\Theta^{\F}_{\uext,\vext} +  \Theta^{\A}_{\uext,\vext} \\
                      &:= \big(\SW' + \ii k \hffg{-} (\Kprimez + \Kz)
                        + \ii k \SK' + \ii k \SKprime'  + k^2 \hffg{+}\Vz  + 1+k^2 \SV'\big) \\
    &\quad+\big(-\gammaoint \AKtilde' + \ii k \lffg{-} (\Kprimez + \Kz)
      + \ii k \gammazint \AKtilde'+ \ii k \gammaoint \AVtilde'  + k^2 \lffg{+}\Vz +k^2\gammazint\AVtilde'\big),
  \end{align*}

  We prove the stated mapping properties, starting with the
  finite-regularity parts. For $\Theta^{\F}_{m,\lambda}=\SV'$, this
  follows directly from Proposition~\ref{prop:decomposition}, estimate~\eqref{eq:estimates-S-a},
  using $s'=1+\sTheta$ and $s=0$.
  More interesting is the operator $\Theta^\F_{m,\vext}$. We use the mapping properties from
  Proposition~\ref{prop:decomposition}
  (bound~\eqref{eq:estimates-S-c} with $s'=1+\sTheta$ and $s=0$, and bound~\eqref{eq:estimates-S-a} with $s'=1+\sTheta$ and $s=1$) and
  Proposition~\ref{prop:ff_gamma1}
  (bound~\eqref{eq:ff_gamma1_2} with $s'=1/2+\sTheta$ and $s=3/2$):
  \begin{align*}
    \|\Theta^\F_{m,\vext} \vext\|_{1/2+\sTheta,\Gamma}
    &=\|\SK'\vext\|_{1/2+\sTheta,\Gamma}+ k \|\hffg{+} \Vz\vext\|_{1/2+\sTheta,\Gamma} +\| \ii k \SV' \vext\|_{1/2+\sTheta,\Gamma} \\
    &\lesssim k^{\sTheta}\|\vext\|_{1/2,\Gamma} + k^{\sTheta}\|\Vz\vext\|_{3/2,\Gamma} + k^{\sTheta}\|  \vext\|_{1/2,\Gamma} \\
    &\lesssim k^{\sTheta}\|\vext\|_{1/2,\Gamma} + k^{\sTheta}\|\vext\|_{1/2,\Gamma} + k^{\sTheta}\|  \vext\|_{1/2,\Gamma}. 
  \end{align*}
  The proof for $\Theta^\F_{\uext,\lambda}$ is similar,
    but we
    need to split the adjoint double layer operator in order
    to get stability in the $H^{-1/2}(\Gamma)$ norm.
    \begin{align*}
      \|\Theta_{\uext,\lambda}\lambda\|_{-1/2+\sTheta,\Gamma}
      &\le
        \|\hffg{-}\SKprime' \lambda\|_{-1/2+\sTheta,\Gamma} + k \|\hffg{-}\Vz \lambda \|_{-1/2+\sTheta,\Gamma} + k\| \SV' \lambda\|_{-1/2+\sTheta,\Gamma} \\
      &\lesssim
        k^{-1+\sTheta}\|\SKprime' \lambda\|_{1/2,\Gamma} + k^{\sTheta}\|\Vz \lambda \|_{1/2,\Gamma} + k\| \SV' \lambda\|_{-1/2+\sTheta,\Gamma}  \\
      &\lesssim
        k^{\sTheta}\|\lambda\|_{-1/2,\Gamma} + k^{\sTheta}\|\lambda \|_{-1/2,\Gamma} + k^{\sTheta}\|\lambda\|_{-1/2,\Gamma}.
    \end{align*}
    Here, in the first step we used
    Proposition~\ref{prop:ff_gamma2}~\ref{it:hffgn_est}
    with $s'=-1/2+\sTheta$ and $s=1/2$ twice,
    and then in the second step we used
    \eqref{eq:estimates-S-b} with
    $s'=2$ and $s=0$, the mapping properties of $\Vz$ and
    and \eqref{eq:estimates-S-a} with $s'=\sTheta$ and $s=0$.

  Finally, we look at $\Theta_{\uext,\vext}^\F$, and bound the different contributions:
  \begin{align*}
    \|\SW' \vext\|_{-1/2+\sTheta,\Gamma}
    &\stackrel{\eqref{eq:estimates-S-d}}{\lesssim} k^{\sTheta}\|\vext\|_{1/2,\Gamma},\\
    k \|\hffg{-} (\Kprimez + \Kz)\vext\|_{-1/2+\sTheta,\Gamma}
    &\stackrel{\eqref{eq:ff_gamma2_2}}{\lesssim} k^{\sTheta}\big(\|\Kprimez \vext\|_{1/2.\Gamma}
      + \|\Kz \vext\|_{1/2,\Gamma}\big)\lesssim k^{\sTheta}\|\vext\|_{1/2,\Gamma},\\    
    \| \ii k \SK' \vext\|_{-1/2+\sTheta,\Gamma} + \| \ii k \SKprime' \vext\|_{-1/2+\sTheta,\Gamma}
    &\stackrel{\eqref{eq:estimates-S-c}\,\eqref{eq:estimates-S-b}}{\lesssim}
       k^{\sTheta}\| \vext\|_{1/2,\Gamma}, \\
   k^2 \|\hffg{+}\Vz \vext\|_{-1/2+\sTheta,\Gamma}
    &\stackrel{\eqref{eq:ff_gamma1_2}}{\lesssim}
      k^{\sTheta}\|\Vz \vext\|_{3/2,\Gamma}
      \lesssim k^{\sTheta}\|\vext\|_{1/2,\Gamma}, \\
    \|\vext\|_{-1/2+\sTheta,\Gamma} &
                                           \stackrel{\phantom{\eqref{eq:ff_gamma1_2}}}{\lesssim} \|\vext\|_{1/2,\Gamma},\\
    k^2\|\SV' \vext\|_{-1/2+\sTheta,\Gamma}
    &\stackrel{\eqref{eq:estimates-S-a}}
      {\lesssim} k^{\sTheta}\|\vext\|_{1/2,\Gamma}.
  \end{align*}

What is left is to show that the remainder terms are in the right analyticity classes.
  This follows directly from the definitions of the operators and the
  mapping properties in Propositions~\ref{prop:decomposition} and~\ref{prop:ff_gamma1}.
  We note that the estimates are not necessarily sharp and instead we
  crudely used the highest power of $k$ and the strongest norm necessary when
  determining the analyticity classes.
\end{proof}

\begin{proposition}[{\cite[Thm.~{A.2}]{MMPR20}}]\label{prop:conforming}
Let $\Gamma$ be analytic and $k \ge k_0 > 0$.  Let $\Theta$ be defined as in Lemma
\ref{lemma:introducing_theta}. Then the following statements are valid:
\begin{enumerate}[label=(\roman*)]
\item \label{item:thm:k-explicit-garding-ii}
(G{\aa}rding inequality) For a constant $c > 0$ depending only on $k_0$, $c_0$, $\nu_{min}$ and $\Gamma$, 
the sesquilinear form $\TTC{k}(\cdot,\cdot)$ defined in~\eqref{form:T:for:Helmholtz} satisfies 
for all
$\bs v=(v,\lambda,\vext) \in \mathcal{V}$
\begin{align*}
\Re \big( \TTC{k}\bigl( \bs v,\bs v\bigr) 
+ \langle \bs v,\Theta \bs v\rangle 
\big)\
  &\ge c
    \smalltriplenorm{}{\bs u}_k^2.
\end{align*}
\item 
\label{item:thm:k-explicit-garding-iii}
(Continuity) For a constant $C_{cont} > 0$ depending only on $k_0$, $c_0$, $\nu_{min}$ and $\Gamma$, the 
sesquilinear form $\T_C(\cdot,\cdot)$ defined
in~\eqref{form:T:for:Helmholtz} satisfies 
for all $\bs u=(u,m,\uext)$, $\bs v=(v,\lambda,\vext) \in \mathcal{V}$
\begin{align}
\label{eq:item:thm:k-explicit-garding-iii-10}
& \big| \TTC{k}\bigl( \bs u,\bs v\bigr) 
+ \langle \bs u,\Theta \bs v\rangle \big|\leq C_{cont} 
\smalltriplenorm{}{\bs u}_k \smalltriplenorm{}{\bs v}_k, \\
\label{eq:item:thm:k-explicit-garding-iii-20}
&  \big| \TTC{k}\bigl( \bs u,\bs v\bigr) 
+ \langle \bs u,\Theta^\A \bs v\rangle \big|\leq C_{cont} 
\smalltriplenorm{}{\bs u}_k \smalltriplenorm{}{\bs v}_k. 
\end{align}
\end{enumerate}
\end{proposition}
\begin{proof}
  This is just a slight modification of~\cite[Thm.~{A.2}]{MMPR20} with minor
  modifications due to the changed definition of $\Theta$.
  Since it provides a prototype
  of how to proceed for the discontinuous Galerkin method, we still include a sketch of the proof.

  \emph{Ad~\ref{item:thm:k-explicit-garding-ii}:}
  We consider the elliptic sesquilinear form
  \begin{align}
    \begin{split}
      \TT_{+}((u,m,\uext),(v,\lambda,\vext))
    &:=
      (\nu \nabla u, \nabla v)_{0,\Omega} +  ((\kn)^2 u, v)_{0,\Omega}
      + i k (u,v)_{0,\Gamma}
      + \langle  u, \lambda \rangle_\Gamma \\
    &\quad -\langle m, v \rangle_\Gamma
      +\langle \big(\nicefrac{1}{2}+ \Kprimez\big) m, \vext \rangle_\Gamma
      +\langle \Vz m, \lambda\rangle_\Gamma  \\
    &\quad+ \langle (\Wz+1)\uext, \vext \rangle_\Gamma
    - \langle  \big(\nicefrac{1}{2} + \Kz\big) \uext, \lambda \rangle_\Gamma.
  \end{split}
      \label{eq:TC_after_split}
  \end{align}
  Note that, except for the terms $((\kn)^2\,u,v)_{0,\Omega}+i k (u,v)_{0,\Gamma}$
  and $\langle  1\, \uext,\vext\rangle_\Gamma$,
  this sesquilinear form 
  corresponds to the sesquilinear form $\TT_{0}$.
  Taking the difference of $\TTC{k}$ and $\TT_+$, we therefore get:
  \begin{align}
  \TTC{k}(	 &(u, m, \uext) , (v,\lambda,\vext) )
  - \TT_+(	 (u, m, \uext) , (v,\lambda,\vext) )
\label{eq:diff_of_blfs}\\
  &=
  - 2 ((\kn)^2 \,u, v)_{0,\Omega}  \nonumber \\
	&\; -\langle  (-(\Wk-\Wz -1 )- \ii k ( \nicefrac{1}{2} - \Kk  )  + \ii k   ( \nicefrac{1}{2} +\Kprimek + \ii k \Vk ) )\uext 
	 -  ( (\Kprimek-\Kprimez) + \ii k \Vk ) m, \vext   \rangle_\Gamma \nonumber \\     
  &\;-  \langle  (\Kk-\Kz) \uext
  - (\Vk-\Vz) m + \ii k \Vk \uext), \lambda \rangle_\Gamma \nonumber\\
  &= \langle (u,m,\uext), \Theta (v,\lambda,\vext)\rangle_\Gamma \nonumber
  \end{align}
  where in the last step we inserted the definition of $\Theta$ from Lemma~\ref{lemma:introducing_theta}. 
 Taking $\bs u:=(v,\lambda,\vext)=(u,m,\uext)$ and only considering the real part, we get:
  \begin{align*}
    \Re\big(\TTC{k}( \bs u,\bs u)    
    +\langle \bs u, \Theta \bs u\big)
      &= \Re\big(\TT_{+}( (u, m, \uext),  (u, m, \uext))\big)\\
    &=
      \|\nu^{1/2} \nabla u\|_{0,\Omega}^2 + \|\kn \,u\|^2_{0,\Omega} 
     + \langle  (\Wz+1) \uext, \uext \rangle + \langle \Vz m , m \rangle_\Gamma,
  \end{align*}
  where we used that $- \langle m, u \rangle_\Gamma + \langle u, m \rangle_\Gamma$
  and $\langle ( \frac{1}{2} +\Kprimez  ) m, \uext   \rangle_\Gamma    -  \langle  (\frac{1}{2} + \Kz) \uext ,m \rangle_\Gamma$  are purely imaginary, since $\Kz^* = \Kz'$.
  The statement then follows because  $\Wz+1$ and $\Vz$ are 
  coercive with respect to the $H^{1/2}(\Gamma)$ norm and the
    $H^{-1/2}(\Gamma)$ norm, respectively.

  \emph{Ad~\ref{item:thm:k-explicit-garding-iii}:}
  The continuity follows directly from the explicit form of $\TTC{k}(\cdot,\cdot) +\langle\cdot, \Theta\,\cdot\rangle=\T_{+}$  
  and the boundedness of all the boundary operators in~\eqref{eq:TC_after_split}. We only remark that
  we use the multiplicative trace estimate and Young's inequality to estimate
  boundary terms of the form
  \begin{align*}
  k \|u\|^2_{0,\Gamma} &\lesssim
  k^2 \| u\|_{0,\Omega}^2 + \|u\|_{1,\Omega}^2
  =(1+k^2)\| u\|_{0,\Omega}^2+\|\nabla u\|_{0,\Omega}^2\\
   & \leq  (k_0^{-2}+1)n_{min}^{-2} \|\kn\, u\|_{0,\Omega}^2 
   + \nu_{min}^{-1} \|\Ad^{\frac12}\nabla u\|_{0,\Omega}^2
   \lesssim\|\kn\, u\|_{0,\Omega}^2 +\|\Ad^{\frac12}\nabla u\|_{0,\Omega}^2. 
  \end{align*}
The continuity estimate (\ref{eq:item:thm:k-explicit-garding-iii-20}) follows from 
(\ref{eq:item:thm:k-explicit-garding-iii-10}) and a triangle inequality in view of 
the uniform-in-$k$ boundedness of $\Theta^\F$ asserted in (\ref{eq:continuity-ThetaF}).
\end{proof}

Finally, we have a boundedness result for the full sesquilinear form with 
polynomial growth in $k$ of the continuity constant.
\begin{corollary}
  \label{cor:TDC_is_bounded_with_k}
  Let $k \geq k_0 >0$. Then, there exists a constant $C>0$, possibly depending on
  $\Gamma$, $k_0$, $\nu_{min}$,  
  and there exists a constant $\mu_{stab} \in [0,4]$
  such that for all $\bs u=(u,m,\uext)$, $\bs v=(v,\lambda,\vext) \in \mathcal{V} = 
H^1(\Omega) \times H^{-1/2}(\Gamma) \times H^{1/2}(\Gamma)$
  \begin{align*}
    |\T_{C} ( \bs u , \bs v )|
    &\leq C k^{\mu_{stab}}\smalltriplenorm{}{\bs u}_k\smalltriplenorm{}{\bs v}_k.
  \end{align*}
  \end{corollary}
  \begin{proof}
    This bound follows from Proposition~\ref{prop:conforming} and
    the mapping properties of $\Theta$ as spelled out in Lemma~\ref{lemma:introducing_theta}.
    The powers $k^{\mu_{stab}}$ originate from the analytic remainder terms, which are
    controlled using standard trace estimates:
      \begin{align*}
        |\TTC{k} ( \bs u , \bs v )|
        &\leq |\TTC{k} ( \bs u , \bs v )
        +  \langle \bs u , \Theta \bs v \rangle | 
        + |\langle \bs u , \Theta^{\F} \bs v \rangle|
        + |\langle \bs v , \Theta^{\A} \bs u \rangle|.
      \end{align*}
      The first term is uniformly bounded in $k$ by
      Proposition~\ref{prop:conforming}. The
      term involving $\Theta^{\F}$ is bounded by Lemma~\ref{lemma:introducing_theta}.
      To bound $\Theta^{\A}$, we look at the terms individually using the mapping properties
      from Lemma~\ref{lemma:introducing_theta}. For simplicity, we only
      treat the dominant term in detail. The others are bounded
      analogously:
      \begin{align*}
        \|\Theta_{\uext,\vext} \vext\|_{-1/2,\Gamma}
        &\lesssim
          \|\Theta_{\uext,\vext} \vext\|_{0,\Gamma}
          \lesssim \| \Theta_{\uext,\vext} \vext \|_{1,\Omega} \\
        &\lesssim  k^3 C_{\uext,\vext} \|\vext\|_{-1/2,\Gamma} (1+ k \vartheta_{\uext,\vext}). \qedhere
      \end{align*}      
  \end{proof}

\subsection{Discontinuous Galerkin method}\label{sec:DG}
In this section, we discretize problem~\eqref{eq:model_problem} using a
discontinuous Galerkin (DG) method for the interior subproblem.
For a detailed derivation of the method and numerical experiments, we refer to~\cite{dgfembem}.

We start by defining the penalty 
functions $\alpha$, $\beta$, $\delta$. 
Denoting by $\nu_K$ the analytic extension of $\nu_{|_K}$ up to $\partial K$, we define the function $\widetilde{\nu}: \FhI
\to \mathbb{R}$ 
by $\widetilde{\nu}(\bx)=\max\{\abs{\nu_{K_1}(\bx)},\abs{\nu_{K_2}(\bx)}\}$ if $\bx$ is on a facet $F\in\FhI$ shared by $K_1$ and $K_2$.
We recall that, for matrices, we denote by $\abs{\cdot}$ the spectral norm.
Similarly to~\cite{dgfembem},
we set 
\begin{equation} \label{dG-parameters}
\begin{split}
\alpha(\bx) := \afrak \,
&\frac{p^2 \widetilde{\nu}(\bx)}{k \msf(\bx)}, \qquad \beta(\bx) := \bfrak \,\frac{k\msf(\bx)}{p\,\widetilde{\nu}(\bx)} \qquad  \forall \bx \in \FhI, \quad
\delta(\bx) :=\dfrak \frac{k\msf(\bx)}{p^2} \quad  \forall \bx \in \FhB,
\end{split}
\end{equation}
with chosen constants $\afrak>0, \bfrak\ge 0, \dfrak >0$.  Throughout we assume that
$\afrak$ is sufficiently large and $\dfrak$ is sufficiently small.
We also require $\delta(\bx) \in (0,1/2)$.
Next, we need some notation for the jump and the average functionals on~$\FhI$
for piecewise smooth, scalar functions~$v$ and vector-valued functions~${\boldsymbol{\tau}}$.
At any~$\bx \in \FhI$ shared by the two elements $K_{\bx}^1$
and~$K_{\bx}^2$, 
the jumps $\llbracket v
\rrbracket$ and $\llbracket {\boldsymbol{\tau}}\rrbracket$, 
and the averages $\ldc  v   \rdc$ and $\ldc  {\boldsymbol{\tau}}\rdc$
are defined as 
\[
\begin{split}
& \llbracket v \rrbracket (\bx) := v{}_{|_{K_{\bx}^1}} (\bx) \,\bn{}_{K_{\bx}^1}  +  v{}_{|_{K_{\bx}^2}}(\bx) \, \bn_{K_{\bx}^2}, \quad \quad \ldc  v   \rdc (\bx) = \frac{1}{2} \left(v{}_{|_{K_{\bx}^1}} (\bx) + v{}_{|_{K_{\bx}^2}} (\bx) \right), \\
& \llbracket {\boldsymbol{\tau}} \rrbracket (\bx)  := {\boldsymbol{\tau}}{}_{|_{K_{\bx}^1}} (\bx)  \cdot \bn{}_{K_{\bx}^1} + {\boldsymbol{\tau}}{}_{|_{K_{\bx}^2}} (\bx) \cdot \bn{}_{K_{\bx}^2}, \quad \quad  \ldc  {\boldsymbol{\tau}} \rdc (\bx)  :=\frac{1}{2} \left({\boldsymbol{\tau}}{}_{|_{K_{\bx}^1}} (\bx) + {\boldsymbol{\tau}}{}_{|_{K_{\bx}^2}}(\bx) \right).
\end{split}
\]

The discrete spaces are chosen as $\VhDG:=S^{p,0}(\Omega_h)$, $\Wh=S^{p-1,0}(\Gamma_h)$, and
$\Zh:=S^{p,1}(\Gamma_h)$ and we set 
$$
\mathcal{V}^{DG}_h:= \VhDG \times \Wh \times \Zh. 
$$
We define the interior DG sesquilinear form by
\begin{align*}
  \ahOmega(\uh, \vh)
  &:= \sum_{K \in \taun}{\!\!\big(\int_K \nu\,\nabla \uh \cdot \overline{\nabla \vh} - \int_K (\kn)^2\uh \overline{\vh}\big)}
 -\int_{\FhI}\left(\llbracket \uh \rrbracket \cdot \ldc \overline{\nu\nabla_h \vh} \rdc +\ldc \nu\nabla_h \uh \rdc \cdot \llbracket \overline{\vh} \rrbracket \right) \\
  &-\int_{\FhI}\beta (ik)^{-1} \llbracket \nu\nabla_h \uh \rrbracket  \llbracket \overline{\nu\nabla_h \vh} \rrbracket
            +\int_{\FhI} \alpha\, ik \llbracket  \uh \rrbracket \cdot \llbracket \overline{\vh}\rrbracket,
\end{align*}
and the DG boundary sesquilinear form~$\bhGamma (\cdot, \cdot)$ by 
\[
 \begin{split}
 \bhGamma (\uh, \vh)
 & := - \int_\Gamma \delta (\ii k)^{-1} \nu\nabla_h \uh \cdot \nGamma \,\overline{\nu\nabla_h \vh \cdot \nGamma} 
 - \int_\Gamma \delta \uh \,\overline {\nu\nabla_h \vh \cdot \nGamma}\\
 & \quad \quad - \int_\Gamma \delta \,\nu\nabla_h \uh \cdot \nGamma \overline{\vh} + \int_{\Gamma} (1-\delta)\ii k  \uh \overline{\vh}.
 \end{split}
 \]

The sesquilinear form corresponding to the DG discretization of~\eqref{system:Helmholtz} is then given by 
\begin{equation} \label{trisesquilinear-form}
\begin{split}
& \TTDG{k} ( (\uh,\mh,\uhext) , (\vh,\lambdah, \vhext) )  \\
& := \ahOmega(\uh,\vh) + \bhGamma(\uh,\vh) - (\mh, \delta (\i k)^{-1} \nu\nabla_h \vh \cdot \nGamma + (1-\delta)\vh)_{0,\Gamma} \\
& \quad \quad - \langle (-\Wk -\i k (\nicefrac{1}{2}-\Kk) + \i k (\nicefrac{1}{2}+\Kprimek+\i k\Vk))\uhext - (\nicefrac{1}{2} + \Kprimek + \i k \Vk)\mh, \vhext \rangle_\Gamma\\
& \quad \quad + \langle -\delta (\i k)^{-1} \nu\nabla_h \uh \cdot \nGamma + (1 - \delta) \uh + \delta (\i k)^{-1} \mh, \lambdah \rangle_\Gamma\\
& \quad \quad - \langle (\nicefrac{1}{2} +\Kk) \uhext - \Vk (\mh - \i k \uhext) , \lambdah   \rangle_\Gamma.
\end{split}
\end{equation}
The method can then be written in compact form as follows:
\begin{equation} \label{dgBEM:short-version}
\begin{cases}
\text{Find } \bs u^h:=(\uh, \mh,\uhext) \in {\mathcal V}^{DG}_h  = \VhDG \times \Wh \times \Zh \text{ such that}\\
\TTDG{k}\left( \bs u^h , \bs v^h \right) = (f,\vh)_{0,\Omega} \quad \forall \bs v^h:=(\vh, \lambdah, \vhext) \in {\mathcal V}^{DG}_h.  
\end{cases}
\end{equation}

The natural norms of the DG method for the interior subproblem,
for~$v\in H_{\mathrm{pw}}^{\frac{3}{2} + t} (\Omega_h)$, with~$t>0$ arbitrary,  
are given by
\begin{equation} \label{dG:norm:1}
\begin{split}
\Vert v \Vert ^2_{\dGnorm}	
& := \Vert \nu^{1/2}\,\nabla_h v\Vert^2_{0,\Omega} +  \Vert \kn\, v \Vert^2_{0,\Omega} + k^{-1} \Vert \beta^{1/2} \llbracket \nu\nabla_h v\rrbracket \Vert^2_{0,\FhI}  + k \Vert \alpha^{1/2} \llbracket v \rrbracket \Vert^2_{0,\FhI}\\
& \quad + k^{-1} \Vert \delta ^{1/2} \nu\nabla_h v \cdot \nGamma \Vert^2_{0,\Gamma} + k \Vert (1-\delta) ^{1/2} v \Vert^2 _{0,\Gamma}
\end{split}
\end{equation}
and 
\begin{equation}
\label{dG:norm:foo}
\begin{split}
\Vert v \Vert ^2_{\dGnormp} := \Vert v \Vert ^2_{\dGnorm} + k^{-1} \Vert \alpha^{-1/2} \ldc \nu\nabla_h v  \rdc \Vert^2_{0,\FhI}.
\end{split}
\end{equation}

We further introduce the  following two energy norms, which extend the $\dGnorm$ and $\dGnormp$  norms to the FEM-BEM coupling:
\begin{align*}
\energynorm{(u,m,\uext)}^2   &:= \|{u}\|_{\dGnorm}^2 + \|m\|_{- 1/2,\Gamma}^2 + \|\uext\|_{1/2,\Gamma}^2, \\
\energynormp{(u,m,\uext)}^2 &:= \|{u}\|_{\dGnormp}^2 + \|m\|_{-1/2,\Gamma}^2 + \| \msf^{1/2}p^{-1} \ m\Vert_{0,\Gamma}^2 + \|\uext\|_{1/2,\Gamma}^2.
\end{align*}

The main ingredients in the  convergence proof of \cite{dgfembem}
are a G{\aa}rding  inequality and the boundedness of the
sesquilinear form. We now transfer these estimates to the $k$-explicit
setting.

\begin{lemma}[$k$-explicit G{\aa}rding inequality]
  \label{lemma:dg_garding}
There are constants $\mathfrak{a}_0$, $\mathfrak{d}_0$, $\varepsilon > 0$ independent of $k \ge k_0$ such that, for 
$\mathfrak{a} \geq \mathfrak{a}_0$, 
$0 < \mathfrak{d} \leq \mathfrak{d}_0$ we have  
  for all $\bs u^h:=(\uh,\mh,\uhext) \in \mathcal{V}^{DG}_h$
  \begin{align*}
    \smalltriplenorm{}{\bs u^h}_{\dGnorm}^2
    \lesssim      
        (\Re+\varepsilon \Im) \Big[
      \TT_{DG}\big(\bs u^h,\bs u^h\big) 
      + \big\langle\bs u^h,\Theta \bs u^h\big\rangle \Big].
  \end{align*}
  
\end{lemma}
\begin{proof}
  Just as in the proof of Proposition~\ref{prop:conforming} for the conforming case, we first derive an explicit representation of
  $\TTDG{k} + \Theta$. Most of the boundary terms are the same
  as in the continuous case, see~\eqref{eq:TC_after_split}. We get 
  \begin{align}
    \label{eq:tdg_minus_theta}
    \begin{split}
      & \TTDG{k} ( (\uh,\mh,\uhext) , (\vh,\lambdah, \vhext) )
      +  \langle (\uh,\mh,\uhext) , \Theta (\vh,\lambdah, \vhext) \rangle )\\
      & =\ahOmega (\uh,\vh) + \bhGamma(\uh,\vh)
      + 2 ((\kn)^2\uh,\vh)_{0,\Omega}\\
      & \quad \quad + \langle (\Wz+1) \uhext + (\nicefrac{1}{2} + \Kprimez )\mh, \vhext \rangle_\Gamma
      + \langle -(\nicefrac{1}{2} +\Kz) \uhext + \Vz \mh  , \lambdah   \rangle_\Gamma \\
      & \qquad - (\mh, \delta (\i k)^{-1} \nu\nabla_h \vh \cdot \nGamma + (1-\delta)\vh)_{0,\Gamma}
      \\
& \qquad \qquad \quad + \langle -\delta (\i k)^{-1} \nu\nabla_h \uh \cdot \nGamma + (1 - \delta) \uh + \delta (\i k)^{-1} \mh, \lambdah \rangle_\Gamma.
\end{split}
  \end{align}

  We proceed similarly to \cite{dgfembem}.
  Selecting $(\vh,\lambdah, \vhext):=(\uh,\mh,\uhext)$ and taking $\Re+\varepsilon \Im$
  give after some minor reordering to be more consistent with \cite{dgfembem},
  \begin{align*}
    \begin{split}
      & (\Re+\varepsilon \Im)\Big(\TTDG{k} ( (\uh,\mh,\uhext) , (\vh,\lambdah, \vhext) )
      +  \langle (\uh,\mh,\uhext) , \Theta (\vh,\lambdah, \vhext) \rangle_\Gamma \Big) \\
      & =(\Re+\varepsilon \Im)\big[\ahOmega (\uh,\uh) + \bhGamma(\uh,\uh)\big]
      + 2 \|\kn\, \uh\|^2_{0,\Omega}
      \\ & \qquad  \qquad
      + \langle \Vz \mh,\mh\rangle  + \langle (\Wz+1) \uhext,\uhext \rangle_\Gamma
      + 2\varepsilon \Im \langle(\nicefrac{1}{2} + \Kz )\uhext, \mh \rangle_\Gamma \\
      & \qquad \qquad- \varepsilon   k^{-1} \|\delta^{1/2} \mh\|_{0,\Gamma}^2      
      -2 \Re\langle \mh, \delta (\ii k)^{-1} \nu \nabla_h \uh \cdot \nGamma\rangle_{\Gamma}                  
      +2 \varepsilon \Im\langle(1 - \delta) \uh, \mh\rangle_\Gamma     
    \end{split} \\
      &=: (\Re+\varepsilon \Im)\big[\ahOmega (\uh,\uh) + \bhGamma(\uh,\uh)\big]
        + 2 \|\kn\, \uh\|_{0,\Omega}^2
        +\sum_{i=1}^{6}{T_i}.
    \end{align*}
    We estimate each of these terms individually.    
    \begin{subequations}
      Starting with the standard DG terms, we can use
      \cite[Sec.~{4}, Prop.~{2}]{dgfembem} to get
      for $\mathfrak{a}$ sufficiently large and $\mathfrak{d}$ sufficiently
        small:
      \begin{multline*}
        (\Re + \varepsilon \Im) \big[\ahOmega(\uh, \uh)  + \bhGamma(\uh, \uh)\big]
        \geq 
        \frac{1}{2}\Vert \nu^{1/2}\nabla_h \uh \Vert^2_{0,\Omega} -  \Vert \kn\, \uh \Vert^2_{0,\Omega}  \\
        + \frac{1}{2}\varepsilon \left( k^{-1} \Vert \beta^{1/2} \llbracket \nu\nabla_h\uh \rrbracket \Vert^2_{0,\FhI}    + k \Vert \alpha^{1/2} 
        \llbracket \uh \rrbracket \Vert^2_{0,\FhI}  \right. \\
        \left. + k^{-1} \Vert \delta^{1/2} \nu\nabla_h \uh \cdot \nGamma \Vert^2_{0,\Gamma} +  k \Vert (1-\delta)^{1/2} \uh\Vert^2_{0,\Gamma}  \right).
      \end{multline*}
      The only negative term is $-  \Vert \kn\, \vh \Vert^2_{0,\Omega}$. We can thus
      focus on the BEM contributions $T_i$.
      For $T_1$ and $T_2$, we have the standard coercivity:
      \begin{align}
        \langle \Vz \mh,\mh\rangle  + \langle (\Wz+1) \uhext,\uhext \rangle
        \geq c_{\Vz} \|\mh\|^2_{-1/2,\Gamma} + c_{\Wz} \|\uhext\|^2_{1/2,\Gamma}.
      \end{align}
      For $T_3$, we use Young's inequality to get
      \begin{align}
        T_3=2\varepsilon \Im \langle(\nicefrac{1}{2} + \Kz )\uhext, \mh \rangle
        \geq - \varepsilon (1+2C_{\Kz}) \big(\|\uhext\|_{1/2,\Gamma}^2 + \|\mh\|^2_{-1/2,\Gamma}\big).
      \end{align}
      For $T_4$, we use the  polynomial inverse inequality~\eqref{eq:inverse} to get,
      after inserting the definition of $\delta$ from~\eqref{dG-parameters},
      \begin{align}
        T_4=- \varepsilon   k^{-1} \|\delta^{1/2} \mh\|_{0,\Gamma}^2
        \geq - \varepsilon \mathfrak{d} c_{inv}^G \|\mh\|_{-1/2,\Gamma}^2.
      \end{align}      
      Similarly,  $T_5$ is dealt with again via an inverse estimate: 
      \begin{align*}
          T_5&=-2\Re \langle \mh, \delta (\i k)^{-1} \nu\nabla_h \uh \cdot \nGamma \rangle  \\
     &\geq 
       -\frac{4 c_{inv}^G \mathfrak{d}}{\varepsilon} \|\mh\|_{-1/2,\Gamma}^2
       - \frac{\varepsilon  k^{-1}}{4} \|\delta^{1/2} \nu\nabla_h u_h \cdot \nGamma\|_{0,\Gamma}^2.
      \end{align*}

     The term $T_6$ does not depend on any boundary integral operator and thus
      can be estimated just as in \cite[Eqn.~(4.28)]{dgfembem}.
      For a sufficiently large constant $c_{10}$ depending on $\nu$,
      the inverse estimate constant $c_{inv}^G$ from~\eqref{eq:inverse},
        and the shape regularity constant of~$\taun$
    we have 
      \begin{align}
        T_{6}
        &= -2\varepsilon  \Im  \langle (1-\delta)\mh, \uh \rangle \nonumber \\
        &\geq
          -\frac{4 c_{10}\varepsilon }{\nu_{0}} \Vert \mh \Vert_{-\frac{1}{2}, \Gamma}^2 \\
        & \quad  -\frac{ \varepsilon}{4} \Big(\Vert \nu^{1/2}\nabla_h \uh \Vert_{0,\Omega}^2
          +  \nu_{0}(k_0c_0)^{-2} \Vert \kn\, \uh\Vert_{0,\Omega} ^2+ k \Vert  \alpha^{1/2} \llbracket \uh \rrbracket \Vert_{0,\FhI}^2 \Big).\nonumber
      \end{align}
    \end{subequations}
    
    We put everything together and choose $\varepsilon \in (0,1)$ appropriately.
    We end up with
    \begin{align}
      \begin{split}
      (\Re+\varepsilon \Im)
      &\Big(\TTDG{k} ( (\uh,\mh,\uhext) , (\vh,\lambdah, \vhext) )
      +  \langle (\uh,\mh,\uhext) , \Theta (\vh,\lambdah, \vhext) \rangle \Big) \\
      & \geq
        \frac{1}{4}\Vert \nu^{1/2}\nabla_h \uh \Vert^2_{0,\Omega} + \Big(1-\frac{\varepsilon \nu_{0}(c_0 k_0)^{-2}}{4} \Big) \Vert \kn\, \uh \Vert^2_{0,\Omega}  \\
        &
        \qquad\begin{multlined}[t][12cm]
          + \frac{1}{4}\varepsilon \Big(2 k^{-1} \Vert \beta^{1/2} \llbracket \nu\nabla_h\uh \rrbracket \Vert^2_{0,\FhI}    + k \Vert \alpha^{1/2} \llbracket \uh \rrbracket \Vert^2_{0,\FhI}   \\
          + k^{-1} \Vert \delta^{1/2} \nu\nabla_h \uh \cdot \nGamma \Vert^2_{0,\Gamma} + 2 k \Vert (1-\delta) \uh\Vert^2_{0,\Gamma}  \Big)
        \end{multlined}
        \\
      & \qquad + \Big(c_{\Vz} - \varepsilon(1+2C_{\Kz}) - \varepsilon \mathfrak{d} c_{inv}^{G} -
        \frac{4 c_{inv}^G \mathfrak{d}}{\varepsilon} - \frac{4\varepsilon c_{10}}{\nu_{0} }\Big) \|\mh\|_{-1/2,\Gamma}^2 \\
        & \qquad + \big(c_{\Wz} - \varepsilon (1+2C_{\Kz})\big)\|\uhext\|^2_{1/2,\Gamma}.
      \end{split}
          \label{eq:garding_int1}
    \end{align}
    We pick $\varepsilon$ and $\mathfrak{d}$ as
    \begin{align}
      \varepsilon&:=\min\Big(
                   \frac{c_{\Vz}}{4\big(1+2C_{\Kz} + c_{inv}^{G} + \nu_{0}^{-1} 4c_{10}\big)},
                   \frac{c_{\Wz}}{2(1+C_{\Kz})},
                   2(c_0 k_0)^2 \nu_{0}^{-1},
                   1
                   \Big),  \;\;
      \mathfrak{d}\leq \frac{c_{\Vz} \varepsilon}{16 c_{inv}^{G}}.
    \end{align}
    This ensures that all the terms on the right-hand side of \eqref{eq:garding_int1} are
    positive and the stated result follows.    
  \end{proof}

  \begin{lemma}[$k$-explicit continuity bound]
    \label{lemma:dg_explicit_bound}
    Let $\Gamma$ be analytic and $k \geq k_0 >0$. Then, up to compact terms,
    the DG-sesquilinear form $\TTDG{k}$ is bounded uniformly in $k$. Namely, there exists $C>0$ depending
    only on $\Gamma$, $k_0$, $\nu_{min}$, and the shape-regularity of $\TT_h$,
    such that for all
      $\bs u$, $\bs v \in  \mathcal{V}_{\mathrm{pw}}:= H^{3/2+\varepsilon}_{pw}(\Omega_h) \times H^{-1/2}(\Gamma) \times H^{1/2}(\Gamma)$
      \begin{subequations}
        \label{eq:k_expl_continuities}
    \begin{align}
      |\TTDG{k} (\bs u , \bs v )
      +  \langle \bs u , \Theta^{\A} \bs v \rangle | 
      &\le C \energynormp{\bs u}\energynormp{\bs v}. 
        \label{eq:k_expl_continuities:1} 
  \end{align}
    If $\bs u$ or $\bs v$ are in the discrete space
    $\mathcal{V}^{DG}_h$ then the corresponding $\energynormp{\cdot}$ norm can
    be weakened to the $\energynorm{\cdot}$ norm:   
    \begin{align}
      |\TTDG{k} (\bs u^h , \bs v )
      +  \langle \bs u^h , \Theta^{\A} \bs v \rangle | 
      &\le C \smalltriplenorm{}{\bs u^h}_{\dGnorm}\energynormp{\bs v}
        \qquad \forall \bs u^h \in \mathcal{V}^{DG}_h
        ,\bs v \in \mathcal{V}_{\mathrm{pw}},
        \label{eq:k_expl_continuities:2}
      \\
      |\TTDG{k} (\bs u , \bs v^h )
      +  \langle \bs u , \Theta^{\A} \bs v^h \rangle | 
      &\le C \energynormp{\bs u}\smalltriplenorm{}{\bs v^h}_{\dGnorm}
        \qquad \forall \bs u \in \mathcal{V}_{\mathrm{pw}},
        \bs v^h \in \mathcal{V}^{DG}_h.
        \label{eq:k_expl_continuities:3}
    \end{align}
  \end{subequations}
The three continuity bounds in (\ref{eq:k_expl_continuities}) are also valid with the operator $\Theta^\A$ replaced by $\Theta$ in view of 
the continuity assertion (\ref{eq:continuity-ThetaF}) for $\Theta^\F$.
  \end{lemma}
  \begin{proof}
    We write $\bs u=:(u,m,\uext)$ and $\bs v=:(v,\lambda,\vext)$ for the different components
    and focus on the non-discrete case~\eqref{eq:k_expl_continuities:1}.
    We again go back to the explicit representation of $\TTDG{k}+ \Theta$ in
    \eqref{eq:tdg_minus_theta}: 
      \begin{align}    
        \begin{split}
          & \TTDG{k} ( (u,m,\uext) , (v,\lambda, \vext) )
          +  \langle (u,m,\uext) , \Theta (v,\lambda, \vext) \rangle )\\
          & = \ahOmega (u,v) + \bhGamma(u,v)
          +  2((\kn)^2 u,v)_{0,\Omega}\\
          & \quad \quad + \langle (\Wz+1) \uext + (\nicefrac{1}{2} + \Kprimez )m, \vext \rangle_\Gamma
          + \langle -(\nicefrac{1}{2} +\Kz) \uext + \Vz m  , \lambda   \rangle_\Gamma \\
          & \qquad - \langle m, \delta (\i k)^{-1} \nu\nabla_h v \cdot \nGamma + (1-\delta)v\rangle_\Gamma
          \\
          & \qquad \qquad \quad + \langle -\delta (\i k)^{-1} \nu\nabla_h u \cdot \nGamma + (1 - \delta) u + \delta (\i k)^{-1} m, \lambda \rangle_\Gamma.
        \end{split}
      \end{align}
      Most of the terms can easily estimated because they only contain boundary integral operators
      for the Laplacian.
      We focus on the few terms whose estimates are not obvious: 
      \begin{align*}
        T_1&:=\ahOmega (u,v) + \bhGamma(u,v), \\
        T_2&:=- \langle m, \delta (\i k)^{-1} \nu\nabla_h v \cdot \nGamma + (1-\delta) v \rangle_\Gamma, \\ 
        T_3&:= \langle -\delta (\i k)^{-1} \nu\nabla_h u \cdot \nGamma + (1 - \delta) u,\lambda \rangle_\Gamma, \\
        T_4&:= \langle\delta (\i k)^{-1} m, \lambda \rangle_\Gamma.
      \end{align*}
      Term $T_1$ is a standard DG sesquilinear form and an upper bound in the energy
      norm is derived in~\cite[Prop.~{3.1}]{MPS13}.
      The terms $T_2,\, T_3,\, T_4$ are the same coupling terms that were analyzed
      in \cite[Prop.~{5.1}]{dgfembem}. We only point out that all estimates
      are indeed  $k$-independent and only involve standard inverse estimates
      and the use of a reconstruction operator to split functions
       into a $H^{1/2}$-conforming contribution and smaller remainder.
    \end{proof}

    Finally, we have a boundedness result for the full sesquilinear form,
    as long as we allow for polynomial growth in $k$ of the continuity constant.
    \begin{corollary}
      \label{cor:TDG_is_bounded_with_k}
      Let $k \geq k_0 >0$. Then, there exists a constant $C>0$, possibly depending on
      $\Gamma$, $k_0$, $\nu_{min}$, and the shape-regularity of $\TT_h$,
      and there exists a constant $\mu_{stab} \in [0,4]$ 
      such that for $u$, $v \in H^{3/2+\varepsilon}_{pw}(\Omega)$ and
      $m,\lambda  \in H^{-1/2}(\Gamma)$, $\uext, \vext \in H^{1/2}(\Gamma)$:
    \begin{align*}
      |\TTDG{k} ( \bs u , \bs v )|
      &\leq C k^{\mu_{stab}}\energynormp{\bs u}\energynormp{\bs v}.
    \end{align*}
    If $(u,m,\uext)$ or $(v,\lambda, \vext)$ are in the discrete space
    $\mathcal{V}^{DG}_h$ then the corresponding $\energynormp{\cdot}$ norm can
    be weakened to the $\energynorm{\cdot}$ norm.
  \end{corollary}
  \begin{proof}
    Follows completely analogously to Lemma~\ref{cor:TDC_is_bounded_with_k},
      just replacing the result from Proposition~\ref{prop:conforming} with
      Lemma~\ref{lemma:dg_explicit_bound}.
  \end{proof}
  
\section{The adjoint problem}
In this section, we consider the following problem, which corresponds to the adjoint
of~\eqref{system:Helmholtz}.
This problem was already analyzed in \cite{MMPR20}, but without
working out $k$-explicit estimates.

Given $r \in \big(H^{1}(\Omega)\big)'$, $\Rm \in H^{1/2}(\Gamma)$, and $\Rext \in H^{-1/2}(\Gamma)$, find $\psi \in H^1(\Omega)$, $\psim \in H^{-1/2}(\Gamma)$ and $\psiext \in H^{1/2}(\Gamma)$
such that
  \begin{subequations}
    \label{eq:dual_problem}
\begin{align}
&\begin{cases}
-\fdiv(\Ad \nabla\overline  \psi) - (\kn)^2 \overline \psi =  \overline r \qquad\text{in }\Omega, \\
\njmp{\overline\psi} + i k \overline\psi +
\overline{\psim} = 0
\qquad\text{on }\Gamma,\\
\end{cases}\label{dual:problem:first:equation:strong}\\
& \begin{cases}  -\overline \psi + (\nicefrac{1}{2} + \Kk +i k \Vk) \overline{\psitilde} + \Vk \overline{\psim} = \overline \Rm \qquad\text{on }\Gamma,
\end{cases}\label{dual:problem:second:equation:strong} \\
& \begin{cases} \left(\Wk + i k (\nicefrac{1}{2} - \Kprimek) - i k (\nicefrac{1}{2} + \Kk + i k\Vk)\right) \overline{\psitilde} - \left( (\nicefrac{1}{2} + \Kprimek) + i k \Vk\right)  \overline{\psim} = \overline{\Rext} \qquad\text{on }\Gamma.
\end{cases}\label{dual:problem:third:equation:strong:3}
\end{align}
\end{subequations}

We immediately observe that the adjoint problem
is well-posed and its solution depends continuously on the right-hand
side with a bound that is polynomial in the
wave-number:

\begin{corollary}
  \label{cor:existence_dual}
  Given $r \in L^2(\Omega)$, $\Rm \in H^{1/2}(\Gamma)$ and $\Rext \in H^{-1/2}(\Gamma)$,
  the solution $\psi \in H^1(\Omega)$, $\psim \in H^{-1/2}(\Gamma)$ and $\psiext \in H^{1/2}(\Gamma)$
  to~\eqref{eq:dual_problem} exists, is unique, and satisfies the bound
  
  \begin{align*}
  \enorm{(\psi,\psim,\psiext)} \sim   \|\psi\|_{1,k,\Omega} + \|\psim\|_{-1/2,\Gamma} + \|\psiext\|_{1/2,\Gamma}
    &\lesssim k^{\beta} \Big(\|r\|_{0,\Omega}
      + \|\Rm\|_{1/2,\Gamma} + \|\Rext\|_{-1/2,\Gamma}\Big).
  \end{align*}
  The constant $\beta$ is that of Lemma~\ref{lemma:existence_primal}.  
\end{corollary}
\begin{proof}
  Follows directly from Lemma~\ref{lemma:existence_primal} since 
  $(\overline{\psi}, -\overline{\psim},\overline{\psiext})$ solves the
  primal problem (\ref{system:Helmholtz}) with right-hand side $(\overline{r},\overline{\Rm},\overline{\Rext})$.
\end{proof}

\subsection{$H^2$ regularity of a proxy problem}
We start with a regularity result for a Poisson transmission problem.
\begin{proposition}[$H^2$-regularity for a transmission problem]
  \label{prop:ell_reg_transmission}
  Fix $R >0$ such that
  $\overline{\Omega} \subseteq B_{R/2}$, where
  $B_{R}$ is denotes a ball of radius $R>0$.
  Let $u \in H^1(B_R \setminus \Gamma)$ solve the following transmission problem:
  \begin{align*}
    -\fdiv(\nu \nabla u)
    &= f \; \text{ in $B_R$},
      \qquad u=0\;\text{ on $\partial B_R$}, \qquad \njmp{u} 
      =r_1, \qquad \djmp{u}=r_0,
  \end{align*}
  with $\nu$ symmetric and uniformly positive definite, $\nu{}_{|_{P_j}} \in C^{\infty}(\overline{P_j})$ for all $P_j \in \PP$, $\nu|_{\R^d \setminus \overline{\Omega}}=1$ and $\nu|_{\Omega}$ is smooth in a  neighborhood
  $\mathcal{N}(\Gamma)$ of $\Gamma$.
  Let $f \in L^2(B_{R})$,
  $r_1 \in H^{1/2}(\Gamma)$, and $r_0 \in H^{3/2}(\Gamma)$.
  Then $u{}_{|_{\Omega}}$ is piecewise $H^2$ on $\PP$ and
    $u{}_{|_{B_{R/2} \setminus \overline{\Omega}}}\in H^2(B_{R/2} \setminus \overline{\Omega})$
  with
  \begin{align*}
    \|u\|_{2,\PP} + \|u\|_{2,B_{R/2} \setminus \overline{\Omega}}
    &\lesssim
      \|u\|_{1,B_{R}}+ 
      \|f\|_{0,B_{R}} + \|r_1\|_{1/2,\Gamma} + \|r_0\|_{3/2,\Gamma}.
  \end{align*}
\end{proposition}
\begin{proof}
  The result follows from standard elliptic regularity theory. By lifting
  the Dirichlet jump $r_0$ to a function $R_0 \in H^2(\Omega)$
  with $\operatorname{supp}(R_0) \subseteq \mathcal{N}(\Gamma)$
  and setting $R_0=0$ in $\Omegaext$,
  we can reduce the problem to the case of $r_0=0$ by instead considering $u-R_0$ and modifying
  $f$ and $r_1$ accordingly.
  We then combine \cite[Lem.~{5.5.5}]{melenk_book} away from the interfaces
  and the estimates for transmission problems
  \cite[Lem.~{5.5.8}]{melenk_book} to deal with the interfaces.
  %
  %
\end{proof}

We consider the following auxiliary (``proxy'') problem, which represents a
``positive definite'' version of the original adjoint problem~\eqref{eq:dual_problem}:
  Given $r \in \left(H^1(\Omega)\right)'$, $\Rm \in H^{1/2}(\Gamma)$, and $\Rext \in H^{-1/2}(\Gamma)$
with $(\Rext,1)_{0,\Gamma}=0$, let
$S^+_k(r,\Rm,\Rext):=(\phi,\phim,\phiext)$ solve
 \begin{subequations}
    \label{eq:modified_problem}    
    \begin{align}      
      &\begin{cases}
        -\fdiv(\Ad \nabla\overline  \phi) + (\kn)^2 \overline \phi =  \overline r \quad \qquad\qquad\;\text{in }\Omega, \\
        \gammanu[int]{\overline\phi}   + \overline{\phim} = 0 \qquad\qquad\qquad\qquad\text{on }\Gamma,\\
\end{cases}\label{modified:problem:first:equation:strong}\\
      & \begin{cases}  -\overline \phi + (\nicefrac{1}{2} + \Kz)  \overline{\phiext} + \Vz \overline{\phim} = \overline \Rm \qquad\text{on }\Gamma,
      \end{cases}\label{modified:problem:second:equation:strong} \\
      & \begin{cases}  \Wz \overline{\phiext} - (\nicefrac{1}{2} + \Kprimez)   \overline{\phim} = \overline{\Rext} \qquad\qquad\text{on }\Gamma.
      \end{cases}\label{modified:problem:third:equation:strong:3}
\end{align}
\end{subequations}
\begin{lemma}
  \label{lemma:mod_apriori}
The operator $S^+_k$ is well defined, i.e.,~\eqref{eq:modified_problem} has
a unique solution. Moreover, 
for $0\le\tau\le 1$ and $r\in\left(H^\tau(\Omega)\right)'$,
this solution satisfies the \textsl{a priori} estimate
\begin{align}
  \label{eq:mod_apriori_h1}
  \|\phi \|_{1,k, \Omega} + \|\phim\|_{-1/2,\Gamma} + \|\phiext\|_{1/2,\Gamma}
  &\lesssim k^{-1+\tau}\|r\|_{\left(H^\tau(\Omega)\right)'}   + \|\Rm\|_{1/2,\Gamma} + \|\Rext\|_{-1/2,\Gamma}.
\end{align}
\end{lemma}
\begin{proof}
We consider the weak formulation of
problem~\eqref{eq:modified_problem} and prove that the sesquilinear form
  \begin{align*}
    \TT_{+}((v,\lambda,\vext),(\phi,\phim,\phiext))
    &:=
    (\nu \nabla v, \nabla \phi)_{0,\Omega} +  ((\kn)^2 v, \phi)_{0,\Omega}
      + \langle
      v, \phim \rangle_\Gamma \\
    &-\langle\lambda, \phi \rangle_\Gamma
      +\langle\lambda, \big(\nicefrac{1}{2}+ \Kz\big)\phiext \rangle_\Gamma
      +\langle\lambda, \Vz\phim\rangle_\Gamma  \\
    &+ \langle\vext, \Wz \phiext \rangle_\Gamma
      - \langle \vext, \big(\nicefrac{1}{2} + \Kprimez\big) \phim \rangle_\Gamma
  \end{align*}
  is coercive with respect to the space
  $$
  H^1(\Omega) \times H^{-1/2}(\Gamma) \times H_\star^{1/2}(\Gamma)
  \quad \text{with} \quad H_\star^{1/2}(\Gamma):=\big\{ \zeta \in H^{1/2}(\Gamma):\; (\zeta,1)_{0,\Gamma}=0 \big\}.
  $$
  (We used the fact that, since the Green function of the
  Laplacian is real valued,  all the boundary operators
  satisfy $\overline{\Vz\phim}
  =\Vz \overline{\phim}$ etc.)
  Choosing $(v,\lambda,\vext)=(\phi,\phim,\phiext)$ we get:
  \begin{align*}
    \TT_{+}((\phi,\phim,\phiext),(\phi,\phim,\phiext))
    &=
      \|\nu^{1/2}\nabla \phi\|^2_{0,\Omega} + k^2 \|n \phi\|^2_{0,\Omega}
      + \langle 
      \phi, \phim \rangle_\Gamma \\
    &-\langle\phim, 
      \phi\rangle_\Gamma
      +\langle\phim, \big(\nicefrac{1}{2}+ \Kz\big)\phiext \rangle_\Gamma
      +\langle\phim, \Vz\phim\rangle_\Gamma  \\
    &+ \langle\phiext, \Wz \phiext \rangle_\Gamma
      - \langle \phiext, \big(\nicefrac{1}{2} + \Kprimez\big) \phim \rangle_\Gamma.
  \end{align*}
  Using the fact that $\Kprimez$ is the adjoint of $\Kz$, we take the real part and are left
  with
  \begin{align*}
    \Re\big(\TT_{+}((\phi,\phim,\phiext),(\phi,\phim,\phiext))\big)
    &=\|\nu^{1/2}\nabla \phi\|^2_{0,\Omega} + k^2 \|n\phi\|_{0,\Omega}^2
      + \langle \phim,\Vz \phim\rangle +
      \langle \phiext, \Wz \phiext \rangle \\
    &\gtrsim \|\phi\|_{1,k,\Omega}^2 + \|\phim\|^2_{-1/2,\Gamma} + \|\phiext\|_{1/2,\Gamma}^2,
  \end{align*}
  where in the last step we used the coercivity of $\Vz$ and $\Wz$ and
  the fact that $\phiext \in H^{1/2}_\star(\Gamma)$.

The continuity of the sesquilinear form
    $\TT_{+}(\cdot,\cdot)$ in the considered space follows from standard estimates and mapping properties of the operators.

  What remains to be shown in order to conclude the proof is that
  the right-hand side
  of the weak formulation of~\eqref{eq:modified_problem} is a bounded linear functional.
Let $v\in H^1(\Omega)$. For $\tau=0$, we have
\[
\big|(v,r)_{0,\Omega}\big|\le
k^{-1}\|r\|_{0,\Omega}\left(k\|v\|_{0,\Omega}\right)\le
k^{-1}\|r\|_{0,\Omega}\|v\|_{1,k,\Omega},
\]
and for $\tau=1$, we have
\[
\big|_{H^1(\Omega)}\langle v,r\rangle_{\left(H^1\Omega)\right)'}\big|\le
\|r\|_{\left(H^1\Omega)\right)'}\|v\|_{H^1(\Omega)}\le
\|r\|_{\left(H^1\Omega)\right)'}
\|v\|_{1,k,\Omega}.
\]
By interpolation, we get
\[
\big|_{H^\tau(\Omega)}\langle v,r\rangle_{\left(H^\tau\Omega)\right)'}\big|\le
k^{-1+\tau}
\|r\|_{\left(H^\tau\Omega)\right)'}
\|v\|_{1,k,\Omega}.
\]
This, together with
the Cauchy-Schwarz inequality, gives
  \begin{multline*}
    \big|{_{H^\tau(\Omega)}\langle v,r\rangle_{\left(H^\tau\Omega)\right)'}} + \langle\lambda,\Rm\rangle + \langle\vext,\Rext\rangle\big|\\
    \leq k^{-1+\tau} \|r\|_{\left(H^\tau\Omega)\right)'} \|v\|_{1,k,\Omega}
    + \|\Rm\|_{1/2,\Gamma} \|\lambda\|_{-1/2,\Gamma} + \|\Rext\|_{-1/2,\Gamma} \|\vext\|_{1/2,\Gamma}.
  \end{multline*}
  From this, the \textsl{a priori} estimate for the weak solution follows readily. The
  standard compatibility condition $(\Rext ,1)_{0,\Gamma}=0$  
allows us to extend the test functions from $H^{1}(\Omega) \times H^{-1/2}(\Gamma) \times H^{1/2}_\star(\Gamma)$ to 
$H^{1}(\Omega) \times H^{-1/2}(\Gamma) \times H^{1/2}(\Gamma)$ so that we indeed have a weak solution of 
    (\ref{eq:modified_problem}).
\end{proof}

\begin{lemma}
  \label{lemma:mod_aprioi_h2}
  Let $S^+_k \bs r
  := S^+_k(r,\Rm,\Rext)=\boldsymbol{\phi}:=(\phi,\phim,\phiext) $ solve
  \eqref{eq:modified_problem}
    For $r \in L^2(\Omega)$,
  $\Rm \in H^{3/2}(\Gamma)$, and $\Rext \in H^{1/2}(\Gamma)$
  with $\langle \Rext,1\rangle_\Gamma=0$.

Then, the following shift theorem is valid:
  \begin{align*}
    \smalltriplenorm{}{\boldsymbol{\phi}}_{k,\mathcal{V},1}
    &\lesssim \smalltriplenorm{}{\bs r}_{k,\mathcal{V}',1}.
  \end{align*}
 The implied constant is independent of $k$.
\end{lemma}
\begin{proof}
We first note that, for the lowest order terms present in the norms,
  we can use Lemma~\ref{lemma:mod_apriori} and obtain
  \begin{align*}
    k^{2}\|\phi\|_{0,\Omega} + k\|\phim\|_{-1/2,\Gamma} + k\|\phiext\|_{1/2,\Gamma}
    \lesssim \smalltriplenorm{}{\bs r}_{k,\mathcal{V}',1}.
  \end{align*}
  Thus, we can focus on the highest-order terms only, i.e., 
  we need to estimate
  $$
  \|\phi\|_{2,\PP} + \|\phim\|_{1/2,\Gamma} + \|\phiext\|_{3/2,\Gamma}.
  $$ 
  We proceed similarly to~\cite[Thm.~{3.11}]{MMPR20}, aiming to
  represent $\phim$ and $\phiext$ as traces of functions on $\R^d \setminus \Gamma$
  and then decomposing these functions into parts that can be analyzed more easily.
  In order to get control of the boundary traces $\phim$ and $\phiext$
  we lift them to the volume using the Laplace potentials:
    \begin{align*}
      \daleth&:=\Vtildez \overline{\phim}  + \Ktildez \overline{\phiext}.
    \end{align*}
  By the jump conditions~\eqref{eq:jump_conditions}, we deduce that
  \begin{align*}
    \djmp{\gammao \daleth}=\overline{\phim}, \qquad \djmp{\daleth}=-\overline{\phiext}.
  \end{align*}
  Taking the interior Neumann trace and using~\eqref{modified:problem:third:equation:strong:3} give
  \begin{align}
    \gammaoint \daleth
    &=\big(\nicefrac{1}{2}  + \Kprimez \big) \overline{\phim}  - \Wz \overline{\phiext}
      = -\overline{\Rext}.
      \label{eq:mod_daleth_neumann_trace}
  \end{align}
  Since the Laplace potential is harmonic,
  we can determine $\daleth$ in $\Omega$ as the solution of a standard
  Laplace-Neumann problem with data $\Rext \in H^{1/2}(\Gamma)$.
The assumption $\langle\Rext,1\rangle_\Gamma=0$ guarantees
  compatibility of the data.
  This gives
  by elliptic regularity and a trace estimate
  \begin{align}\label{eq:dalet2}
    \|\daleth\|_{2,\Omega}  + \|\gammanu[int]{\daleth}\|_{1/2,\Gamma}+\|\gammazint\daleth\|_{3/2,\Gamma}\lesssim
    \|\Rext\|_{1/2,\Gamma}.
  \end{align}
  In order to obtain control of the exterior contribution,
  we compute $\gammazext \daleth$,
    using~\eqref{modified:problem:second:equation:strong}:  
  \begin{align}\label{eq:daleth_ext}
    \gammazext \daleth
    &=\Vz \overline{\phim} + (\nicefrac{1}{2} + \Kz)\overline{\phiext}
    = \overline{\phi}+ \overline{\Rm}.
  \end{align}

Next, we want to eliminate the dependence on $\overline{\phi}$
  and only rely on the given data $(r, \Rm,\Rext)$.
  This can be achieved by using the auxiliary function
  \begin{align*}
          \LL:=\overline{\phi} \mathds{1}_\Omega+ \daleth \,\chi_{\mathcal{N}(\Gamma)},
  \end{align*}
  where $\mathds{1}_{\Omega}$ denotes the characteristic function of the set $\Omega$
  and $\chi_{\mathcal{N}(\Gamma)}$ is a smooth cutoff function
  that is equal to $1$ in a (sufficiently small) neighborhood of $\Gamma$ so that on $\Omega$ its support is restricted to where $\nu|_{\Omega}$ is smooth by assumption.
  
  From estimate~\eqref{eq:mod_apriori_h1} and the mapping properties of $\Vtildez$ and
  $\Ktildez$ from Proposition~\ref{prop:mapping_properties_laplace_potentials}, we directly get,
  for any fixed ball $B_R \supseteq \overline{\Omega}$,
  \begin{align}
    \label{eq:LL_est_h1}
  \|\LL\|_{1,B_R \setminus \Gamma}
  \lesssim  k^{-1}\|r\|_{0,\Omega} +   \|\Rm\|_{1/2,\Gamma} + \|\Rext\|_{-1/2,\Gamma}.
  \end{align}
  
  By the identity~\eqref{eq:daleth_ext}, we 
  get~$\djmp{\LL}=\gammazint \daleth-\overline{\Rm}$.
  To compute the Neumann jump, we take the exterior Neumann trace of $\daleth$,
use~\eqref{modified:problem:third:equation:strong:3} and~\eqref{modified:problem:first:equation:strong}, and obtain
  \begin{align*}
    \gammaoext \daleth
    &=\big(-\nicefrac{1}{2}  + \Kprimez \big) \overline{\phim}  - \Wz \overline{\phiext} 
    = -\overline{\Rext} - \overline{\phim} \\
    &= -\overline{\Rext} + \gammanu[int]{\overline{\phi}} .
  \end{align*}
  Using this and~\eqref{eq:mod_daleth_neumann_trace}, we get 
  \[
  \njmp{\LL}=\gammanu[int]{\overline{\phi}}
  +\gammanu[int]{\daleth}
  -\gammaoext \daleth
  =
  \gammanu[int]{\overline{\phi}}
  +\gammanu[int]{\daleth}
  +\overline{\Rext} -\gammanu[int]{\overline{\phi}} 
  =\gammanu[int]{\daleth}
  +\overline{\Rext}.
\]
  Overall, the function $\LL$ solves
  \begin{align*}
    -\fdiv(\nu \nabla \LL)
    &= \big(\overline{r} - {(kn)}^2 \overline{\phi} - 
    \fdiv \big(\nu \nabla ( \daleth \,\chi_{\mathcal{N}(\Gamma)})\big) \big) \mathds{1}_{\Omega},\\
    \djmp{\LL} &= \gammazint\daleth - \overline{\Rm}, \quad \njmp{\LL}=
                 \gammanu[int]{\daleth}
  +\overline{\Rext}
  \end{align*}

  Thus, we can apply the elliptic regularity result from Proposition
  \ref{prop:ell_reg_transmission} to get
  \begin{align*}
    \| \LL \|_{2,\PP} + \| \LL \|_{2,B_R \setminus \overline{\Omega}}
    &\lesssim  \|\LL\|_{1,B_R}+ \|r\|_{0,\Omega} + k^2 \|n^2\phi\|_{0,\Omega}
      + \|\nabla \daleth \|_{0,\Omega} 
      + \| \Rext \|_{1/2,\Gamma}
      + \| \gammanu[int]{\daleth}\|_{3/2,\Gamma}
    \\
    &\qquad\qquad+ \| \gammazint \daleth\|_{3/2,\Gamma}
      + \| \Rm \|_{3/2,\Gamma}
    \\
    &\!\!\!\!\!\!\!\!\!\!\!\!\!\!\stackrel{\eqref{eq:mod_apriori_h1},\,\eqref{eq:LL_est_h1},\,\eqref{eq:dalet2}}
      {\lesssim}
      \|r\|_{0,\Omega} + (1+k) \big(k^{-1}\|r\|_{0,\Omega}      
      + \|\Rm\|_{1/2,\Gamma}
      + \|\Rext\|_{-1/2,\Gamma}\big) \\
    &\qquad \qquad
      +\|\Rext\|_{1/2,\Gamma}
      +\|{\Rm}\|_{3/2, \Gamma}.
  \end{align*}
Thus, in $\Omega$, we get
  \begin{align*}
    \|\phi\|_{2,\PP}
    &\lesssim \|\LL\|_{2,\PP} + \|\daleth\|_{2,\Omega} \\
      &\lesssim
        \|r\|_{0,\Omega}      
        + k\|\Rm\|_{1/2,\Gamma}
        + \|\Rm\|_{3/2,\Gamma}
      + k\|\Rext\|_{-1/2,\Gamma}
      + \|\Rext\|_{1/2,\Gamma}.
  \end{align*}
  In $\Omegaext$, we have $\LL=\daleth$ and thus we get
  control of the Dirichlet and Neumann jumps of~$\daleth$---and hence of $\phi$, $\phim$---from standard
  trace estimates for $\LL$. 
\end{proof}

\subsection{Analytic regularity}
In this section, we show that for analytic data, the solution to the adjoint
  problem is also analytic.  

The following lemma analyzes how analycity classes are transformed if functions are multiplied by analytic functions or if derivatives are taken: 
\begin{lemma}
  \label{lemma:product_and_derivative_of_analytics}
  Let $\mathcal{O}$ be an open set,
  $u \in \aclass\big(C_u,\vartheta_u,\mathcal{O}\big)$, and $\eta \in \aclass^{\infty}\big(C_{\eta},\vartheta_\eta,\mathcal{O}\big)$.

  Then $\nabla u$ and $\eta u$ are both analytic and
  \begin{align*}
    \nabla u &\in \aclass\big({C_{\nabla u},\vartheta_{\nabla u}\mathcal{O}}\big), \qquad \text{with } \qquad C_{\nabla u } \leq \max(k,1) e^2 \vartheta_{u} C_u \text{ and } \vartheta_{\nabla u} \leq e \vartheta_u,\\
    \eta u &\in \aclass\big({C_{\eta u},\vartheta_{\eta u},\mathcal{O}}\big) \qquad \text{with } \qquad \;C_{\eta u } \leq C_{\eta} C_{u} \text{ and } \vartheta_{\eta u} \leq \vartheta_\eta+\vartheta_u.
  \end{align*}
\end{lemma}
\begin{proof}
  The product of analytic functions is again analytic (see~{\cite[Lem.~{2.6}]{MS21}}). The precise behavior of the constants
  for $\eta u$  can be found by tracking the proof.
  For $\eta u$, we compute
  \begin{align*}
    \|\nabla^{p}(\nabla u \big)\|_{0,\mathcal{O}}
    &\leq
      \|\nabla^{p+1}(\nabla u \big)\|_{0,\mathcal{O}}
      \leq C_u \vartheta_u^{p+1} \max(p+2,k)^{p+1}\\
    &\leq \big[C_{u} \vartheta_{u} e\big] \max(k,1) 2^p\vartheta_u^{p} \max(p+1,k)^{p},
  \end{align*}
  where in the last step we used the crude estimate $(p+2)\leq e^{p+1}$ to estimate
  \begin{align*}
  (p+2)^{p+1} &= (p+1)^p (p+2) \Big(1+\frac{1}{p+1}\Big)^{p+1} \leq
  (p+1)^p \,e^{p+2}. \qedhere
  \end{align*}
  \end{proof}

\begin{proposition}[Regularity of a transmission problem]
  \label{prop:analytic_reg_transmission}
  Let $\mathcal{N}(\Gamma)$ be a neighborhood of $\Gamma$.
  Let $u \in H^1(\mathcal{N}(\Gamma) \setminus \Gamma)$ solve the following transmission problem:
  \begin{align*}
    -\fdiv(\nu \nabla u)-(kn)^2u
    &= f \; \text{ in $\mathcal{N}(\Gamma)$},
    \quad \njmp{ u+i k u}=\gammazint R_1, \quad \djmp{u}=\gammazint R_0
  \end{align*}
  with $\nu{} \in
  \aclass^{\infty}(C_{\nu},\vartheta_\nu,
  \mathcal{N}(\Gamma) \setminus \Gamma)$, $\nu{}_{|_{\R^d\setminus \overline{\Omega} }}=1$,
  and $f \in \aclass(C_f,\vartheta_f,\mathcal{N}(\Gamma) \setminus \Gamma)$ 
as well as $n \in \aclass^{\infty}(C_{n}, \vartheta_n,\mathcal{N}(\Gamma)\setminus \Gamma)$
  Assume $R_j \in \aclass(C_{R_j},\vartheta_{R_j},\mathcal{N}(\Gamma),\Gamma)$ for $j=0$, $1$.
  Then $u \in \aclass(C,\vartheta,\mathcal{N}(\Gamma) \setminus \Gamma)$ with
  \begin{align*}
    C\leq k^{-1} \|\nabla u\|_{0,B_R\setminus\Gamma}+ \|u\|_{0,B_R\setminus \Gamma}+C_{R_0}+k^{-1} C_{R_1}+k^{-2} C_{f}
  \end{align*}
  and $\vartheta >0$ depending only on $\nu$, $n$ and $\vartheta_f$, $\vartheta_{R_0}$, $\vartheta_{R_1}$.
\end{proposition}
\begin{proof}
  We start with the simple case of a localized version on a half-ball $ B_{R}^+:=\{x \in B_{R}: x_{d}>0\}$ and
  $\Gamma=B_{R} \cap \{x_d=0\}$, and also assume that
  the data and $\nu$ are analytic in $B_{R} \setminus \Gamma$.
We further assume that
  both $R_0$ and $R_1$ are traces of functions on all of $B_R^+$. 
  We decompose $u$ into two functions. First, we subtract $R_0 \mathds{1}_{B_{R}^+}$
  to get that $u_2:=u-R_0 \mathds{1}_{B_{R}^+}$ solves

  \begin{align*}
  -\fdiv(\nu \nabla u_2) - (kn)^2 u_2&= f + \big(\fdiv(\nu \nabla R_0) + (kn)^2 R_0 \big)\mathds{1}_{B_R^+},  \\
    \llbracket \partial_{x_d} u_2 \rrbracket_\Gamma&=\gammazint R_1- i k \gammazint R_0 +  \partial_{x_d}{R_0}, \qquad \djmp{u_2}=0.
  \end{align*}
  We can then apply \cite[Prop.~{5.5.4}]{melenk_book} with $\varepsilon=k^{-1}$,
  $G_2:=0$ and
  $
  G_1=k^{-1} \big(R_1 - i k R_0  -\gammanu[int]{R_0}\big).
  $
  The right-hand side becomes
  $$\widetilde{f}:=k^{-2} f + k^{-2}\fdiv(\nu \nabla R_0)\mathds{1}_{B_{R}^+} + n^2 R_0 \mathds{1}_{B_R^+}.$$
  By Lemma~\ref{lemma:product_and_derivative_of_analytics}, we get that
  $  \widetilde{f} \in \aclass{\big(C_{\widetilde{f}},\vartheta_{\widetilde{f}},\mathcal{N}(\Gamma) \setminus \Gamma\big)}$
  with
  \begin{align*}
          C_{\widetilde{f}}\lesssim k^{-2}C_f + C_{\nu} C_{R_{0}} + C_{n}^2 C_{R_0}
  \quad &\text{and} \quad \vartheta_{\widetilde{f}} \leq \max(\vartheta_f, \vartheta_{R_0}e+\vartheta_\nu,2\vartheta_{n}+\vartheta_{R_0}).
  \end{align*}
  
  This gives the following translation table for the constants involved
  when applying \cite[Prop.~{5.5.4}]{melenk_book}:
  \begin{alignat*}{5}
    C_A&:=C_\nu,&  \quad C_b&:=0, &\quad  C_c&:=C_{n^2}, \\
    C_f &:= C_{\widetilde{f}}, & \quad
    \gamma_f&:=\vartheta_{\widetilde{f}}, &\quad
    C_{G_1}&:=k^{-1} C_{R_1} + C_{R_0} + k^{-1}C_{\gammanu[int]{R_0}},\\
     C_{G_2}&:=0, & \quad
    \gamma_{G_1}&:=\max(\vartheta_{R_1},\vartheta_{R_0},\vartheta_{\gammanu[int]{R_0}}), &\quad
    \gamma_{G_2}&:=0.
  \end{alignat*}
  The resulting estimate can be summarized as
  \begin{align*}
    \|\nabla^p u\|_{0,B_R \setminus \Gamma}
    &\leq C_u K^p \max(p+1,k)^p \qquad
      C_u\lesssim k^{-1}\|\nabla u\|_{0,\Omega}
      + (\|u\|_{0,\Omega}+C_{\widetilde{f}})
      + C_{G_1}.
  \end{align*}
  where the implied constants depends only on the geometry and on the coefficients $\nu$ and $n$; the constant 
  $K$ depends additionally on $\gamma_{f}$ and $\gamma_{G_1}$.
   
  For the case of general $\mathcal{N}({\Gamma})$, we cover $\mathcal{N}({\Gamma})$ with sufficiently small open sets $\mathcal{O}_j, j=1,\dots,N$.
  If $\mathcal{O}_j$ intersects $\Gamma$, we assume it is small enough so that $R_0$, $R_1$ are defined 
  on all of $\mathcal{O}_j$ and $\nu|_{\Omega \cap \mathcal{O}_j}$ is analytic. In addition we require that $\mathcal{O}_j$ is analytically mapped to the
  unit sphere $B_1$ such that $\Gamma \cap \mathcal{O}_j$ is mapped to $B_1 \cap \{x_d = 0\}$.
  By \cite[Lem.~{2.6}]{MS21}, we can then transform the problem to the special case already
  covered.
  For points away from the interfaces,  $\mathcal{O}_j$ is taken as a simple sphere and we use
  standard interior analytic regularity results of \cite[Prop.~{5.5.1}]{melenk_book}.
  Summing up over all such sets $\mathcal{O}_j$, we get the stated general result.
\end{proof}

\begin{lemma}
  \label{lemma:analytic_regularity}
  Let $r\in \aclass(C_r,\vartheta_r,\PP)$. 
  Let $\mathcal{O}$ be a neighborhood of $\Gamma$ such that the normal vector map
  $x \mapsto \nGamma(x)$ has an analytic continuation to $\mathcal{O}$ and $\nu$ is
    analytic in $\mathcal{O}$.
    
  Let $\Rm \in \aclass(C_{\Rm},\vartheta_{\Rm}, \mathcal{O},\Gamma)$,
  $\Rext \in \aclass(C_{\Rext},\vartheta_{\Rext},\mathcal{O},\Gamma)$.
  Then, the solution $(\psi,\psim,\psiext)$ to the adjoint problem~\eqref{eq:dual_problem} satisfies
  $$
  \psi \in \aclass(C,\vartheta, \PP), 
\quad
 \psim \in \aclass(kC,\vartheta, \mathcal{O}, \Gamma),
 \quad \psiext \in \aclass(C,\vartheta, \mathcal{O}, \Gamma)
  $$
  with $C \lesssim \max(k^{\beta},1)( C_r+ C_{\Rm}+C_{\Rext})$ and $\vartheta>0$ depending only on $\nu$, $n$, $\vartheta_r$,
  $\vartheta_{\Rm}$, $\vartheta_{\Rext}$, $k_0$ but independent of $k$. 
\end{lemma}
\begin{proof} 
  We follow the proof of Lemma~\ref{lemma:existence_primal},
  but using analytic regularity results
  instead of only showing  estimates in the energy norm. 
  First we note that away from $\Gamma$ we can use standard 
  analytic regularity theory to get the stated estimates. 
(See, e.g., \cite[Sec.~{5.1}]{bernkopf}, which in turn is based on \cite[Prop.~{5.5.1}]{melenk_book} for interior regularity and \cite[Prop.~{5.5.4}]{melenk_book} 
for transmission problems in conjunction 
with the invariance of the analyticity classes $\A$ under analytic changes of variable, \cite[Lem.~{2.6}]{MS21}.)
  We therefore focus on a neighborhood $\mathcal{N}(\Gamma)$ of the boundary $\Gamma$.
  We again use the auxiliary functions
  \begin{align*}
    \daleth&:=\Vtildek (\overline{\psim} +  \ii k \overline{\psiext}) + \Ktildek \overline{\psiext}
             \qquad \text{and  } \qquad
             \LL:=\overline{\psi} \mathds{1}_{\Omega}+ \daleth.
  \end{align*}
  In Lemma~\ref{lemma:existence_primal}, it was shown that $\daleth|_{\Omega}$ is the
  solution to  a Robin-type boundary value problem.
  By combining~\cite[Proof of Lem.~{4.13}]{MS11} and the stability bound~\cite[Thm.~{1.8}]{BSW16},
  there exist constants $C_{\daleth}$ and $\vartheta_{\daleth}$ such that
  \begin{align*}
    \daleth|_{\Omega} \in \aclass(C_\daleth,\vartheta_{\daleth},\Omega) \qquad \text{with} \quad
    C_\daleth \lesssim k^{-1} \|\Rext\|_{1/2,\Gamma} \leq  C_{\Rext}.
  \end{align*}

  As in Lemma~\ref{lemma:existence_primal}, we establish
  that the function $\LL$ solves the following transmission problem:
  \begin{align*}
    - \fdiv(\nu \nabla \LL) - (\kn)^2 \LL
    &= \overline{r}-\big[\fdiv(\nu \nabla \daleth)  - (\kn)^2
      \daleth  \big] \mathds{1}_{\Omega} =: \tilde{f}\\
    \njmp{\LL+ i k \LL}&=\gammanu[int]{\daleth} - \gammaoint{\daleth}
        \quad \text{and} \quad
    \djmp{\LL}=\gammazint \daleth - \overline{\Rm},
  \end{align*}
  
  For the new right-hand side $\tilde{f}$ we compute, since
  we have already established control of $\daleth{}_{|_{\Omega}}$
  by Lemma~\ref{lemma:product_and_derivative_of_analytics}:
  \begin{align*}
    \widetilde{f} \in \aclass \Big( C( C_r + k^2 C_{\daleth} ),
    e\vartheta_{\nu}+(1+e^2)\vartheta_{\daleth}+\vartheta_n +\vartheta_r,\PP\Big),
  \end{align*}
  with a constant $C$ only depending on $C_{\nu}, C_{n}$ and $\vartheta_{\nu}$, $\vartheta_n$.

  Since the analyticity classes are closed under multiplication (see Lemma~\ref{lemma:product_and_derivative_of_analytics}),
  and we can extend $\nGamma$ to $\mathcal{O}$,
  the impedance jump
  $R_1:=\gammanu[int]{\daleth} - \gammaoint{\daleth}$
  and Dirichlet jump $R_0:=\overline{\Rm} + (1-\nu) \gammazint{\daleth}$
  satisfy
  \begin{align*}
    R_1 &\in \aclass \big( C( C_{\nu}C_{\daleth} k), c_{\vartheta}(\vartheta_{\daleth}+\vartheta_{\nu}),\mathcal{O},\Gamma\big). \\
    R_0 &\in \aclass \big(C(C_{\Rm}+  C_{\daleth}), c_{\vartheta}\max(\vartheta_{\Rm},\vartheta_{\daleth}),\mathcal{O},\Gamma\big).
  \end{align*}
  We can thus apply Proposition~\ref{prop:analytic_reg_transmission}.
  Inserting the easily  derived 
  stability estimate $\|\LL\|_{1,k,\R^d \setminus\Gamma} \lesssim k^{\beta-1}(C_r + k C_{\Rext}+k C_{\Rm})$ (see Corollary~\ref{cor:existence_dual}), we get:
  $$
  \LL \in \aclass \big(C_{\LL}, \vartheta_{\LL},B_R \setminus \partial \PP\big)
  \qquad \text{ with } \qquad C_{\LL} \leq C k^{\beta} ( C_r + C_{\Rm}+C_{\Rext}),
  $$
  where the constant $C$ depends on $\nu$, $n$ and $k_0$ but is independent
  of $k$.
  The stated results then follow by writing $\overline{\psi}=\LL|_{\Omega} - \daleth|_{\Omega}$ and
  standard trace estimates. Note that the extra power $k$ for bounding $\psi_m$ comes from
    the fact that it involves a normal derivative.
\end{proof}

\subsection{A decomposition result}
In this section, we prove our main result concerning the
adjoint problem. Namely, we prove that its solution admits a
decomposition into
a finite regularity part with good $k$-dependence and an analytic remainder.
This strategy closely follows what was done in previous works~\cite{bernkopf,MS11,MS21}.

We start with some preliminary lemmas.
\begin{lemma}
  \label{lemma:op_decomposition_1}
  We have the following operator decompositions:
 \begin{align*}
  (\Kk - \Kz + \ii k \Vk) \psi &=  \mathcal{S}_{1} \psi
                                 + {\mathcal{A}}_{1} \psi,
  \\
   (\Vk -\Vz) \lambda&=
    \mathcal{S}_{2}\lambda
                       +{\mathcal{A}}_{2} \lambda,\\
    ik \Vk\psi &= \mathcal{S}_{3} \psi
    +{\mathcal{A}}_{3} \psi,
 \end{align*}
 with the following estimates:
 \begin{align*}
   &\|\mathcal{S}_{1}\|_{H^{1/2}(\Gamma) \to H^{3/2}(\Gamma)}
   + \|\mathcal{S}_{2}\|_{H^{-1/2}(\Gamma) \to H^{3/2}(\Gamma)}
   +\|\mathcal{S}_{3}\|_{H^{1/2}(\Gamma) \to H^{3/2}(\Gamma)}
   \lesssim k,    \\[0.2cm]
   &\mathcal{A}_1 \psi \in \mathcal{A}\big(C_{1}k^2\|\psi\|_{-1/2,\Gamma},\vartheta_{1},\Omega, \Gamma\big),
   \quad
   \mathcal{A}_2 \lambda\in \mathcal{A}\big(C_{2}k\|\lambda\|_{-3/2,\Gamma},\vartheta_{\mathcal{V}},\Omega, \Gamma\big),
   \quad
   \mathcal{A}_3 \psi \in \mathcal{A}\big(C_{3}k^2\|\psi\|_{-1/2,\Gamma},\vartheta_{\mathcal{V}},\Omega, \Gamma\big).
 \end{align*}
\end{lemma}

 \begin{proof}
   The fact that we can decompose $\Kk-\Kz$  and
   $\Vk-\Vz$  as stated follows from
   Proposition~\ref{prop:decomposition}, using decompositions~\eqref{eq:differences-c}
     and~\eqref{eq:differences-a},
     estimates~\eqref{eq:estimates-S-c} and~\eqref{eq:estimates-S-a}
     with $s=0$ and $s^\prime=2$,
     as well as~\eqref{eq:estimates-S-f} and~\eqref{eq:estimates-S-e}.
   
   For the term $\ii k \Vk\psi$,
   using again~\eqref{eq:differences-a}, we obtain
   \begin{align*}
     \ii k \Vk \psi 
     &= \ii k  \Vz \psi 
       + \ii k (\Vk-\Vz) \psi 
       = \ii k  \Vz \psi 
       + \ii k \SV \psi 
       + \ii k \gamma_0\AVtilde \psi.
   \end{align*}
Owing to Proposition~\ref{prop:mapping properties_laplace_smooth}
  with $s=3/2$, estimate~\eqref{eq:estimates-S-a} with $s=1$ and
  $s^\prime=2$, and~\eqref{eq:estimates-S-e},
the first two terms have the correct mapping properties and the
third one is analytic.
 \end{proof}

 \begin{lemma}
   \label{lemma:op_decomposition_2}
   We have the following operator decompositions:
   \begin{subequations}
   \begin{align}
     \big((\Wz - \Wk) + i k (\Kprimek + \Kk) - k^2 \Vk\big) \psi
     &=
       \mathcal{S}_{4} \psi
       + {\mathcal{A}}_{4} \psi,
       \label{eq:op_splitting_3}
      \\
     \big((\Kprimek - \Kprimez) + i k \Vk )  \lambda
     &=    \mathcal{S}_{5}\lambda
       + {\mathcal{A}}_{5} \lambda,
       \label{eq:op_splitting_4}
       \\
     -\big(ik(\nicefrac{1}{2}+\Kprimek)-k^2\Vk\big)\psi
     &=\mathcal{S}_{6} \psi
       + {\mathcal{A}}_{6} \psi,
       \label{eq:op_splitting_5}
   \end{align}
 \end{subequations}
 with the following estimates:
 \begin{align*}
   &\|\mathcal{S}_{4}\|_{H^{1/2}(\Gamma) \to H^{1/2}(\Gamma)}
   + \|\mathcal{S}_{5}\|_{H^{-1/2}(\Gamma) \to H^{-1/2}(\Gamma)}
   +\|\mathcal{S}_{6}\|_{H^{1/2}(\Gamma) \to H^{1/2}(\Gamma)}
   \lesssim k,    \\[0.2cm]
   &\mathcal{A}_4 \psi \in \mathcal{A}\big(C_{4} k^3\|\psi\|_{-1/2,\Gamma},\vartheta_{4},\mathcal{O}, \Gamma\big),
   \quad
   \mathcal{A}_5\in \mathcal{A}\big(C_{5}k^2\|\lambda\|_{-3/2,\Gamma},\vartheta_{5},\mathcal{O},\Gamma\big),
   \quad
   \mathcal{A}_6 \psi \in \mathcal{A}\big(C_{6} k^3\|\psi\|_{-1/2,\Gamma},\vartheta_{6},\mathcal{O},\Gamma\big)
 \end{align*}
for a fixed (unilateral) tubular neighborhood $\mathcal{O} \subset \Omega$ of $\Gamma$ that depends solely on $\Omega$.
\end{lemma}

\begin{proof}
  We prove decomposition~\eqref{eq:op_splitting_3} with the
    corresponding estimates. 
   The operator $\Wk-\Wz$ can be split according to the statement of this lemma
   via Proposition~\ref{prop:decomposition},
     using decomposition~\eqref{eq:differences-d}, estimate~\eqref{eq:estimates-S-d}
     with $s=0$ and $s^\prime=2$, as well as~\eqref{eq:estimates-S-f} and~\eqref{eq:estimates-S-e}.
   For $\Kprimek$ and $\Kk$, we can split
   \begin{align*}
     \Kprimek \psi&= \Kprimez \psi + (\Kprimek-\Kprimez) \psi,\qquad
     \Kk \psi= \Kz \psi + (\Kk-\Kz) \psi.  
   \end{align*}
   The 
   operators $\Kprimez$ and $\Kz$ are bounded in $H^{1/2} \to H^{1/2}$
   by Proposition~\ref{prop:mapping properties_laplace_smooth}. The
   operators $\Kprimek-\Kprimez$ and $\Kk-\Kz$ 
   can be split into a finite regularity part and an analytic remainder
   by Proposition~\ref{prop:decomposition}
   (see~\eqref{eq:differences-b}
     and~\eqref{eq:differences-c} and note that the presence of the normal vector $\nGamma$ mandates the restriction
to a sufficiently small tubular neighborhood $\mathcal{O} \subset \Omega$ of $\Gamma$ where $\nGamma$ is analytic). The estimates follow
from~\eqref{eq:estimates-S-b} with $s=1$ and $s^\prime=2$,~\eqref{eq:estimates-S-c} with $s=0$
  and $s^\prime=1$, and from~\eqref{eq:estimates-S-f} and~\eqref{eq:estimates-S-e}.

It remains to decompose $k^2 \Vk\psi$. We split   
   \begin{align*}
     \Vk \psi
     &
       =\Vz \psi+(\Vk-\Vz) \psi
       =  \hffg{+}\Vz \psi + \lffg{+}\Vz \psi + \mathcal{S}_{\mathcal{V}} \psi
       +\gammazint \AVtilde \psi,
   \end{align*}
   where we have used again Proposition~\ref{prop:decomposition}
     for $(\Vk-\Vz) \psi$, and $\hffg{+}$, $\lffg{+}$ are the filters
     defined in Section~\ref{sect:filters}.
   Again, $\lffg{+}\Vz \psi$ and $\gammazint \AVtilde \psi$ are 
in the stated analyticity class since $\lffg{+}\Vz \psi$ is entire (cf.\ Prop.~\ref{prop:ff_gamma1}) and $\AVtilde \psi \in \A(C k,\vartheta_\V,\Omega)$
(cf.\ Prop.~\ref{prop:decomposition}).
   For the finite regularity part,
from~\eqref{eq:ff_gamma1_2} with $s=3/2$ and~$s^\prime=1/2$, and from~\eqref{eq:estimates-S-a} with $s=s^\prime=1$,
   we get
   \begin{align*}
     \|\hffg{+}\Vz \psi\|_{1/2,\Gamma} +  \|\mathcal{S}_{\mathcal{V}} \psi  \|_{1/2,\Gamma}
     &\lesssim k^{-1}\|\Vz\psi\|_{3/2,\Gamma} +  k^{-1} \|\psi  \|_{1/2,\Gamma} \\
     &\lesssim k^{-1}\|\psi\|_{1/2,\Gamma} +  k^{-1} \|\psi  \|_{1/2,\Gamma},
   \end{align*}
   where in the last step we have used Proposition~\ref{prop:mapping properties_laplace_smooth} with
     $s=3/2$.
     
Decomposition~\eqref{eq:op_splitting_4} with the corresponding estimates follow along the same lines. 
Decomposition~\eqref{eq:op_splitting_5} with the corresponding estimates follow along the same lines as well, with the observation that the operator $ik(\nicefrac{1}{2})$ satisfies $\|ik(\nicefrac{1}{2})\|_{H^{1/2}(\Gamma)\to H^{1/2}(\Gamma)}\lesssim k$, and thus
can be absorbed in $\mathcal{S}_{6}$.
This completes the proof.
 \end{proof}

 Next, we prove a preliminary decomposition, which allows for small remainder
 terms.
 \begin{lemma}
   \label{lemma:decomposition_single_step}
  Let $\boldsymbol{\psi}:=(\psi, \psim,\psiext)=S_{k}^- (r, \Rm,\Rext)=S_k^-(\bs r)$ be the
    solution to~\eqref{eq:dual_problem}
  for $r\in L^2(\Omega)$, $\Rm\in H^{3/2}(\Gamma)$, and $\Rext\in H^{1/2}(\Gamma)$; see Corollary~\ref{cor:existence_dual}.
  Then $\Psi$ can be decomposed into finite regularity, analytic,
    and remainder terms as
  \begin{align*}
    \psi = \psiF+ \psiA + \widetilde{\psi},
    \qquad
    \psim = \psimF + \psimA+ \widetilde{\psi}_m,
    \qquad
    \psiext = \psiextF +  \psiextA + \widetilde{\psi}_{\mathrm{ext}},
  \end{align*}
  with the following properties:
  \begin{enumerate}[label=(\roman*)]
  \item 
   \label{item:lemma:decomposition_single_step-i}
The finite regularity terms satisfy,  for a constant $C>0$ independent of $k$, 
\begin{align*}
   \smalltriplenorm{}{\boldsymbol{\psi}^{\mathcal{F}}}_{k,\mathcal{V},1}
  &\leq C \smalltriplenorm{}{\bs r}_{\mathcal{V}',1}. 
  \end{align*}

  \item
   \label{item:lemma:decomposition_single_step-ii}
The analytic terms satisfy 
  \begin{align*}
    \psiA \in \aclass(M_{\Psi},\vartheta,\PP),
    \qquad \psimA \in \aclass(k M_{\Psi},\vartheta,\mathcal{O},\Gamma),
    \qquad \psiextA \in \aclass(M_{\Psi},\vartheta,\mathcal{O},\Gamma)
    \end{align*}
     with 
  $$
  M_{\Psi} \leq C k^{\beta+3}  \big(\|r\|_{0,\Omega} + \|\Rm\|_{3/2,\Gamma} + \|\Rext\|_{1/2,\Gamma} \big),
  $$
where the (unilateral) tubular neighborhood $\mathcal{O}\subset\Omega$ of $\Gamma$ and $\vartheta$ depend solely on $\Omega$, $\PP$, $\nu$, and $n$.
\item
   \label{item:lemma:decomposition_single_step-iii}
  The remainder
  $\widetilde{\boldsymbol{\psi}}=(\widetilde{\psi},\widetilde{\psi}_m,\widetilde{\psi}_{\mathrm{ext}})$ 
  solves a problem with
  modified right-hand side, namely,
  $\widetilde{\psi}=S_k^{-} \widetilde{\bs r}$
   with
     \begin{align*}
               \smalltriplenorm{}{\widetilde{\bs r}}_{\mathcal{V}^\prime,1}
     &\leq  q \smalltriplenorm{}{\bs r}_{\mathcal{V}^\prime,1}
             \end{align*}
  and $0<q<1$ independent of $k$, $r$, $\Rm$, and $\Rext$.
  \end{enumerate}
\end{lemma}
\begin{proof}
  We construct the decomposition in multiple steps.
  Using the frequency filters from Section~\ref{sect:filters},
  we start by defining
  \begin{align*}
   & \boldsymbol{\psi}_{I}^{\mathcal{F}}:=(\psiF,\widehat{\psimF},\psiextF)=S_{k}^{+}(\hff{\Omega}r,\hffg{+}\Rm,\hffg{-}\Rext), \qquad
     \boldsymbol{\psi}^{\mathcal{A}}_I:=S_{k}^{-}(\lff{\Omega}r,\lffg{+}\Rm,\lffg{-}\Rext),\\
    &\boldsymbol{\psi}^{\mathcal{F}}=(\psiF,\psimF,\psiextF):=(\psiF, \widehat{\psimF}+ik \gammazint(\psiF),\psiextF),
  \end{align*}
  for some $\eta \in (0,1)$ to be fixed later on. Recall that
  $S_{k}^{+}$ is the solution operator associated with the coercive auxiliary problem~\eqref{eq:modified_problem}, 
  and $S_{k}^-$ is the solution operator associated with the adjoint problem~\eqref{eq:dual_problem}.

  We start by noting that,  by Lemma~\ref{lemma:mod_apriori},
  $\boldsymbol{\psi}_{I}^{\mathcal{F}}$ satisfies the following estimate:
  \begin{align}
    \smalltriplenorm{}{\boldsymbol{\psi}_{I}^{\mathcal{F}}}_{k,\mathcal{V},0}
    &\lesssim k^{-3/4}\|\hff{\Omega} r\|_{\left(H^{1/4}(\Omega)\right)'} + \|\hffg{+} \Rm\|_{1/2,\Gamma}
      +\|\hffg{-} \Rext\|_{-1/2,\Gamma}  \nonumber \\
    &\lesssim \eta^{1/4}k^{-1}\big(\|r\|_{0,\Omega}
      + \|R_{m}\|_{3/2,\Gamma} + \|\Rext\|_{1/2,\Gamma}\big)
      =\eta^{1/4} k^{-1} \smalltriplenorm{}{\bs r}_{\mathcal{V'},1},
      \label{eq:phip_energy_est}
  \end{align}
  where we have used
  Proposition~\ref{prop:ff_volume} with $\tau=1/4$,
  Proposition~\ref{prop:ff_gamma1} with $s=3/2$ and $s^\prime=1/2$,
    and Proposition~\ref{prop:ff_gamma2} with $s=1/2$ and $s^\prime=-1/2$.
Note that $\hffg{-}\Rext$ has vanishing integral mean by
Proposition~\ref{prop:ff_gamma2}.

Due to Lemma~\ref{lemma:mod_aprioi_h2}, we also have that $\Psi_{I}^{\mathcal{F}}$ satisfies
  \begin{align*}
    \smalltriplenorm{}{\boldsymbol{\psi}_{I}^{\mathcal{F}}}_{k,\mathcal{V},1}
    &\lesssim
      \smalltriplenorm{}{(\hff{\Omega}r,\hffg{+}\Rm,\hffg{-}\Rext)}_{k,\mathcal{V}',1},\\
      &=
      \|\hff{\Omega} r\|_{0,\Omega} + k \|\hffg{+}  \Rm \|_{1/2,\Gamma} + \|\hffg{+} \Rm\|_{3/2,\Gamma}
      +  k \|\hffg{-} \Rext\|_{-1/2,\Gamma} + \|\hffg{-} \Rext\|_{1/2,\Gamma}
    \\
    &\lesssim \|r\|_{0,\Omega} + \|\Rm\|_{3/2,\Gamma} + \|\Rext\|_{1/2,\Gamma}
      =\smalltriplenorm{}{\bs r}_{\mathcal{V}',1}.
  \end{align*}
This, together with
$\|ik \psiF\|_{1/2,\Gamma}\lesssim k\|\psiF\|_{1,\Omega}$ and~\eqref{eq:phip_energy_est},
  implies that the components of $\Psi^{\mathcal{F}}$ satisfy 
   \ref{item:lemma:decomposition_single_step-i}. 
  
  We will show that the remainder
  \[
    \Delta=(\delta,\delta_m,\delta_{\mathrm{ext}}):=
    \boldsymbol{\psi} - \boldsymbol{\psi}^{\mathcal{F}} - \boldsymbol{\psi}^{\mathcal{A}}_I
    \]
  solves a problem structurally similar to~\eqref{eq:dual_problem},
  with a right-hand side that decomposes into a
  ``finite regularity'' part, which is characterized by good $k$-explicit estimates,
  an ``analytic part'', and a remainder that is strictly smaller than the
  original right-hand side.

  It is easy to see that $\delta$ solves the following equation (see~\eqref{dual:problem:first:equation:strong}): 
  $$
  -\fdiv(\nu \nabla \overline{\delta}) - (\kn)^2
  \overline{\delta}
  =-2 (\kn)^2 \overline{\psiF},
  \qquad \text{in }\Omega
    $$
    For the boundary condition as in~\eqref{dual:problem:first:equation:strong}, a simple calculation reveals
  \begin{align*}  
    &\nu \nabla \overline{\delta} \cdot \nGamma + \ii k
    \overline{\delta}
    +\overline{\delta_m}\\
    &\quad = (\gammanu[int]{\overline{\psi}} + \ii k \overline{\psi} )
      -(\gammanu[int]{\overline{\psiF}} + \ii k \overline{\psiF} )
      -(\gammanu[int]{\overline{\psiA}} + \ii k
      \overline{\psiA} )
+\overline{\psim} - \overline{\widehat{\psimF}} - \overline{\psimA} +ik \gammazint(\overline{\psiF})\\
    &\quad =0.
  \end{align*}
  Using a similar calculation, the coupling equation as in~\eqref{dual:problem:second:equation:strong} reads
  \begin{multline*}
    -\overline \delta + (\nicefrac{1}{2} + \Kk +i k \Vk) \overline{\delta_{\mathrm{ext}}} + \Vk \overline{\delta_m} \\
    \begin{aligned}[t]
      &= \overline \Rm - \hffg{+}\overline{\Rm}  - \lffg{+}\overline{\Rm} 
    - (\Kk - \Kz + \ii k \Vk ) \overline{\psiextF} 
    - (\Vk -\Vz)\overline{\widehat{\psimF}}+ik \Vk \gammazint(\overline{\psiF})\\
    &= -(\Kk - \Kz + \ii k \Vk ) \overline{\psiextF} 
    - (\Vk -\Vz)\overline{\widehat{\psimF}}+ik \Vk \gammazint(\overline{\psiF}).
  \end{aligned}
\end{multline*}
The operators on the right-hand side can be split by Lemma~\ref{lemma:op_decomposition_1}: 
  \begin{align*}
    -\overline \delta + (\nicefrac{1}{2} + \Kk +i k \Vk) \overline{\delta_{\mathrm{ext}}} + \Vk \overline{\delta_m}
    &=
      \mathcal{S}_{1} \overline{\psiextF} + \mathcal{S}_{2}\overline{\widehat{\psimF}}
      + {\mathcal{A}}_{1} \overline{\psiextF}
      + {\mathcal{A}}_{2} \overline{\widehat{\psimF}}
      +\mathcal{S}_{3}\gammazint(\overline{\psiF}) +{\mathcal{A}_{3}}\gammazint(\overline{\psiF}).
  \end{align*}

  The exterior problem as in~\eqref{dual:problem:third:equation:strong:3} becomes
  \begin{multline*}
   \left(\Wk + i k (\nicefrac{1}{2} - \Kprimek) - i k (\nicefrac{1}{2} + \Kk + i k\Vk)\right) \overline{\delta_{\mathrm{ext}}} -\left( (\nicefrac{1}{2} + \Kprimek) + i k \Vk\right)  \overline{\delta_m}\\
    = \left((\Wz - \Wk) + i k (\Kprimek + \Kk) - k^2 \Vk\right) \overline{\psiextF}
    +(\Kprimek - \Kprimez + i k \Vk ) \overline{\widehat{\psimF}}\\
    -\left(ik(\nicefrac{1}{2}+\Kprimek)-k^2\Vk\right) \gammazint(\overline{\psiF}).
  \end{multline*}
 The operators on the right-hand side can be split by Lemma~\ref{lemma:op_decomposition_2},
  giving 
  \begin{multline*}
    \left(\Wk + i k (\nicefrac{1}{2} - \Kprimek) - i k (\nicefrac{1}{2} + \Kk + i k\Vk)\right) \overline{\delta_{\mathrm{ext}}} - \left( (\nicefrac{1}{2} + \Kprimek) + i k \Vk\right)  \overline{\delta_m} \\
    =
    \mathcal{S}_{4}\overline{\psiextF} +     \mathcal{S}_{5} \overline{\widehat{\psimF}}
    + {\mathcal{A}}_{4} \overline{\psiextF}
    + {\mathcal{A}}_{5}  \overline{\widehat{\psimF}}
    +\mathcal{S}_{6}\gammazint(\overline{\psiF}) +{\mathcal{A}_{6}}\gammazint(\overline{\psiF}).
  \end{multline*}

  Collecting all the ``finite regularity contributions'', we get the new right-hand sides
  $\widetilde{\bs r}=(\widetilde{r},\widetilde{R}_m,\widetilde{R}_{\mathrm{ext}})$:
  \begin{align*}
    \overline{\widetilde{r}}:=-2 (\kn)^2 \overline{\psiF},
    \qquad \overline{\widetilde{R}_m}:=\SS_{1}\overline{\psiextF} +
    \SS_{2} \overline{\widehat{\psimF}}
    +\SS_{3}\gammazint(\overline{\psiF}),
    \qquad
    \overline{\widetilde{R}_{\mathrm{ext}}}:=\SS_{4}\overline{\psiextF}
    + \SS_{5} \overline{\widehat{\psimF}}
    +\SS_{6}\gammazint(\overline{\psiF}).
  \end{align*}
  From the mapping properties of the operators $\SS_i$ 
  and $k\|\psiF\|_{1/2,\Gamma}\lesssim
    k\|\psiF\|_{1,k,\Omega}$, we get
  with the properties of the high-frequency filters, 
  \begin{align*}
    \smalltriplenorm{}{\widetilde{\bs r}}_{\mathcal{V}',1}&=
    \|\widetilde{r}\|_{0,\Omega}+ \|\widetilde{R}_m\|_{3/2,\Gamma} +
    \|\widetilde{R}_{\mathrm{ext}}\|_{1/2,\Gamma} \\
    &\lesssim k^2\|\psiF\|_{0,\Omega} + k\|\psiextF\|_{1/2,\Gamma} + k
      \|\widehat{\psimF}\|_{-1/2,\Gamma}
      +k\|\psiF\|_{1/2,\Gamma}\\
    &\lesssim k\left(\|\psiF\|_{1,k,\Omega}
       +\|\psiextF\|_{1/2,\Gamma}+\|\widehat{\psimF}\|_{-1/2,\Gamma}
      \right)\\
    &\stackrel{\eqref{eq:phip_energy_est}}{\lesssim} 
      \eta^{1/4}\smalltriplenorm{}{\bs r}_{\mathcal{V}^\prime,1}.
  \end{align*}
  Thus, the by choosing $\eta$ sufficiently small to compensate for the implied
  constant, we get \ref{item:lemma:decomposition_single_step-iii} for the remainder $\widetilde{\Psi}=(\widetilde{\psi},\widetilde{\psi}_m,\widetilde{\psi}_{\mathrm{ext}})
    :=S_{k}^{-}(\widetilde{r},\widetilde{R}_m,\widetilde{R}^{\mathrm{ext}})$.

In order to prove \ref{item:lemma:decomposition_single_step-iii}
  we
  collect the remaining analytic terms as
  $$
  \boldsymbol{\psi}^{\mathcal{A}}:=\boldsymbol{\psi}^{\mathcal{A}}_I +  \boldsymbol{\psi}^{\mathcal{A}}_{II},
  $$
  where $\boldsymbol{\psi}^{\mathcal{A}}_{II}$ is defined as
  \begin{equation}
   \label{eq:lemma:decomposition_single_step-100}
   \boldsymbol{\psi}^{\mathcal{A}}_{II}:=S_{k}^{-}\left(0,\quad 
    \overline{{\A}_1(\overline{\psiextF}) +
      {\A}_{2}(\overline{\widehat{\psimF}})+\A_{3}\gammazint(\overline{\psiF})},\quad 
    \overline{{\A}_4(\overline{\psiextF}) + {\A}_{5}(\overline{\widehat{\psimF}}) +\A_{6}\gammazint(\overline{\psiF})}
    \right).
  \end{equation}
  This clearly gives that
  $\boldsymbol{\psi}^{\mathcal{A}}=\boldsymbol{\psi}-\boldsymbol{\psi}^{\mathcal{F}}-\widetilde{\Psi}$, as
  required. Next, we assert that the components of
  $\boldsymbol{\psi}^{\mathcal{A}}$ belong to the stated classes of analyticity.
For $\boldsymbol{\psi}^{\A}_{I}$, by the properties of the filter operators in
Propositions~\ref{prop:ff_volume}--\ref{prop:ff_gamma2}, we have 
  \begin{align*}
  \lff{\Omega}r &\in \aclass(C\|r\|_{0,\Omega},\vartheta_\eta 
  , \PP),
  \;\; \qquad \qquad \qquad \lffg{+}\Rm \in \aclass(C \|\Rm\|_{1/2,\Gamma}, 
  \vartheta_\eta, 
  \R^d), \\
  \lffg{-}\Rext &\in \aclass(C
  k^{d/2}\|\Rext\|_{-1/2,\Gamma}, \vartheta_\eta 
  , \mathcal{O},\Gamma),
  \end{align*}
for a suitable (unilateral) tubular neighborhood $\mathcal{O} \subset \Omega$ of $\Gamma$.
This, via Lemma~\ref{lemma:analytic_regularity}, gives that $\Psi^{\mathcal{A}}_{I}=:\left((\psiA)_{I},(\psimA)_I,(\psiextA)_I\right)$ satisfies
  $$
  (\psiA)_I \in \aclass(M_{I},\vartheta,\Omega),
  \quad (\psimA)_I \in \aclass( k  M_{I},\vartheta,\mathcal{O},\Gamma),
  \quad (\psiextA)_I \in \aclass(M_{I},\vartheta,\mathcal{O}, \Gamma)
  $$
  for a suitable $\vartheta>0$ independent of $k$, a possibly adjusted tubular neighborhood $\mathcal{O}\subset\Omega$ of $\Gamma$, 
  and
  \[
    \begin{split}
    M_{I} &\leq C \max(1,k^{\beta})\big(\|r\|_{0,\Omega} +
  \|\Rm\|_{1/2,\Gamma} 
  + k^{d/2}\|\Rext\|_{-1/2,\Gamma}\big)\\
  &\le C \max(1,k^{\beta}) k^{d/2}\big(\|r\|_{0,\Omega} +
  \|\Rm\|_{1/2,\Gamma} 
  + \|\Rext\|_{-1/2,\Gamma}\big).
  \end{split}
  \]
  For $\Psi^{\A}_{II}$, we note that, by Lemmas~\ref{lemma:op_decomposition_1} and \ref{lemma:op_decomposition_2},
  the three arguments of $S_{k}^{-}$ in (\ref{eq:lemma:decomposition_single_step-100}) are analytic with leading constant $\bigO(k^3)$. From
    Lemma~\ref{lemma:analytic_regularity} and from \ref{item:lemma:decomposition_single_step-i}, we infer that
    the components of $\Psi^{\A}_{II}$ belong to analyticity classes
    with leading constant
    \[
M_{II}\le C \max(1,k^{\beta})k^3\big[\|r\|_{0,\Omega} +
  \|\Rm\|_{3/2,\Gamma} 
  +\|\Rext\|_{1/2,\Gamma}\big].
\]
Since $d\leq 3$,
we have that
  $\max(M_{I},M_{II})=M_{II}=:M_{\Psi}$. Therefore, the components
  of $\Psi^{\A}=\Psi^{\A}_{I}+\Psi^{\A}_{II}$
  satisfy~\ref{item:lemma:decomposition_single_step-ii}. This completes the proof.
\end{proof}

Finally, by iterating the argument of Lemma~\ref{lemma:decomposition_single_step}, we can eliminate the remainder term and prove the main theorem of this section.
\begin{theorem}[Regularity splitting]
  \label{thm:regularity_splitting}
  Let $\boldsymbol{\psi}:=(\psi,\psim,\psiext)$ solve the adjoint problem~\eqref{eq:dual_problem}
  for right-hand sides $\bs r=(r,\Rm,\Rext)$ satisfying
  $$
  r \in L^2(\Omega), \quad \Rm\in H^{3/2}(\Gamma), \quad \Rext \in H^{1/2}(\Gamma).
  $$
  Then, we can decompose $\boldsymbol{\psi}$ as 
  \begin{align*}
    \psi=\psiF + \psiA, \quad
    \psim=\psimF + \psimA, \quad
    \psiext=\psiextF + \psiextA,
  \end{align*}
  with the following properties:
  \begin{enumerate}[label=(\roman*)]
   \item
\label{item:thm:regularity_splitting-i}
The finite regularity terms satisfy, for a constant $C>0$ independent of $k$,
\begin{align*}
  \smalltriplenorm{}{\boldsymbol{\psi}^{\mathcal{F}}}_{k,\mathcal{V},1} 
  &\leq
    C   \smalltriplenorm{}{\bs r}_{\mathcal{V}',1}. 
  \end{align*}
  \item
\label{item:thm:regularity_splitting-ii}
The analytic terms satisfy
  \begin{align*}
    \psiA\in \aclass\big(M,\vartheta,\PP), \quad
    \psimA \in \aclass\big( k M,\vartheta,\mathcal{O},\Gamma),\quad
    \psiextA \in \aclass\big(M,\vartheta,\mathcal{O},\Gamma)
  \end{align*}
  with $$
  M \leq C k^{\beta+3}
  \big(\|r\|_{0,\Omega}  +\|\Rm\|_{3/2,\Gamma} 
  + \|\Rext\|_{1/2,\Gamma} 
  \big),
  $$
where the tubular neighborhood $\mathcal{O}\subset\Gamma$ of $\Gamma$ and $\vartheta$ depend solely on $\Omega$, $\PP$, $\nu$, and $n$.
  \end{enumerate}
\end{theorem}
\begin{proof}
  We apply iteratively Lemma~\ref{lemma:decomposition_single_step}.
  We start by writing $\bs r^{(0)}=(r^{(0)},\Rm^{(0)},\Rext^{(0)}):=(r,\Rm,\Rext)$
  and decompose the
    solution $\bs \psi$ according to Lemma~\ref{lemma:decomposition_single_step}:
  \begin{align*}
    \boldsymbol \psi=\boldsymbol \psi^{\mathcal{F}}_{(0)} + \boldsymbol \psi^{\A}_{(0)} + \boldsymbol{\widetilde{\psi}} _{(0)}.
  \end{align*}
  Since the remainder $\boldsymbol{\widetilde{\psi}} _{(0)}$ solves the same
  problem with a new right-hand side
  $
  \bs r^{(1)}:=({r}^{(1)},{R}^{(1)}_m,{R}^{(1)}_{\mathrm{ext}})
  $,
  namely, $\bs{\widetilde{\psi}} _{(0)}=S^-_k(\bs r^{(1)})$,
  we can again split this as
  \begin{align*}
    \boldsymbol {\widetilde{\psi}} _{(0)}=\boldsymbol \psi^{\mathcal{F}}_{(1)} + \boldsymbol \psi^{\A}_{(1)} + \boldsymbol{ \widetilde{\psi}}_{(1)}
  \end{align*}
  and so on, defining sequences $\boldsymbol\psi^{\mathcal F}_{(\ell)}$, $\boldsymbol \psi^{\A}_{(\ell)}$, $\boldsymbol{\widetilde\psi}_{(\ell)}$
of functions by
  \begin{align*}
    \boldsymbol{\widetilde{\psi}} _{(\ell-1)} 
    = \boldsymbol \psi^{\mathcal{F}}_{(\ell)} + \boldsymbol\psi^{\A}_{(\ell)} + \boldsymbol{\widetilde{\psi}} _{(\ell)}
    \qquad \text{with} \qquad \boldsymbol{\widetilde{\psi}} _{(\ell)} =
    S_{k}^- \big(\boldsymbol r^{(\ell+1}\big), \qquad \ell\ge 1.
  \end{align*}
  From Lemma~\ref{lemma:decomposition_single_step}\ref{item:lemma:decomposition_single_step-iii}, for any $\ell\ge 1$,
  the right-hand sides $\bs r^{(\ell)}$ 
  satisfy
  \begin{align*}
    \smalltriplenorm{}{\bs r^{(\ell)}}_{\mathcal{V}^\prime,1}
    \leq q \smalltriplenorm{}{\bs r^{(\ell-1)}}_{\mathcal{V}^\prime,1}
    \leq q^{\ell} \smalltriplenorm{}{\bs r^{(0)}}_{\mathcal{V}^\prime,1}
  \end{align*}
  Since $0<q<1$, we get that these right-hand sides converge to zero as $\ell \to \infty$.
  This implies that we can write
  \begin{align*}
    \boldsymbol \psi:=\boldsymbol \psi^{\mathcal{F}} + \boldsymbol \psi^\A, \qquad \text{with} \quad
    \boldsymbol \psi^{\mathcal{F}}:=\sum_{\ell=0}^{\infty}{\boldsymbol \psi^{\F}_{(\ell)}} \quad\text{and}\quad
    \boldsymbol \psi^\A:=\sum_{\ell=0}^{\infty}{\boldsymbol \psi^\A_{(\ell)}}.
  \end{align*}
  For the finite regularity part $\boldsymbol \psi^{\mathcal{F}}=:(\psiF,\psimF,\psiextF)$,
  we get the estimate
  \begin{align*}
    \smalltriplenorm{}{\boldsymbol \psi^{\mathcal{F}}}_{k,\mathcal{V},1}
    &\lesssim
      \sum_{\ell=0}^{\infty}{
      \smalltriplenorm{}{\boldsymbol \psi^{\mathcal{F}}_{(\ell)}}_{k,\mathcal{V},1}}
      \lesssim
          \sum_{\ell=0}^{\infty}{
      \smalltriplenorm{}{\bs r^{(\ell)}}_{\mathcal{V}',1}}
    \lesssim
      \smalltriplenorm{}{\bs r^{(0)}}_{\mathcal{V}',1}
      \big( \sum_{\ell=0}^{\infty}{ q^{\ell}} \big). 
%
  \end{align*}
  Since the geometric series converges, we get the stated estimate for
  the finite regularity terms.

  For the analytic part, we focus on the interior contribution
  and denote it as $\psi^{\A}$.
  For any $\ell\ge 0$, we get  
  $
  \psi^{\A}_{(\ell)} \in \aclass(M_{(\ell)},\vartheta,\Omega) 
  $
(note that $\vartheta$ is independent of the right-hand side
$\bs r^{(\ell)}$ 
and is thus the same in each iteration), with the leading constant
  \begin{align*}
    M_{(\ell)}
    &\lesssim k^{\beta+3}
    \smalltriplenorm{}{\bs r^{(\ell)}}_{\mathcal{V}',1}
      \lesssim k^{\beta+3} q^{\ell} \smalltriplenorm{}{\bs r^{(0)}}_{\mathcal{V}',1}.
%
  \end{align*}
  We can therefore again argue as before that the geometric series converges, and we
  get
  \begin{align*}
    \|\nabla^p \psi^{\A}\|_{0,\Omega \setminus \partial \PP}
    &\lesssim \sum_{\ell=0}^{\infty}{\|\nabla^p\psi^{\A}_{(\ell)}\|_{0,\Omega \setminus \partial \PP}} 
    \lesssim k^{\beta+3} \sum_{\ell=0}^{\infty}{M_{(\ell)} \vartheta^p \max(p+1,k)^p} \\
    &\lesssim k^{\beta+3} \vartheta^p \max(p+1,k)^p  
      \smalltriplenorm{}{\bs r^{(0)}}_{\mathcal{V}',1}
    \sum_{\ell=0}^{\infty}{ q^\ell}.
  \end{align*}
  In other words $\psi^{\A} \in \aclass\big(M,\vartheta,\Omega\big)$
  with
  $$M:=\frac{k^{\beta+3}}{1-q}\Big( \|r\|_{0,\Omega} + \|\Rm\|_{3/2,\Gamma} + \|\Rext\|_{1/2,\Gamma}\Big).$$
  The result for the other contributions $\psim^\A$ and $\psiext^\A$ follows along the same lines.
\end{proof}

\section{Convergence of the discrete schemes}
\label{sec:convergence}
With the results of the previous two sections in place, we can prove that both the
conforming FEM and the DG method lead to $k$-independent quasi-optimal error estimates 
provided a weak ($k$-dependent) resolution condition is satisfied.
For $s \in \N_0$, $M>0$, and a fixed tubular neighborhood $\mathcal{O}$ of $\Gamma$
sufficiently enough so that $\nu$, $n$ are smooth on $\mathcal{O}$ and $\nGamma$ has a smooth extension
  to $\mathcal{O}$, we also introduce the following product spaces:
\begin{align*}
  \mathcal{V}^{s}&:=H^{1+s}_\PP(\Omega) \times H^{-1/2+s}(\Gamma) \times H^{1/2+s}(\Gamma),\\
  \mathcal{A}^{M}&:=\aclass(M, \vartheta_u,\PP) \times
  \aclass(k M, \vartheta_m,\mathcal{O},\Gamma) \times
  \aclass(M, \vartheta_{\uext},\mathcal{O},\Gamma).
\end{align*}
The space $\mathcal{V}^{s}$ is equipped with the norm 
$\|\cdot\|_{k,\mathcal{V},s}$.
Most arguments of the convergence analysis will be presented in detail only for the DG case,
as the proofs for the conforming FEM would be essentially the
same, just with 
less technicalities.

\subsection{Discontinuous Galerkin method}
For $s \in \N$, we define the approximation quantities
\begin{align} \label{eq:etaDG}
  \eta^{(s)}_{DG}:= \adjustlimits\sup_{0\ne \x \in \mathcal{V}^{s}}\inf_{\xh \in \mathcal{V}^{DG}_h} \frac{\energynormp{\x-\xh}}{
  \smalltriplenorm{}{\x}_{k,\mathcal{V},{s}}},
  \qquad \eta^{(\A)}_{DG}:=\adjustlimits\sup_{\x \in \mathcal{A}^{1}}
  \inf_{\xh \in \mathcal{V}^{DG}_h}
  \energynormp{\x - \xh}.
\end{align}
In our final convergence result, we will  assume these quantities to be sufficiently small.
We remark that
\begin{equation}\label{eq:AM}
\forall \x\in \mathcal{A}^{M} \colon \qquad 
\inf_{\xh \in \mathcal{V}^{DG}_h}
\energynormp{\x - \xh}\le M
\eta^{(\A)}_{DG}.
\end{equation} 
Our convergence analysis hinges on duality arguments. These duality arguments are possible for the chosen DG-discretization
since it is designed to be ``adjoint consistent''.
\begin{lemma}[adjoint consistency] 
\label{lemma:adjoint-consistency}
Given $(r,\Rm, \Rext) \in L^2(\Omega) \times H^{1/2}(\Gamma) \times H^{-1/2}(\Gamma)$, let $\bs \Psi$ be the solution to the adjoint 
problem (\ref{eq:dual_problem}). Then, for any $t > 3/2$ and for all $(\Phi,m,\uext) \in H^{t}_{pw}(\taun) \times H^{-1/2}(\Gamma) \times H^{1/2}(\Gamma)$,
we have
\begin{equation}
\label{eq:lemma:adjoint-consistency}
\TTDG{k}((\Phi,m,\uext),\bs \Psi) = (\Phi,r)_{0,\Omega} + \langle m,\Rm\rangle_\Gamma + \langle \uext,\Rext\rangle_\Gamma = \T((\Phi,m,\uext),\bs\Psi). 
\end{equation}
\end{lemma}
\begin{proof}
The result follows as in the proof of \cite[Prop.~{4}]{MMPR20}.
\end{proof}
We have the following abstract quasi-optimality result for the DG-FEM: 
\begin{theorem}
  \label{thm:dg_is_quasioptimal}
  Let the solution $(u,m,\uext)$ to~\eqref{system:Helmholtz} be in
  $H^{\frac{3}{2} +t} (\taun)\times  L^2(\Gamma) \times H^{\frac{1}{2}}(\Gamma)$
  for some $t>0$,
  and let $(\uh,\mh, \uhext)\in
  \mathcal{V}^{DG}_h$ be the solution of
  method~\eqref{dgBEM:short-version} with penalty parameters $\mathfrak{a} \ge\mathfrak{a}_0$, $\mathfrak{b}\ge 0$, 
$0<\mathfrak{d} \leq \mathfrak{d}_0$ (see~\eqref{dG-parameters} and Lemma~\ref{lemma:dg_garding}).    
  Let $\mu_{stab}  \in [0,4]$ be as in
  Corollary~\ref{cor:TDG_is_bounded_with_k} and
  $\mu_{\A}:=\beta+3$, with $\beta$ as in Assumption~\ref{ass:existence_global}.
There is a positive $q < 1$ independent of $k \ge k_0$ such that, under the resolution condition
  \begin{equation}
\label{eq:scale-resolution-condition-abstract}
  k \eta^{(1)}_{DG} + k^{1+\mu_{stab}+\mu_{\A}} \eta^{(\A)}_{DG} \leq q, 
  \end{equation}
  the following estimate holds true: 
  \begin{multline*}
 \Vert u-\uh \Vert_{\dGnorm} + \Vert m-\mh\Vert_{-1/2,\Gamma} +\Vert \uext-\uhext\Vert_{1/2, \Gamma}\\
\lesssim \Vert u-\vh \Vert_{\dGnormp}+\Vert m-\lambdah\Vert_{-1/2,\Gamma} +\Vert \uext-\vhext\Vert_{1/2, \Gamma} +  \Vert \msf^{1/2} p^{-1} (m-\lambdah)\Vert_{0,\Gamma}
\end{multline*}
for any $(\vh,\lambdah,\vhext) \in \mathcal{V}^{DG}_h$, with
hidden constant independent of $k$, $h$, and $p$.
\end{theorem}
\begin{proof}
  To shorten the expressions, we write $\x:=(u,m,\uext)$ and $\xh=(\uh,\mh,\uhext)$.
  For arbitrary $\yh \in \mathcal{V}_h$, using the
  G{\aa}rding inequality from Lemma~\ref{lemma:dg_garding},
  Galerkin orthogonality, and the boundedness from
  Lemma~\ref{lemma:dg_explicit_bound}, we obtain
  \begin{align}
    \energynorm{\xh - \yh}^2
    &= (\Re + \varepsilon \Im) \big[\T_{DG}(\xh - \yh, \xh - \yh) + \langle \xh -\yh, \Theta(\xh-\yh)\rangle  \big] \nonumber\\
    &= \begin{multlined}[t][11cm]
      (\Re + \varepsilon \Im)\big[\T_{DG}(\x - \yh, \xh - \yh) + \langle \x -\yh, \Theta(\xh-\yh)\rangle \big]  \\
      + (\Re + \varepsilon \Im)\langle \xh - \x, \Theta(\xh-\yh)\rangle  \nonumber
    \end{multlined}
    \\
    &\stackrel{\mathclap{\eqref{eq:k_expl_continuities:3}}}{\lesssim} \;\,\energynormp{\x - \yh} \energynorm{\xh-\yh} + |\langle \x -\xh, \Theta(\xh-\yh)\rangle  |.
      \label{eq:energy_estimate1}
  \end{align}
  We now focus on the last term and split $\Theta$ into two parts according to Lemma~\ref{lemma:introducing_theta}:
  \begin{align}
    \label{eq:split_the_theta_terms}
    \langle \x -\xh, \Theta(\xh-\yh)\rangle 
    &=\langle \x -\xh, \Theta^\F(\xh-\yh)  \rangle + \langle \x -\xh,
      \Theta^\A(\xh-\yh)\rangle=: T_1+T_2.
  \end{align}
  For the 
  term $T_1$, we define the adjoint solution $\bs \Psi^1$ by 
$$
  \TTC{k}(\bs \Phi,\bs \Psi^1)=\langle \bs \Phi, \Theta^\F(\xh-\yh) \rangle
  \qquad \forall \bs \Phi \in \mathcal{V} = H^{1}(\Omega) \times H^{-1/2}(\Gamma) \times H^{1/2}(\Gamma).
$$
By adjoint consistency (Lemma~\ref{lemma:adjoint-consistency}), we have 
  $$
  \TTDG{k}(\bs \Phi,\bs \Psi^1)=\langle \bs \Phi, \Theta^\F(\xh-\yh) \rangle
  \qquad \forall \bs \Phi \in H^{3/2+t}_{pw}(\taun) \times H^{-1/2}(\Gamma) \times H^{1/2}(\Gamma).
  $$
  We take $\bs \Phi=\x -\xh$ and, using Galerkin orthogonality, we
  can write, for arbitrary $\bs \Psi_h^1 \in \mathcal{V}^{DG}_h$, 
\begin{align}
  T_1&=  \langle \x -\xh, \Theta^\F(\xh-\yh)  \rangle
    =\TTDG{k}(\x-\xh,\bs \Psi^1) =  \TTDG{k}(\x-\xh,\bs \Psi^1 - \bs \Psi_h^1) \nonumber\\
&=\big[\TTDG{k}(\x-\xh,\bs \Psi^1 - \bs \Psi_h^1)
+\langle \x-\xh,\Theta^\A(\bs \Psi^1 - \bs \Psi_h^1)\rangle\big]
-\langle \x-\xh,\Theta^\A(\bs \Psi^1 -\bs \Psi_h^1)\rangle.
      \label{eq:dG:est_of_thetaf}
  \end{align}
  For the first term, we can  use the boundedness of $\TTDG{k}$ given in (\ref{eq:k_expl_continuities:2})
  from Lemma~\ref{lemma:dg_explicit_bound} 
  to get
  \begin{multline*}
\big[\TTDG{k}(\x-\xh,\bs \Psi^1 - \bs \Psi_h^1)
+\langle \x-\xh,\Theta^{\A}(\bs \Psi^1 - \bs \Psi_h^1)\rangle\big]
    \\ \lesssim  
\left(
\bigtriplenorm{\big}{\x-\yh}\bigtriplenorm{\big}_{\dGnormp}
+
\bigtriplenorm{\big}{\xh-\yh}\bigtriplenorm{\big}_{\dGnorm}
\right)
\energynormp{\bs \Psi^1 - \bs \Psi_h^1}.
  \end{multline*}
  
  We need to understand the best approximation of $\bs\Psi^1$. Using
  the splitting
$\bs \Psi^1=(\bs \Psi^1)^\F+(\bs \Psi^1)^\A$
  from Theorem~\ref{thm:regularity_splitting}, we get
  \begin{align}
\nonumber 
    \energynormp{\bs \Psi^{1} -\bs \Psi_h^{1}}
    &\stackrel{\!\!\mathclap{\eqref{eq:etaDG},\eqref{eq:AM}}}{\lesssim}  \quad  
      \eta^{(1)}_{DG}
      \|(\bs \Psi^{1})^\F\|_{k,\mathcal{V},1}
      +\eta^{(\A)}_{DG}  k^{\beta+3}
      \|\Theta^\F(\xh-\yh)\|_{\mathcal{V}',1}
      \\
\nonumber 
    &\lesssim \quad \big(\eta^{(1)}_{DG}+ \eta^{(\A)}_{DG}k^{\beta+3}\big)\|\Theta^\F(\xh-\yh)\|_{\mathcal{V}',1}\\
    & 
       \lesssim
      \quad k\big(\eta^{(1)}_{DG}+ \eta^{(\A)}_{DG}k^{\beta+3}\big)
      \energynorm{\xh-\yh},
\label{eq:thm.7.1-10}
  \end{align}
  where, in the second inequality, we have used Theorem~\ref{thm:regularity_splitting}\ref{item:thm:regularity_splitting-i},
    and in the last inequality, we have used
    Lemma~\ref{lemma:introducing_theta}. Inserting the previous two 
    estimates
into~\eqref{eq:dG:est_of_thetaf} and taking into account the
    definition of $T_2$ in~\eqref{eq:split_the_theta_terms}, we get 
  \begin{equation}\label{eq:T1T2}
  \begin{split}
    T_1+T_2\lesssim &\ 
      k\big(\eta^{(1)}_{DG}+ \eta^{(\A)}_{DG}k^{\beta+3}\big)
      \big[\bigtriplenorm{\big}{\x-\yh}\bigtriplenorm{\big}_{\dGnormp}
      +\bigtriplenorm{\big}{\xh-\yh}\bigtriplenorm{\big}_{\dGnorm}
      \big]\energynorm{\xh-\yh}\\
  &  +\big|\big\langle \x -\xh, \Theta^\A(\xh-\yh - (\bs \Psi^1 - \bs \Psi_h^1)) \big\rangle\big|.
    \end{split}
    \end{equation}
  We still need to estimate the last term on the right-hand side of~\eqref{eq:T1T2}, which we again
  represent via a dual solution $\bs \Psi^2$:
  \begin{align*}
    \big\langle \bs v , \Theta^\A(\xh-\yh -( \bs \Psi^1 - \bs \Psi_h^1)) \big \rangle
    =: \T(\bs v, \bs \Psi^{2}) 
\qquad \forall \bs v \in \mathcal{V}.
  \end{align*}
By adjoint consistency (Lemma~\ref{lemma:adjoint-consistency}), we have in particular
  \begin{align*}
    \big\langle \x -\xh, \Theta^\A(\xh-\yh -( \bs \Psi^1 - \bs \Psi_h^1)) \big \rangle
    = \TTDG{k}(\x -\xh, \bs \Psi^{2}).
  \end{align*}
  To analyze $\bs \psi^2$, we start by observing that $\Theta^{\A}$ maps
  into the space of analytic functions. 
    Namely, we have by Lemma~\ref{lemma:introducing_theta} that 
    $$
    \Theta^\A\big(\xh-\yh -( \bs \Psi^1 -\bs\Psi_h^1)\big) \in \mathcal{A}^{M_1}
    \qquad \text{with} \qquad 
    M_1\lesssim k^{3} 
    \energynorm{\xh-\yh -( \bs\Psi^1 -\bs\Psi_h^1)},$$
    where we used that the $\energynorm{\cdot}$ norm is stronger than
    all the norms required in Lemma~\ref{lemma:introducing_theta}.
    The analytic regularity result from Lemma~\ref{lemma:analytic_regularity} yields 
    with $\mu_{\A}:=\beta+3$ 
   \begin{align}
     \bs\Psi^2 \in \mathcal{A}^{M_{\bs \Psi_2}}
   \qquad \text{with}
   \qquad
   M_{\bs \Psi_2}\lesssim
   k^{\beta} M_{1} 
   \lesssim k^{\mu_{\A}}\energynorm{\xh-\yh -( \bs\Psi^1 -\bs\Psi_h^1)}.
   \label{eq:Psi2_is_analytic}
   \end{align}
  To estimate the last term in \ref{eq:T1T2}, we use again
  the Galerkin orthogonality and get 
  \begin{align*}
    \langle \x -\xh, \Theta^\A(\xh-\yh + \bs\Psi^1 -\bs\Psi_h^1) \rangle
    = \TTDG{k}(\x -\xh, \bs\Psi^{2}-\bs\Psi^{2}_h)
  \end{align*}
  for any $\bs\Psi^{2}_h \in \mathcal{V}_h$.
  Using the analyticity of $\bs \Psi$, and
  the bound from Corollary~\ref{cor:TDG_is_bounded_with_k},
  we get 
  \begin{align*}
   & |\langle \x -\xh, \Theta^\A(\xh-\yh + \bs\Psi^1 -\bs\Psi_h^1) \rangle|
    = |\T_{DG}(\x -\xh, \bs\Psi^{2}-\bs\Psi^{2}_h)| \\
&\qquad\lesssim k^{\mu_{stab}}\big(\energynormp{\x -\yh} 
    +\energynorm{\xh -\yh}
    \big)\bigtriplenorm{}{\bs\Psi^{2}-\bs\Psi^{2}_h}\bigtriplenorm{}_{\dGnormp} \\
&\qquad\lesssim k^{\mu_{stab}}
\big(\energynormp{\x -\yh} 
    +\energynorm{\xh -\yh}
    \big)
    k^{\mu_{\A}}\eta^{(\A)}_{DG}
    \bigtriplenorm{}{\xh-\yh + \bs \Psi^1 -\bs\Psi_h^1}\bigtriplenorm{}_{\dGnorm} \\
&\qquad\lesssim
    k^{\mu_{stab}+\mu_\A}\eta^{(\A)}_{DG}
    \big[
    k\big(\eta^{(1)}_{DG}+\eta^{(\A)}_{DG}
   k^{\beta+3}
    \big)+1
    \big]\\
    &\qquad\qquad\qquad\qquad\qquad\qquad\qquad \cdot
    \big(\energynormp{\x -\yh} 
    +\energynorm{\xh -\yh}
    \big)
    \energynorm{\xh -\yh}.
  \end{align*}
In the penultimate inequality, 
we have used~\eqref{eq:etaDG} and~\eqref{eq:Psi2_is_analytic}, and in the last inequality, we have used the estimate 
(\ref{eq:thm.7.1-10}) for $\energynormp{\bs\Psi^{1} - \bs\Psi_h^{1}}$ derived above.

Inserting the previous estimate into~\eqref{eq:T1T2}, and taking
into account~\eqref{eq:split_the_theta_terms}
and~\eqref{eq:energy_estimate1}, we obtain
  \begin{align*}
    \energynorm{\xh - \yh}^2
    \lesssim&
      \energynormp{\x - \yh} \energynorm{\xh-\yh} \\
      &+
A\big(\energynormp{\x - \yh}+
\energynorm{\xh-\yh}\big)
\energynorm{\xh-\yh}\\
&+B(A+1)\big(\energynormp{\x - \yh}+
\energynorm{\xh-\yh}\big)
\energynorm{\xh-\yh},
  \end{align*}
  where we have set
  \[
  A:=k\big(\eta^{1)}_{DG}+\eta^{(\A)}_{DG}
    k^{\beta+3}
    \big),\qquad
  B:=  k^{\mu_{stab}+\mu_A} \eta^{(\A)}_{DG}.
  \]
We divide both sides of the above inequality by $\energynorm{\xh-\yh}$  
and collect the terms. Writing $C$ for the implied constant, we get
\[
\big(1-CA-CB(A+1)\big)\energynorm{\xh-\yh}\leq C(1+A+B(A+1))
\energynormp{\x - \yh}.
\]
Selecting $q$ in the statement of the theorem sufficiently small implies that the terms 
$1 - CA - CB (A+1)$ and 
$1 + A + B (A+1)$ are close to $1$.
This, together with the triangle inequality, concludes the proof:
\[
\energynorm{\x-\xh}\le
\energynorm{\x-\yh}+\energynorm{\xh-\yh}
\lesssim 
\energynormp{\x - \yh}. \qedhere
\]
%
\end{proof}

For the specific case of regular meshes whose element maps satisfy the conditions of Assumption~\ref{def:element-maps}, 
we arrive at the following corollary: 
\begin{corollary}
\label{cor:dg_explicit_convergence}
Let the element maps satisfy Assumption~\ref{def:element-maps}.
Let $f \in L^2(\Omega)$ with $f|_{P_j} \in H^{\ell-1}(P_j)$ for all $P_j \in \PP$, $j =0,\ldots,L$.  
Let $(u,m,\uext) \in \mathcal{V}$ be the solution to~\eqref{system:Helmholtz},
  and let $(\uh,\mh, \uhext)\in
  \mathcal{V}^{DG}_h$ be the solution of
  method~\eqref{dgBEM:short-version} with flux parameters defined
  in~\eqref{dG-parameters} as in Theorem~\ref{thm:dg_is_quasioptimal}. 
Then: given $c_2 > 0$, there is $c_1>0$
independent of $h$, $p$, $k$ such that, under the scale resolution condition
  \begin{align}
\label{eq:scale-resolution-condition-hp}
     \frac{ k h}{p} &\leq c_1  \qquad \text{and} \qquad p \geq \max (1,c_2 \ln k ), 
  \end{align}
the resolution condition (\ref{eq:scale-resolution-condition-abstract}) is satisfied so that the quasi-optimality result
of Theorem~\ref{thm:dg_is_quasioptimal} holds true. Furthermore, there exists a constant $\sigma >0$ independent of $h$, $p$, $k$ such that
  \begin{align*}
 \Vert u-\uh \Vert_{\dGnorm} + \Vert m-\mh\Vert_{-1/2,\Gamma} +\Vert \uext-\uhext\Vert_{1/2, \Gamma}
 &\lesssim \Bigg( \Big(\frac{h}{p}\Big)^{\ell}
 + k^{\beta_0+1}\Big(\Bigl(\frac{h}{h+\sigma}\Bigr)^p + k\Big(\frac{ k h }{\sigma p}\Big)^p\Bigg),
  \end{align*}
with hidden constant independent of $k$, $h$, and $p$.
  Here $\beta_0 \geq 0$ is the constant in Assumption~\ref{ass:existence_global}.
\end{corollary}

\begin{proof}
  By Theorem~\ref{thm:dg_is_quasioptimal},
  it is key to estimate the approximation quantities $\eta^{(1)}_{DG}$ and $\eta^{(\A)}_{DG}$.
  For the interior DG-contributions, we use the arguments of \cite[Thm.~{4.11}]{MPS13}. We 
  point out that we explicitly introduce the quantities $\eta^{(1)}_{DG}$ and $\eta^{(\A)}_{DG}$, whereas
  \cite{MPS13} estimates these two terms implicitly in \cite[Thm.~{4.11}]{MPS13}. 

\emph{Step 1 (volume contributions):}
  There exists $\sigma>0$, independent of $h$ and $p$, such that 
  the following estimates are valid for $1\leq \ell\leq p$ (under the assumption that $kh/p$ is bounded): 
  \begin{subequations}
    \label{eq:approx_int}
  \begin{align}
    \label{eq:approx_int-a}
    \inf_{\vh \in \Vh^{C}} \|u - \vh\|_{\dGnormp}&\leq C_\ell  \Big(\frac{h}{p}\Big)^{\ell}  \|u\|_{{1+\ell},\PP}
                            \qquad \text{for $u \in H^{1+\ell}_{\PP}(\Omega)$},\\
    \label{eq:approx_int-b}
    \inf_{\vh \in \Vh^{C}}\|u -\vh \|_{\dGnormp} &\lesssim 
                              k^{\mu}\Big[\Big(\frac{h}{h+\sigma}\Big)^p + k\Big(\frac{ k h }{\sigma p}\Big)^p\Big]
                             \qquad \text{for $u \in \aclass\big(C _uk^{\mu}, \vartheta,\PP)$.}                            
  \end{align}
\end{subequations}
Estimate (\ref{eq:approx_int-a}) is a consequence of the proof of \cite[Thm.~{4.11}]{MPS13}. Key is \cite[Cor.~{7.4}]{MPS13},
which constructs an $H^1(\Omega)$-conforming approximation in an element-by-element fashion that is optimal in a broken $H^2$-norm. 
Details of the present generalization to piecewise $H^{\ell+1}$-functions are given in Lemma~\ref{lemma:broken-H2-approximation} in Appendix~\ref{app:apprDG}. 

Estimate (\ref{eq:approx_int-b}) is taken from \cite[Thm.~{4.11}]{MPS13}, and follows with the arguments presented in the proof of 
\cite[Thm.~{5.5}]{MS10} (cf.\ in particular the formula above \cite[(5.9)]{MS10}).

\emph{Step 2 (boundary contributions):}
For the boundary functions, we recall from \cite[App.~{C}]{dgfembem}
that the following estimates are valid:
  \begin{align*}
    \inf_{\lambdah \in \Wh}\|m - \lambdah\|_{-1/2,\Gamma} + \|\msf^{1/2} p^{-1}(m - \lambdah)\|_{0,\Gamma}
    &\lesssim  \
      \Big(\frac{h}{p}\Big)^{\ell} \|m\|_{{-1/2+\ell},\Gamma}
                            \qquad \text{for $m \in H^{-1/2+\ell}(\Gamma)$,}\\
   \inf_{\vhext \in \Zh} \|\uext - \vhext\|_{1/2,\Gamma}&\lesssim  \Big(\frac{h}{p}\Big)^{\ell} \|\uext\|_{1/2+\ell,\Gamma}
                            \qquad \text{for $\uext \in H^{1/2+\ell}(\Gamma)$.}
\end{align*}
For a fixed tubular neighborhood $\mathcal{O}$ of $\Gamma$ and $h$ so small that the elements of $\taun$ touching~$\Gamma$ are 
contained in this neighborhood, by approximating on these element using the results of \cite[proof of Thm.~{5.5}]{MS10} 
and appropriate trace (or multiplicative trace) estimates, and additionally estimating generously 
$\|m - \lambdah\|_{-1/2,\Gamma} \lesssim \|m - \lambdah\|_{0,\Gamma}$, we obtain  
  \begin{align*}
    \inf_{\lambdah \in \Wh}\|m -\lambdah \|_{-1/2,\Gamma} + \|\msf^{1/2} p^{-1}(m - \lambdah)\|_{0,\Gamma}
    &\lesssim
      k C_m\Bigg[\Big(\frac{h}{h+\sigma}\Big)^{p+1/2} + k^{1/2}\Big(\frac{ k h }{\sigma p}\Big)^{p+1/2} \Bigg], \\
    \inf_{\vhext \in \Zh} \|\uext -\vhext \|_{1/2,\Gamma} &\lesssim
                                       C_{\uext}\Bigg[\Big(\frac{h}{h+\sigma}\Big)^p + k\Big(\frac{ k h }{\sigma p}\Big)^p\Bigg]
  \end{align*}
for $m \in \aclass\big( k C_m , \vartheta,\mathcal{O},\Gamma)$ and 
                             $\uext \in \aclass\big(C_{\uext} , \vartheta,\mathcal{O},\Gamma)$.

\emph{Step 3 (estimating $\eta^{(\ell)}_{DG}$, $\eta^{(\A)}_{DG}$):}
  From the estimates in Steps~1 and 2, we get (with a possibly modified $\sigma$)
  \begin{align}
  \label{eq:est_on_etas}
  \eta^{(\ell)}_{DG}\lesssim \left(\frac{h}{p}\right)^{\ell}
  \qquad \text{and} \qquad \eta^{(\A)}_{DG}\lesssim  k^{1/2} \Big(\frac{h}{h+\sigma}\Big)^p + k^{3/2}\Big(\frac{ k h }{\sigma p}\Big)^p.
  \end{align}
\emph{Step 4 (quasi-optimality):}  We assert that, for sufficiently small $c_1$, assumption (\ref{eq:scale-resolution-condition-hp}) implies the resolution condition (\ref{eq:scale-resolution-condition-abstract}). That is, we assert that 
$$
\frac{kh}{p} + k^{1 + \mu_{stab} + \mu_{\A}} \eta^{(\A)}_{DG} 
$$
is sufficiently small under the constraint $kh/p \leq c_1$. This is shown in \cite[Lem.~{9.5}]{MS23} by noting that the core of the proof of 
\cite[Lem.~{9.5}]{MS23} is the assertion that the left-hand side of 
\cite[(9.35)]{MS23} can be made small by tuning $c_1$. 

\emph{Step 5 (error estimates for $f \in  H_{\PP}^{\ell-1}(\Omega)$)}:
\cite[Thm.{3.5}]{bernkopf} (see also \cite[Sec.~{5.1}]{bernkopf}, which asserts the validity of the assumptions 
in \cite[Thm.~{3.5}]{bernkopf}) shows that the exact solution $u$ to \eqref{eq:model_problem}
  can be split as
  \begin{align*}
    u = u_{\F} + u_{\A} \qquad
    \text{ with } \quad\|u_{\F}\|_{{1+\ell},\PP} \lesssim \|f\|_{-1+\ell,\PP}
    \quad \text{ and } \quad
    u_{\A} \in \aclass \big( C k^{\beta_0} \|f\|_{0,\Omega}, \vartheta, \PP \big).
  \end{align*}
  By using trace estimates, we easily get that $m = \gammaoint u + \ii k u$ and
  $\uext=\gammazint u$ can also be split~as
  \begin{align*}
    m &= m_{\F} + m_{\A} \quad
    \text{with} \quad \|m_{\F}\|_{{-1/2+\ell},\Gamma} \lesssim \|f\|_{{-1+\ell},\PP}
        \; \text{ and } \;
    m_{\A} \in \aclass \big( Ck^{\beta_0+1} \|f\|_{0,\Omega}, \vartheta, \mathcal{O},\Gamma\big), \\
    \uext &= \uext_{\F} + \uext_{\A} \quad
            \text{with} \quad \|\uext_{\F}\|_{{1/2+\ell},\Gamma} \lesssim \|f\|_{{-1+\ell},\PP}
            \;\text{ and }\;
    \uext_{\A} \in \aclass \big( C k^{\beta_0}\|f\|_{0,\Omega}, \vartheta, \mathcal{O}, \Gamma \big),
  \end{align*}
for some tubular neighborhood $\mathcal{O}$ of $\Gamma$. 
  Since we have already established the best approximation estimates
  on $\eta^{(\ell)}_{DG}$ and $\eta^{(\A)}_{DG}$ in (\ref{eq:est_on_etas}), the statement follows.
\end{proof}
%
%
%

\begin{remark}
The scale resolution condition (\ref{eq:scale-resolution-condition-hp}) hinges on the fact that the meshes $\Omega_h$ are 
regular, i.e., do not have hanging nodes, so that the conforming $hp$-FEM subspace is sufficiently rich. (A reflection of
this is that Lemma~\ref{lemma:broken-H2-approximation} constructs an approximant from the $H^1$-conforming $\Vh^{C}$ 
instead of $\Vh^{DG}$.)
The condition
``$kh/p$ sufficiently small'' has to be replaced with ``$kh/\sqrt{p}$ sufficiently small'' for more general meshes,
as discussed in \cite{MPS13}. 
\eremk
\end{remark}

\subsection{Conforming FEM}
\label{sec:conforming-fem-convergence}
For completeness, we also include the convergence analysis of the
conforming FEM.

Similarly to the DG case but with different norms, for $m \in \N$, we define the approximation quantities
\begin{align}  
  \eta^{(s)}_C:= \adjustlimits\sup_{0 \ne \x \in \mathcal{V}^{s}}\inf_{\xh \in \mathcal{V}_h^C}
  \frac{\smalltriplenorm{}{\x-\xh}_{k}}{\smalltriplenorm{}{\x}_{k,\mathcal{V},s}},
  \qquad \eta^{(\A)}_{C}:=\adjustlimits\sup_{\x \in \mathcal{A}^{1}}
  \inf_{\xh \in \mathcal{V}_h^C}
  \smalltriplenorm{}{\x - \xh}_{k}
\end{align}
and assume them to be sufficiently small.

We start with a quasi-optimality result, assuming weak resolution conditions.
\begin{theorem}  
  \label{thm:fem_is_quasioptimal}
  Let $(u,m,\uext) \in \mathcal{V}$ be the solution to~\eqref{system:Helmholtz}
  and let $(\uh,\mh, \uhext)\in \mathcal{V}_h^C$ be the solution of method~\eqref{equivalent:weak:formulation}.
  Let $\mu_{\A}:=\beta+3$ with $\beta$ as in Lemma~\ref{lemma:existence_primal}
  and let $\mu_{stab}  \in [0,4]$ be as in  Corollary~\ref{cor:TDG_is_bounded_with_k}.
  There is a positive $q <1$ 
  independent of $k\ge k_0$
  such that, 
under the resolution condition 
 \begin{equation} 
  \label{eq:thm:fem_is_quasioptimal-10}
  k \eta_C^{(2)} + k^{1+\mu_{stab} + \mu_{\A}} \eta_C^{(\A)} \leq q, 
 \end{equation} 
  the following estimate holds true: 
  \begin{multline*}
        \Vert \nu^{1/2} \nabla(u-\uh) \Vert_{0,\Omega}
        +\Vert \kn (u-\uh) \Vert_{0,\Omega}
 + \Vert m-\mh\Vert_{-{1/2},\Gamma} +\Vert \uext-\uhext\Vert_{{1/2}, \Gamma}\\
\lesssim \Vert u-\vh \Vert_{0,k,\Omega}+\Vert m-\lambdah\Vert_{-{1/2} ,\Gamma} +\Vert \uext-\vhext\Vert_{{1/2}, \Gamma}
\end{multline*}
for any $(\vh,\lambdah,\vhext) \in \mathcal{V}_h^C$,
with hidden constant independent of $k$, $h$, and $p$.
\end{theorem}
\begin{proof}
  The proof of Theorem~\ref{thm:dg_is_quasioptimal} could be
    repeated almost \emph{verbatim},
  as all the main building blocks (G{\aa}rding inequality, boundedness of the
  sesquilinear form) are valid also
  in this case  (Lemma~\ref{lemma:introducing_theta}, Proposition~\ref{prop:conforming}), and the adjoint problem is the same.
  The result analogous to Corollary~\ref{cor:TDG_is_bounded_with_k} follows
  from Proposition~\ref{prop:conforming} with the mapping properties
  of $\Theta$ found in Lemma~\ref{lemma:introducing_theta}. 

  Here, one could even simplify the argument in two ways:
  \begin{enumerate}[label=(\roman*)]
    \item Both the  $\energynorm{\cdot}$ and $\energynormp{\cdot}$ norms can be replaced
      by the standard norm $\|\cdot\|_{1,k,\PP}$. 
    \item The  G{\aa}rding inequality holds true also for non-discrete functions.
      Thus, one could start directly with $\x-\xh$ instead of
      $\xh -\yh$ in~\eqref{eq:energy_estimate1} and avoid some terms.
    \end{enumerate}
    None of these differences has a significant impact on the overall argument and
    thus we omit the full proof for brevity.
  \end{proof}

\begin{corollary}
  \label{cor:conforming_explicit_convergence}
  Let $f \in L^2(\Omega)$ with $f|_{P_j} \in H^{\ell-1}(P_j)$, $P_j \in \PP$, $j=0,\ldots,L$. 
  Let $(u,m,\uext) \in \mathcal{V}$ be the solution to~\eqref{system:Helmholtz},
  and let $(\uh,\mh, \uhext)\in \mathcal{V}_h^C$ be the solution of method~\eqref{equivalent:weak:formulation}.
Given $c_2 > 0$, there is $c_1>0$ independent of $h$, $p$, $k$ such that, under the scale resolution condition 
  \begin{align*}
     \frac{ k h}{p} &\leq c_1,  \qquad \text{and} \qquad p \geq \max (1,c_2 \ln k),
  \end{align*}
the Galerkin error is quasi-optimal with respect to $\enorm{\cdot}$ (see (\ref{eq:enorm})), with constants independent of $k$. 
Furthermore, there exists a constant $\sigma >0$ independent of $h$, $p$, $k$ such that 
  \begin{multline*}
    \Vert \nu^{1/2} \nabla(u-\uh) \Vert_{0,\Omega}
    +\Vert \kn (u-\uh) \Vert_{0,\Omega}
    + \Vert m-\mh\Vert_{-{1/2},\Gamma} +\Vert \uext-\uhext\Vert_{{1/2}, \Gamma} \\
    \lesssim \Bigg( \Big(\frac{h}{p}\Big)^{\ell}
 + k^{\beta_0+1}\Big(\Bigl(\frac{h}{h+\sigma}\Bigr)^p + k\Big(\frac{ k h }{\sigma p}\Big)^p\Bigg),
\end{multline*}
with hidden constant independent of $k$, $h$, and $p$. 
Here $\beta_0$ is the constant in Assumption~\ref{ass:existence_global}.
\end{corollary}
\begin{proof}
  The result follows from the best approximation result in Theorem~\ref{thm:fem_is_quasioptimal}
  and standard approximation results for FEM spaces, see, e.g., \cite{MS11}.
  The argument is very similar to Corollary~\ref{cor:dg_explicit_convergence}
  and is thus omitted here.
\end{proof}

\appendix
\section{Approximation in DG-norms}
We generalize the approximation results of \cite[Thm.~{4.11}]{MPS13} to piecewise smooth functions. 
\subsection{Norm equivalences and liftings for weighted spaces} 

We start with the following preparatory result concerning
interpolation spaces with weighted norms. 
In the following, for any two continuously embedded Banach spaces $X_1\subset X_0$ and $0 < \theta < 1$, 
the space $[X_0,X_1]_{\theta,2}$ is defined with the $K$-method of real interpolation as introduced in, e.g., 
\cite[Sec.~{22}]{tartar07}. 
\begin{lemma}
\label{lemma:interpolation-in-weighted-spaces-abstract}
    Let $X_1 \subseteq X_0$ be Banach spaces with continuous inclusion.  
    Let $|\cdot|_{X_1}$ be a (semi)-norm on $X_1$. Introduce the interpolation (semi)-norm
    $$
    |u|^2_{X_\theta}:=\int_{0}^{\infty}\left(t^{-\theta}
    \inf_{v\in X_1}{\|u-v\|_{X_0} + t |v|_{X_1}}\right)^2\,\frac{dt}{t}.
    $$
    For $\tau>0$, define the following weighted norms:
    \begin{align*}
    \|u\|_{X_1,\tau}&:=|u|_{X_1} + \tau^{-1} \|u\|_{X_0}, & 
    \|u\|_{X_0,\tau}&:=\|u\|_{X_0}, \\
       \|u\|_{X_\theta,\tau}& :=|u|_{X_\theta} + \tau^{-\theta} \|u\|_{X_0}, \quad 0 < \theta < 1. 
    \end{align*}
    Then for $0 < \theta  <1$  these weighted norms can be characterized as 
    interpolation norms:
    For $0\leq \theta_0 < \theta < \theta_1 \leq 1$ there holds
    \begin{align*}
    \big\|u\big\|_{X_\theta,\tau} \sim \big\|u\big\|_{
    [(X_{\theta_0},\|\cdot\|_{X_{\theta_0},\tau})
,(X_{\theta_1},\|\cdot\|_{X_{\theta_1},\tau})]_{\eta,2}}
\qquad\text{with} \qquad \eta:=\frac{\theta-\theta_0}{\theta_1-\theta_0}
    \end{align*}
with implied constants depending only on $\theta_0$, $\theta_1$, and $\theta$ but independent of $\tau$.
\end{lemma}
\begin{proof}
\emph{Step 1:}
By \cite[Lemma~{4.1}]{karkulik-melenk-rieder20}, the equivalence
 $ \|u\|_{X_{\theta,\tau}} \sim 
 \|u\|_{[(X_{0},\|\cdot\|_{X_0,\tau}),(X_{1},\|\cdot\|_{X_1,\tau})]_{\theta,2}} $
holds with implied constants depending only on $\theta$. 

\emph{Step 2:} In the notation of \cite[Thm.~{26.3}]{tartar07} the Reiteration Theorem takes the following form: 
If for $E_1 \subset E_0$ and spaces $F_0$, $F_1$ one has the continuous embeddings, 
$[E_0,E_1]_{\theta_i,1,J} \subset F_i \subset [E_0,E_1]_{\theta_i,\infty,K}$, $i=0$, $1$, 
then $[E_0,E_1]_{\eta,2} = [F_0,F_1]_{\theta,2}$ with equivalent norms. Inspection of the 
proof of \cite[Thm.~{26.3}]{tartar07} shows that the equivalence constants depend only
on $\theta_0$, $\theta_1$, $\theta$, and the embedding constants 
of $[E_0,E_1]_{\theta_i,1,J} \subset F_i \subset [E_0,E_1]_{\theta_i,\infty,K}$. 
We specify now $F_i$: If $\theta_i \in (0,1)$, then we select $F_i:= [E_0,E_1]_{\theta_i,2}$. 
By \cite[Lemmata~{22.2}, {24.3}]{tartar07}, we have the continuous embeddings 
$[E_0,E_1]_{\theta_i,1,J} \subset F_i \subset [E_0,E_1]_{\theta_i,\infty,K}$ with embedding
constants depending solely on $\theta$. In the limiting cases $\theta_0 = 0$ (for $i = 0$) 
and $\theta_1 = 1$ (for $i=1$) we set $F_0:= E_0$ and $F_1 := E_1$. By definition, 
it follows $\|u\|_{[E_0,E_1]_{0,\infty,K}} \leq \|u\|_{E_0}$ for all $u \in E_0$ and 
$\|u\|_{[E_0,E_1]_{1,\infty,K}} \leq \|u\|_{E_1}$ for all $u \in E_1$. For the remaining 
embeddings,  inspection of the proof of \cite[Lemma~{25.2}]{tartar07} shows that 
the embedding constants $C_1$ and $C_2$ in the estimates $\|u\|_{E_0} \leq C_1 \|u\|_{[E_0,E_1]_{0,1,J}}$ and 
$\|u\|_{E_1} \leq C_2 \|u\|_{[E_0,E_1]_{1,1,J}}$ are the best constants in the multiplicative
estimates $\|u\|_{E_0} \leq C_1 \|u\|^{1-0}_{E_0} \|u\|^0_{E_1}$ and 
$\|u\|_{E_1} \leq C_2 \|u\|^{1-1}_{E_0} \|u\|^1_{E_1}$ for all $u \in E_1$. 
Clearly, they are $C_1 = C_2 = 1$.  

\emph{Step 3:} 
We now apply the Reiteration Theorem \cite[Thm.~26.3]{tartar07} as discussed in Step~2 
with $E_0 = (X_{0},\|\cdot\|_{X_0,\tau})$, $E_1 = (X_{1},\|\cdot\|_{X_1,\tau})$ 
and $F_i$ as in Step~2. We obtain
\begin{align}
\label{eq:lemma:interpolation-in-weighted-spaces-abstract-100}
\big\|u\big\|_{[(X_{0},\|\cdot\|_{X_{0},\tau}),(X_{1}, \|\cdot\|_{X_{1},\tau})]_{\theta,2}}
\sim  \big\|u\big\|_{[(X_{\theta_0},\|\cdot\|_{X_{\theta_0},\tau})
,(X_{\theta_1},\|\cdot\|_{X_{\theta_1},\tau})]_{\eta,2}}
\end{align}
with implied constants depending only on $\theta_0, \theta_1$ and $\theta$, but independent of $\tau$. 
The proof of the lemma is completed by noting that the left-hand side of 
\eqref{eq:lemma:interpolation-in-weighted-spaces-abstract-100} is equivalent to $\|u\|_{X_{\theta},\tau}$ by Step~1. 
\end{proof}
\begin{corollary}
\label{cor:interpolation-in-weighted-sobolevspaces}
Let $\mathcal{O}$, $\widetilde {\mathcal O} \subset {\mathbb R}^d$ be bounded Lipschitz domains and 
$\partial\widetilde{\mathcal{O}}$ be smooth. Let $k \ge k_0 > 0$, 
$0 \leq s_0 < s_1$. Then for every 
$0<\theta< 1$ there hold the norm equivalences 
\begin{align*}
\|u\|_{(1-\theta)s_0 + \theta s_1,k,\mathcal{O}} &\sim \|u\|_{[(H^{s_0}(\mathcal{O}),\|\cdot\|_{{s_0},k,\mathcal{O}}),(H^{s_1}(\mathcal{O}),\|\cdot\|_{s_1,k,\mathcal{O}})]_{\theta,2}}, \\
\|u\|_{(1-\theta) s_0 + \theta s_1,k,\partial\widetilde{\mathcal{O}}} &\sim \|u\|_{[(H^{s_0}(\partial{\widetilde{\mathcal{O}}}),\|\cdot\|_{s_0,k,\partial{\widetilde{\mathcal{O}}}}),(H^{s_1}(\partial{\widetilde{\mathcal{O}}}),\|\cdot\|_{s_1,k,\partial{\widetilde{\mathcal{O}}}})]_{\theta,2}} 
\end{align*}
with implied constants independent of $k$. 
\end{corollary}
\begin{proof}
Follows from Lemma~\ref{lemma:interpolation-in-weighted-spaces-abstract}. 
An alternative proof can be obtained by eigenfunction expansions. 
For example, for the first estimate, consider the variational eigenvalue problem: Find $(u,\lambda) \in H^{s_1}(\mathcal{O})\setminus \{0\} \times {\mathbb R}$
such that $(u,v)_{H^{s_1}(\mathcal{O})} = \lambda^{s_1} (u,v)_{L^2(\mathcal{O})}$ for all $v \in H^{s_1}(\mathcal{O})$ where $(\cdot,\cdot)_{H^{s_1}(\mathcal{O})}$ denotes the scalar product on 
$H^{s_1}(\mathcal{O})$. One can normalize the eigenpairs $(\varphi_j,\lambda_j)_{j \in {\mathbb N}_0}$ 
such they form an orthonormal basis of $L^2(\mathcal{O})$ and an orthogonal basis of $H^{s_1}(\mathcal{O})$. This, together with interpolation, results in the norm equivalences 
$\|u\|^2_{H^{\theta s_1}(\mathcal O)} \sim \sum_{j\in \N_0} (1+\lambda_j^{\theta s_1}) u_j^2$ for any $\theta \in [0,1]$, where $u_j = (u,\varphi_j)_{L^2({\mathcal O})}$. 
Hence, 
$\|u\|^2_{s',k,\mathcal{O}} \sim \sum_{j \in \N_0} \left(\lambda_j^{s'} + k^{2s'}\right) u_j^2$ for $0 \leq s' \leq s_1$. 
From \cite[Sec.~{23}]{tartar07}, one gets, with constants depending only on $\theta$, 
\begin{align*}
\|u\|^2_{[(H^{s_0}(\mathcal{O}),\|\cdot\|_{s_0,k,\mathcal{O}}),(H^{s_1}(\mathcal{O}),\|\cdot\|_{s_1,k,\mathcal{O}})]_{\theta,2}} 
& \sim \sum_{j\in \N_0} \left(\lambda_j^{s_0} + k^{2s_0}\right)^{1-\theta}\left(\lambda_j^{s_1} + k^{2s_1}\right)^{\theta}  u_j^2 \\
& \sim \sum_{j\in \N_0} (\lambda_j^{(1-\theta) s_0 + \theta s_1} + k^{2(1-\theta)s_0 +\theta s_1})  u_j^2 
\sim \|u\|^2_{(1-\theta) s_0 + \theta s_1,k,\mathcal{O}}, 
\end{align*}
where the second equivalence follows from distinguishing between the cases $\lambda_j \leq k^2$ and $\lambda_j > k^2$. 
\end{proof}

We are in position to prove a lifting result. 

\begin{lemma}
\label{lemma:lifting} Let $\omega \subset {\mathbb R}^d$, $d \in {\mathbb N}$, be a bounded domain with a smooth boundary. 
For $\tau  \in (0,1]$ and $g \in H^{s-1/2}(\partial\omega)$ with $s\geq 1$, let ${\mathcal G}$ be the solution to
\begin{align*}
-\Delta {\mathcal G} + \tau^{-2} {\mathcal G} = 0 \quad \mbox{ in $\omega$}, \qquad {\mathcal G}{}_{|_{\partial\omega}} = g. 
\end{align*}
Then, 
there is $C_s > 0$ independent of $\tau$ and $g$ such that 
\begin{align*}
\|{\mathcal G}\|_{s,\omega} + \tau^{-s} \|{\mathcal G}\|_{0,\omega} \leq C_s \left[ \|g\|_{{s-1/2},\partial\omega} + \tau^{-(s-1/2)} \|g\|_{0,\partial\omega}\right]. 
\end{align*}
\end{lemma}
\begin{proof} 
For $\sigma > 0$, we introduce the norms 
$\|{\mathcal G}\|_{\sigma,\tau,\omega}:= \|{\mathcal G}\|_{\sigma,\omega} + \tau^{-\sigma} \|{\mathcal G}\|_{0,\omega}$,
and 
$\|g\|_{\sigma,\tau,\partial\omega}:= \|g\|_{\sigma,\partial\omega} + \tau^{-\sigma} \|g\|_{0,\omega}$. 

\emph{Step 0 (a norm equivalence):}
Before we tackle the PDE, we note that we can estimate lower order norms by a weighted linear combination
of a higher order norm and the $L^2$-norm. Namely, for $\sigma\geq 2$, the following estimate holds:
\begin{align}
\label{eq:norm_eq_interpol}
 \|g\|_{\sigma-2,\partial \omega} \lesssim \tau^{2} \big(\|g\|_{\sigma,\partial \omega} + \tau^{-\sigma} \|g\|_{0,\partial \omega}\big).
\end{align}
This follows by a multiplicative interpolation estimate and Young's inequality.
Namely, for $\theta:={(\sigma-2)}/{\sigma}$, we have
\begin{align*}
\|g\|_{\sigma-2,\partial \omega}
&\lesssim \tau^{2}  \|g\|_{\sigma,\partial \omega}^{\theta} \tau^{-2} \|g\|_{0,\partial \omega}^{1-\theta} 
\lesssim \tau^2\big(\|g\|_{\sigma,\partial \omega} + \tau^{-2/(1-\theta)} \|g\|_{0,\partial \omega}\big).
\end{align*}
By inserting the value of $\theta$, estimate~\eqref{eq:norm_eq_interpol} follows.

\emph{Step 1 (shift theorem for $s \in {\mathbb N}$):} The case $s = 1$ is given in \cite[Lemma~{4.5}]{melenk-rieder21}. For $s \in {\mathbb N}$, the estimate 
$\|{\mathcal G}\|_{s,\tau,\omega} \leq C_s \|g\|_{s-1/2,\tau,\partial\omega}$ follows inductively from
the standard shift theorem for the Laplacian see, e.g., \cite[Thm. 8.13]{gilbarg_trudinger} by writing $-\Delta {\mathcal G} = -\tau^{-2} {\mathcal G}$ on $\omega$ with Dirichlet conditions ${\mathcal G} = g$ on $\partial\omega$: 
\begin{align*}
\|\mathcal{G}\|_{s,\omega} + \tau^{-s}\|\mathcal{G}\|_{0,\omega}
&\lesssim \tau^{-2}\|\mathcal{G}\|_{s-2,\omega} + \|g\|_{s-1/2,\partial \omega} + \tau^{-s}\|\mathcal{G}\|_{0,\omega} \\
&\lesssim \tau^{-2}\left(\|\mathcal{G}\|_{s-2,\omega}+\tau^{-(s-2)}\|\mathcal{G}\|_{0,\omega}\right)+ \|g\|_{s-1/2,\partial \omega}\\
&\lesssim \tau^{-2} \Big(\|g\|_{s-2-1/2,\partial \omega} + \tau^{-(s-2-1/2)} \|g\|_{0,\partial \omega}\Big)
+ \|g\|_{s-1/2,\partial \omega} \\
&\lesssim \|g\|_{s-1/2,\partial \omega} + \tau^{-(s-1/2)} \|g\|_{0,\partial_\omega}, 
\end{align*}
where in the last step we used the norm-equivalence in \eqref{eq:norm_eq_interpol}.

%

\emph{Step 2 (shift theorem for $s \not\in {\mathbb N}$):} The desired shift theorem for $s \not\in {\mathbb N}$ follows by interpolation from the case of 
$s \in {\mathbb N}$ and the norm equivalences for interpolation spaces given in Corollary~\ref{cor:interpolation-in-weighted-sobolevspaces}. 
\end{proof}
\subsection{Approximation in DG-norms}\label{app:apprDG}
\begin{lemma}
\label{lemma:broken-H2-approximation}
Let Assumption~\ref{def:element-maps} be valid. Then, for every $v \in H^{s+2}(\PP)$, $s \ge 0$ and $p \ge s+1$, there is $v_h \in \Vh^C \subset H^1(\Omega)$ 
with 
\begin{equation*}
\|v - v_h\|_{\dGnormp} \lesssim \left(\frac{h}{p}\right)^{s+1}\left(1 + \frac{kh}{p}\right) \|v\|_{s+2,\PP}. 
\end{equation*}
\end{lemma}
\begin{proof}
In~\cite[Lemma~{4.7}]{MPS13}, it is proven that, for element maps satisfying Assumption~\ref{def:element-maps}, usual scaling arguments lead to the expected powers of the mesh size. The approximation $v_h$ is constructed elementwise with
an operator defined on the reference element in such a way that an~$H^1$-conforming interpolant is obtained for $v \in H^{t+2}_{pw}(\taun) \cap H^1(\Omega)$;
see~\cite[Thm.~{7.3}, Cor.~{7.4}]{MPS13} 

As in the proof of \cite[Thm.~{4.11}]{MPS13}, we use the operator $I:H^t_{pw}(\taun)\cap H^1(\Omega) \rightarrow S^{p,1}(\taun)$
of \cite[Thm.~{7.3}, Cor.~{7.4}]{MPS13}
with the approximation property 
\begin{align}
\label{eq:lemma:broken-H2-approximation-5}
\left(\frac{h_K}{p}\right)^2 \|\nabla^2 (v - Iv)\|_{L^2(K)} + \left(\frac{h_K}{p}\right) \|\nabla (v - Iv)\|_{L^2(K)} + \|v - Iv\|_{L^2(K)} 
\lesssim \left(\frac{h_K}{p}\right)^t \|v\|_{H^t(K)}, 
\end{align}
which is valid under the \textsl{proviso} $t > t^*_d$, where 
$t^*_d := 5/2$ for $d = 2$ and $t^*_d := 5$ for $d = 3$ as well as the assumption $p \ge t-1$. 

If $s+2 > t^*_d$, then we take the approximation $v_h:= I v$ and the result follows from elementwise multiplicative trace inequalitites
as shown in more detail in~\eqref{eq:lemma:broken-H2-approximation-90} below. 

If $s + 2 \leq t^*_d$, fix $\sigma > t^*_d$
and use the characterization of Sobolev spaces as interpolation spaces with the ``real method'' (see, e.g., \cite{tartar07}) 
to write $H^{s+2}(P_i) = (L^2(P_i), H^{\sigma}(P_i))_{\theta,2}$ with $\theta = (s+2)/\sigma$ for each $i=0,\ldots,L$. 
As in the proof of \cite[Thm.~{4.11}]{MPS13}, this characterization allows us to find, for each $\tau > 0$, a function $v_{i,\tau} \in H^\sigma(P_i)$ 
with 
\begin{align}
\label{eq:lemma:broken-H2-approximation-10}
\tau^{\sigma-(s+2)} \|v_{i,\tau}\|_{\sigma,P_i} 
&\lesssim \|v\|_{s+2,P_i}, 
& 
\|v - v_{i,\tau}\|_{\rho,P_i} & \lesssim \tau^{s+2-\rho}\|v\|_{s+2,P_i}, \quad 0 \leq \rho \leq s+2. 
\end{align}
The piecewise defined function $\widetilde v_\tau$ given by $\widetilde v_\tau|_{P_i} = v_{i,\tau}$ is piecewise smooth but not necessarily in $H^1(\Omega)$. 
This is corrected with a lifting. We note that the multiplicative trace inequality and the trace inequality yield for the 
jump $\kappa_{ij}:= v_{i,\tau} - v_{j,\tau}$ on $\partial P_i \cap \partial P_j$, in view of 
(\ref{eq:lemma:broken-H2-approximation-10}) and the fact that $v \in H^1(\Omega)$, 
\begin{align}
\label{eq:lemma:broken-H2-approximation-20}
\|\kappa_{ij}\|_{0,\partial P_i \cap \partial P_j} &\lesssim \tau^{s+2-1/2} \|v\|_{s+2,\PP}, 
&
\|\kappa_{ij}\|_{\sigma-1/2,\partial P_i \cap \partial P_j} &\lesssim \tau^{s+2-\sigma} \|v\|_{s+2,\PP}, 
\end{align}
which implies 
\begin{align}
\label{eq:lemma:broken-H2-approximation-30}
\tau^{-(\sigma-1/2)} \|\kappa_{ij}\|_{0,\partial P_i \cap \partial P_j} + \|\kappa_{ij}\|_{\sigma-1/2,\partial P_i \cap \partial P_j} 
&\lesssim \tau^{s+2-\sigma} \|v\|_{s+2,\PP}. 
\end{align}
Note that $\partial P_i \cap \partial P_j$ is a smooth manifold without boundary by our assumptions on $\PP$. 
The lifting of Lemma~\ref{lemma:lifting} allows us to correct the jump across $\partial P_i \cap \partial P_j$ with a 
function $J_{ij}$ (supported by either $P_i$ or $P_j$) with 
$$
\|J_{ij} \|_{r,\PP} \lesssim \tau^{s+2-r} \|v\|_{s+2,\PP}, \qquad 0 \leq r \leq s+2. 
$$
Using liftings for all interfaces leads to a function $v_\tau:= \widetilde v_\tau + \sum_{ij} J_{ij} \in H^{\sigma}(\PP) \cap H^1(\Omega)$ with 
\begin{align}
\label{lemma:broken-H2-approximation-40}
\|v_\tau\|_{\sigma,\PP} & \lesssim \tau^{s+2-\sigma} \|v\|_{s+2,\PP}, 
&
\|v - v_\tau\|_{r,\PP} & \lesssim \tau^{s+2-r} \|v\|_{s+2,\PP}, \quad 0 \leq r \leq s+2. 
\end{align}
As in the proof of \cite[Thm.~{4.11}]{MPS13}, we select $\tau = h/p$. Elementwise appropriate multiplicative trace inequalities 
yield 
\begin{equation*}
\|v - v_\tau\|_{\dGnormp} \lesssim \left( \left(\frac{h}{p} \right)^{s+1} + k \left(\frac{h}{p} \right)^{s+2} +\sqrt{k} \left(\frac{h}{p}\right)^{s+3/2}\right)\|v\|_{s+2,\PP},
\end{equation*}
and~\eqref{eq:lemma:broken-H2-approximation-5} with $t = \sigma$ results in 
\begin{align}
\label{eq:lemma:broken-H2-approximation-90}
\|v_\tau - I v_\tau\|_{\dGnormp} &\lesssim \left( 
\left(\frac{h}{p}\right)^{\sigma-1} + 
k \left(\frac{h}{p}\right)^{\sigma} + \sqrt{k} \left(\frac{h}{p}\right)^{\sigma-1/2}\right) \|v_\tau\|_{\sigma,\PP} \\
\nonumber
&
\stackrel{ (\ref{lemma:broken-H2-approximation-40})} {\lesssim }
\left(\frac{h}{p}\right)^{s+1} \left( 1 + \frac{kh}{p} + \left(\frac{kh}{p}\right)^{1/2}\right) \|v\|_{s+2,\PP}. 
\end{align}
This shows the statement under the assumption that $p \ge \sigma-1$. The remaining case $s+1 \leq p < \sigma -1$ is a pure $h$-version
statement, which follows by usual scaling arguments and a polynomial reproducing interpolation operator (note: $s+2 \ge 2$ and $d \leq 3$ so that 
Lagrange interpolation is admissible).
\end{proof}

\bibliographystyle{alphaabbr}
\bibliography{literature}

\end{document}